\documentclass[reqno, 11pt]{amsart}

\usepackage[all]{xy}
\usepackage{url}
\usepackage{amssymb}
\usepackage{hyperref}
\usepackage[margin=1 in,marginpar=0.75in]{geometry}
\usepackage{color,soul}
\usepackage{amsthm}
\usepackage{enumitem}
\usepackage{comment}
\usepackage{varwidth}
\usepackage[textwidth=0.75in]{todonotes}

\input{kmacros3.sty}
\usepackage{mabliautoref}

\usepackage{mathtools}

\DeclareMathOperator{\td}{tr.deg}
\DeclareMathOperator{\pos}{pos}
\DeclareMathOperator{\Int}{Relint}
\DeclareMathOperator{\Aff}{Aff}
\DeclareMathOperator{\Lin}{Lin}

\numberwithin{equation}{subsubsection}
\numberwithin{equation}{subsubsection}
\numberwithin{figure}{section}
\normalfont
\usepackage[T1]{fontenc}{}
\usepackage{setspace}

\usepackage{tikz-cd}

\makeatletter
\def\@tocline#1#2#3#4#5#6#7{\relax
  \ifnum #1>\c@tocdepth 
  \else
    \par \addpenalty\@secpenalty\addvspace{#2}%
    \begingroup \hyphenpenalty\@M
    \@ifempty{#4}{%
      \@tempdima\csname r@tocindent\number#1\endcsname\relax
    }{%
      \@tempdima#4\relax
    }%
    \parindent\z@ \leftskip#3\relax \advance\leftskip\@tempdima\relax
    \rightskip\@pnumwidth plus4em \parfillskip-\@pnumwidth
    #5\leavevmode\hskip-\@tempdima
      \ifcase #1
       \or\or \hskip 1em \or \hskip 2em \else \hskip 3em \fi%
      #6\nobreak\relax
    \dotfill\hbox to\@pnumwidth{\@tocpagenum{#7}}\par
    \nobreak
    \endgroup
  \fi}
\makeatother
\usepackage{hyperref}

\title{\vspace{-5em}Finite generation of split-$F$-regular monoid algebras}
\author{Rankeya Datta}
\address{Department of Mathematics, University of Missouri, Columbia, MO, USA}
\email{rankeya.datta@missouri.edu}
\thanks{Datta was supported by a grant from the Simons Foundation \# MP-TSM-00002400}
\author{Karl Schwede}
\address{Department of Mathematics, University of Utah, Salt Lake City, UT, USA}
\email{schwede@math.utah.edu}
\thanks{Schwede was supported by NSF Grants \#2101800, \#1952522, \#1801849 and a Simons Foundation Fellowship.}
\author{Kevin Tucker}
\address{Department of Mathematics, University of Illinois at Chicago, Chicago, IL, USA}
\email{kftucker@uic.edu}
\thanks{Tucker was supported by NSF Grant DMS \#2200716.}

\begin{document}
\begin{abstract}
    Let $S$ be a submonoid of a free Abelian group of finite rank. We show that if $k$ is a field of prime characteristic such that the monoid $k$-algebra $k[S]$ is split-$F$-regular, then $k[S]$ is a finitely generated $k$-algebra, or equivalently, that $S$ is a finitely generated monoid. Split-$F$-regular rings are possibly non-Noetherian or non-$F$-finite rings that satisfy the defining property of strongly $F$-regular rings from the theories of tight closure and $F$-singularities. Our finite generation result provides evidence in favor of the conjecture that split-$F$-regular rings in function fields over $k$ have to be Noetherian. The key tool is Diophantine approximation from convex geometry.
\end{abstract}

 \vspace*{40pt}
\maketitle
\vspace*{-27pt}
\renewcommand{\baselinestretch}{0.75}\normalsize
{\small
\tableofcontents}
\renewcommand{\baselinestretch}{1.0}\normalsize
\newpage
\section{Introduction}
  
 \emph{Strongly $F$-regular} rings form an important class of singularities in prime characteristic $p>0$ commutative algebra arising from Hochster and Huneke's celebrated theory of \emph{tight closure} \cite{HochsterHunekeTightClosureAndStrongFRegularity}.
 Their prominence stems partly from their connection to KLT singularities in characteristic zero \cite{HaraWatanabeFRegFPure}. While originally defined in the Noetherian $F$-finite setting, natural extensions to Noetherian excellent rings have since been introduced and studied \cite{SmithThesis,HochsterFoundations, HashimotoF-pure-homomorphisms,DattaSmithFrobeniusAndValuationRings,DattaMurayamaF-Solid, HochsterYaoSFR, HochsterYaoGenericLocal, DattaEpsteinTucker} (see \cite{HochsterTightClosureandStronglyFregularRings2022} for a recent overview). In this article, we look to impose similar conditions on a priori non-Noetherian rings -- with a view towards ultimately deducing Noetherianity in the setting we work in. 
 
Suppose $R$ is a possibly non-Noetherian integral domain of prime characteristic $p > 0$.  We say that $R$ is \emph{split-$F$-regular} if for every $0 \neq c \in R$ there exists an integer $e > 0$ and an $R$-linear map $\phi : R^{1/p^{e}} \to R$ with $\phi(c^{1/p^e})=1$.  When $R$ is Noetherian and $F$-finite, this condition is called \emph{strong $F$-regularity}. This latter notion was introduced by Hochster and Huneke as a localization-stable variant of the tight closure notion of weak $F$-regularity \cite{HochsterHunekeTightClosureAndStrongFRegularity}. While there are certainly split-$F$-regular rings that are not Noetherian, such as $\bF_p[x_1, x_2, \dots]$, the following conjecture has been discussed informally among experts for over a decade.
\begin{conj*}
    Suppose $K$ is a function field over a perfect field $k$ of characteristic $p > 0$.  If a $k$-algebra $R \subseteq K$ is split-$F$-regular, then $R$ is Noetherian.  
\end{conj*}

Note that, if one weakens split-$F$-regular to simply Frobenius split above, the analogous statement is false (see for instance \autoref{prop.SemiGroupAlgebrasAreFSplit} or \cite{DattaFrobeniusSplittingAbhyankar}). To date, the best direct evidence for this conjecture stems from work of the first author and K.~E.~Smith in \cite{DattaSmithFrobeniusAndValuationRings} who proved it when $R$ is a valuation ring. In this paper, we give additional evidence for the conjecture above by verifying it for certain monoid algebras -- namely those arising as submonoids of a lattice in analogy with the toric setting.

\begin{mainthm*}[\autoref{thm:split-F-regular-monoid-finite-generation}] Let $k$ be a field of characteristic $p > 0$ and $S$ be a submonoid of a free Abelian group $M$ of finite rank.  If the monoid algebra $k[S]$ is split-$F$-regular, then $k[S]$ is a finitely generated $k$-algebra. Equivalently, $S$ is a finitely generated monoid.
\end{mainthm*}

For geometric rings, \textit{i.e.} those of finite type over an $F$-finite field $k$, a central motivation in studying the above conjecture stems from the following two observations.

First, by \cite{HashimotoGFRgraded-splitting}, for a globally $F$-regular projective variety, we know that the Cox ring is split-$F$-regular.  
Hence, if the above conjecture holds then the Cox ring would be Noetherian and $X$ would be a Mori dream space under mild assumptions \cite{HuKeelMoriDreamSpaces}, \cf \cite[Theorem 3.3]{LeePandeFSignatureAmpleCone}.       

Second, and similarly, suppose $R$ is the symbolic Rees algebra of a pure height one (that is a divisorial) ideal in a strongly $F$-regular Noetherian local ring. It is not difficult to verify that $R$ is split-$F$-regular (\textit{e.g.}  \cite{WatanabeInfiniteCyclicCoversOfStronglyFRegular,DeStefaniMontanoNunezBetantcourtBlowupAlgebrasOfDeterminantalInPrimeChar}). Thus, if the above conjecture holds, then for a strongly $F$-regular ring essentially of finite type over $k$, any divisorial symbolic Rees algebra is Noetherian. In light of the exciting recent work of Aberbach-Huneke-Polstra, centrally important cases of the famous ``weak implies strong $F$-regularity conjecture'' would follow from finite generation of divisorial symbolic Rees algebras for strongly $F$-regular rings \cite{AberbachHunekePolstraEqualityOfTestIdeals}. 

By appealing to the minimal model program for threefolds in characteristic $p > 5$ \cite{HaconXuThreeDimensionalMinimalModel,Birkar16},  finite generation of divisorial symbolic Rees algebras is known to hold over strongly $F$-regular rings essentially of finite type over $k$ in dimension $\leq 3$ and characteristics $p > 5$, via \cite{HaraWatanabeFRegFPure,SchwedeSmithLogFanoVsGloballyFRegular}. Moreover, garnering intuition coming from reduction to characteristic $p>0$, such finite generation in higher dimensions is particularly plausible -- the divisorial symbolic Rees algebras of KLT singularities in characteristic zero are known to be finitely generated \cite[Exercises 108 and 109]{KollarExercisesInBiratGeom} as a consequence of the minimal model program in higher dimensions \cite{BirkarCasciniHaconMcKernan}. Note that related finite generation statements are in fact central to the proof of the minimal model program -- particularly the existence of flips; see also \cite{CasciniLazicMMPrevisited}.  

One way to construct non-Noetherian subrings of a function field is to take $\Gamma(U, \cO_X)$ for $U$ open in an affine variety $X = \Spec R$ (\cite{VakilNonFGNote,MurayamaNonNoetherianMathOverflow}).  If $R$ is split-$F$-regular, it is straightforward to see that so is $S = \Gamma(U, \cO_X)$, hence one could hope that this could give counter-examples to the conjecture.  However, in characteristic zero, if $R$ has KLT singularities then $\Gamma(U, \cO_X)$ is Noetherian, see for instance \cite[Lemma 2.6]{ZhuangDirectSummandsOfKLT}.  Note Zhuang's argument uses the aforementioned finite generation of the symbolic Rees algebras in characteristic zero.  




\subsection{Proof sketch}
We now turn to an overview of the proof of the main theorem. The first step is to translate the problem into the language of convex geometry. Thus, let $\mathbb{Z}S$ be the subgroup of $M$ generated by $S$ and let $\mathbb{R}S \coloneqq \mathbb{Z}S \otimes_{\mathbb Z} \mathbb{R}$. Let $\sigma_S$ be the convex cone in $\mathbb{R}S$ that is generated by $S$. By construction, $\sigma_S$ is full-dimensional in $\mathbb{R}S$. Since $k[S]$ is split-$F$-regular by assumption, one can show that $S$ is \emph{normal} in $\mathbb{Z}S$ \autoref{normal-semigroup}\ref{normal-semigroup.3}. Consequently, using Carath\'eodory's Theorem from convex geometry \autoref{Caratheodory}, one can deduce (see \autoref{integrally-closed-cones}) that
\[
\sigma_S \cap \mathbb{Z}S = S.	
\]
Thus, $k[S] = k[\sigma_S \cap \mathbb{Z}S]$, and the main theorem would follow from the following result about monoid algebras determined by convex cones in finite dimensional $\mathbb R$-vector spaces.

\begin{thm*}[\autoref{thm:SFR-cone-finite-generation}]
Suppose $M = \bZ^d$, and $\sigma \subseteq M \otimes_{\mathbb Z} \bR$ is a full-dimensional convex cone containing the origin.  Suppose that $k$ is a field of characteristic $p > 0$.  If $k[\sigma \cap M]$ is split-$F$-regular, then $\sigma$ is a closed cone generated by finitely many rational vectors.  In other words $\sigma$ is a rational polyhedral cone.  As a consequence, $k[\sigma \cap M]$ is a finitely generated $k$-algebra.
\end{thm*}

We use the following key observation throughout.
\begin{enumerate}
    \item If one can find a map $k[\sigma \cap M]^{1/p^e} \to k[\sigma \cap M]$ sending 
    $ 
        X^{\alpha/p^e} \mapsto 1,
    $ 
    for a lattice point $\alpha$ in the relative interior of $\sigma$, then $k[\sigma \cap M]$ is split-$F$-regular (the converse is clear for $e \gg 0$).  See \autoref{thm:splitting-rel-int-implies-SFR}. \label{itm.IntoA}
\end{enumerate}
We will call the condition in \autoref{itm.IntoA} \emph{splitting an element of the relative interior}.

We can also often restrict ourselves to the case of a strongly convex (aka pointed) cone $\sigma$ containing the origin, although in our general proof we must leave this setting.  Let us explain this reduction.  

\subsubsection*{Reduction to the strongly convex case}  If $L$ is the lineality space\footnote{the largest linear space contained in $\sigma$} of $\sigma$, then by \autoref{extremal-rays-are-rational} \ref{extremal-rays-are-rational.lineality}, we have that $L$ defined over $\bQ$, that is, $L$ has an $\mathbb{R}$-basis consisting of vectors in $M$.  Furthermore, by \autoref{lem:cone-ses}, the quotient cone in the quotient space $(M \otimes_{\mathbb Z} \mathbb{R})/L$ still defines a split-$F$-regular ring and it suffices to check that the quotient cone is generated by finitely many rational rays. Note that by killing the lineality space of $\sigma$ we make the quotient cone in $(M \otimes_{\mathbb Z} \mathbb{R})/L$ strongly convex.

In the strongly convex case, it suffices to show that $\sigma$ is a closed cone generated by finitely many rational extremal rays as then we are in the classical toric case.

\subsubsection*{Extremal rays of the closure are rational}  Let $\overline{\sigma}$ be the closure of $\sigma$ and assume $\overline{\sigma}$ is strongly convex. Strong convexity is needed in order for $\overline{\sigma}$ to have extremal rays. If \autoref{itm.IntoA} holds, then by \autoref{extremal-rays-are-rational}\ref{extremal-rays-are-rational.extremalRay} the extremal rays of $\overline{\sigma}$ are rational.  The main technical tool used in this proof is \autoref{prop:common-setup-rationality-boundary} which says that if the minimal rational subspace containing a ray $\tau$ has a certain property (which we can arrange if the ray is in the lineality space or is extremal), then $\tau$ is rational. \autoref{prop:common-setup-rationality-boundary} in turn uses a Diophatine approximation result from convex geometry (\autoref{BCHM-density}) that features in the seminal paper \cite{BirkarCasciniHaconMcKernan}. 
Note that in the two-dimensional case, one is just checking that the defining rays have rational slope, which is a straightforward exercise reminiscent of arguments found in \cite{DattaSmithFrobeniusAndValuationRings}.  

\subsubsection*{The cone is closed}
We show that $\sigma = \overline{\sigma}$ in \autoref{thm:split-F-regular-arbitrary-cone-closed}. 
A key step in our argument is the following.

\begin{enumerate}
    \setcounter{enumi}{1}
    \item If we can split a monomial coming from a lattice element in the relative interior of $\sigma$ as in \autoref{itm.IntoA}, then every rational ray of $\overline{\sigma}$ is contained in $\sigma$ and hence $k[\sigma \cap M] = k[\overline{\sigma} \cap M]$; see \autoref{rational-rays-closure}.  It follows that $k[\overline{\sigma} \cap M]$ is also split-$F$-regular.
    \label{itm.AlgebraOfSigmaIsAlgebraOfSigmaBar}    %
\end{enumerate}%
Suppose first that $\overline{\sigma}$ is strongly convex.  Since $k[\overline{\sigma} \cap M]$ is split-$F$-regular by \autoref{itm.AlgebraOfSigmaIsAlgebraOfSigmaBar}, we know the extremal rays of $\overline{\sigma}$ are rational by the previous step, and hence are all contained in $\sigma$ again by \autoref{itm.AlgebraOfSigmaIsAlgebraOfSigmaBar}. Then $\sigma = \overline{\sigma}$ because a closed strongly convex cone is generated as a convex cone by its extremal rays.

For a general (non-strongly convex) $\overline{\sigma}$, we again replace $k[\sigma \cap M]$ by $k[\overline{\sigma} \cap M]$ using \autoref{itm.AlgebraOfSigmaIsAlgebraOfSigmaBar}. We decompose $\overline{\sigma}$ into its lineality space $L$ and a strongly convex quotient cone.  As discussed above, the lineality space is rational and the quotient cone is split-$F$-regular.  \autoref{thm:split-F-regular-arbitrary-cone-closed} then shows that $\overline{\sigma}$ is also generated by rational rays, which must all be contained in $\sigma$ by \autoref{itm.AlgebraOfSigmaIsAlgebraOfSigmaBar}. Consequently, $\sigma = \overline{\sigma}$. 

\subsubsection*{Finite generation}

\begin{figure}[h!]
\includegraphics{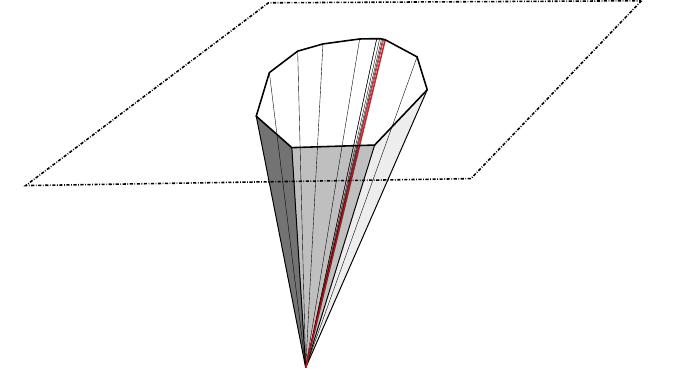}
\caption{A convex cone with accumulating extremal rays, intersected with an affine hyperplane.}
\label{fig.ConeHyperplaneSection}
\end{figure}

Reducing to the closed strongly convex case, it suffices to show that $\sigma = \overline{\sigma}$ has finitely many extremal rays.  If there are infinitely many, some subsequence of them must accumulate when intersected with a \emph{reference} rational affine hyperplane $H$, as depicted in a three-dimensional setting in \autoref{fig.ConeHyperplaneSection}, where the red ray is an accumulation of extremal rays. Here by a reference rational affine hyperplane we mean a translate of a hyperplane $\varphi^\perp$, for some rational $\varphi$ in the relative interior of the dual cone $\sigma^\vee$.

\begin{figure}[h]
\includegraphics{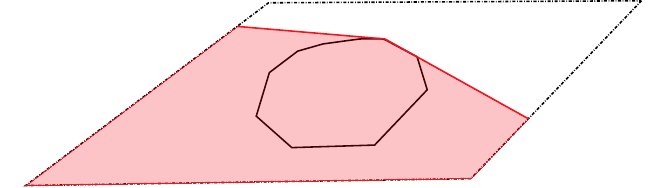}
\caption{A convex cone generated by $\sigma \cap H$ in the affine hyperplane $H$, with new chosen origin.}
\label{fig.ConeGeneratedBySigmaCapH}
\end{figure}

While \autoref{fig.ConeHyperplaneSection} is a good pictorial aid, unlike in the three-dimensional setting, an accumulation of extremal rays need not be extremal \cite[Example 4.2]{ZuffetiConesOfSpecialCyclesOfCodim2OnShimura}.  Regardless, the main accumulation result is found by combining \autoref{thm:dual-cone-rational-SFR} and \autoref{prop:extremal-rays-cannot-accumulate}.  The strategy is induction on dimension.  Indeed, if a ray $\rho$ is an accumulation of extremal rays, after restricting to a reference affine hyperplane $H$ defined over $\bQ$ like the one pictured above, one ends up with a compact polytope.  Picking an extremal point sufficiently near\footnote{In higher dimensions, one choses a vertex of the smallest face of $\sigma \cap H \cap \varphi^\perp$ for a rational supporting hyperplane $\varphi$ of $\sigma$ containing $\rho \cap H$.} the accumulation point and declaring it to be the origin, we take the cone generated by $\sigma \cap H$ in the affine hyperplane $H$.  In \autoref{fig.ConeGeneratedBySigmaCapH}, this is the cone between with two red extremal rays.  This cone is in a lower dimensionsional space!

We show that the monoid algebra defined by the cone generated by $\sigma \cap H$ in $H$ is still split-$F$-regular by \autoref{thm:affine-hyperplane-split-F-regular}. In particular, this implies that this lower-dimensional cone is closed and by induction on dimension, generated by finitely many rational rays.  That finiteness can be used to contradict the accumulation of infinitely many extremal rays in the original cone $\sigma$ via \autoref{prop:extremal-rays-cannot-accumulate}. A key fact that we use in our inductive argument is that the dual cone $\sigma^{\vee}$ is also generated by rational extremal rays when $\sigma$ is a strongly convex full-dimensional cone such that $k[\sigma \cap M]$ is split-$F$-regular (see \autoref{thm:dual-cone-rational-SFR}).

\subsection{$F$-pure-regularity}
One might also ask what happens if one weakens split-$F$-regularity to \emph{$F$-pure-regularity} in the sense of \cite{DattaSmithFrobeniusAndValuationRings} (called \emph{very strong $F$-regularity} in \cite{HashimotoF-pure-homomorphisms}, and studied as early as \cite{SmithThesis}).  Recall that a domain $R$ of characteristic $p > 0$ is called \emph{$F$-pure-regular} if for every nonzero $c \in R$, there exists an $e > 0$ such that $R \xrightarrow{c^{1/p^e}\cdot} R^{1/p^e}$ is a pure (aka universally injective) homomorphism of $R$-modules.

Note that split-$F$-regularity implies $F$-pure-regularity because split maps are pure. The converse is not true even for nice classes of regular rings \cite{DattaMurayamaF-Solid}. Thus, $F$-pure-regularity is in many ways a better adaptation of the notion of strong $F$-regularity to a non-$F$-finite (and non-Noetherian) setting than split-$F$-regularity; see \autoref{rem:difference-splitF-vs-Fpure} for more evidence. Nevertheless, in our setting, $F$-pure-regularity and split-$F$-regularity coincide and are independent of the choice of base field.  Indeed, we obtain the following:

\begin{theoremA*}[{\autoref{cor:splitF-regular-equiv-Fpure-regular}}, \cf \autoref{split-F-regular-cones}]
    Let $\sigma$ be a convex cone in $M_{\mathbb R}$ such that $0 \in \sigma$. 
	Let $p > 0$ be a prime number. Then the following are equivalent:
	\begin{enumerate}[label=(\arabic*)]
		\item There exists a field $k$ of characteristic $p$ such that $k[\sigma \cap M]$ is split-$F$-regular. 
		\item For all fields $K$ of characteristic $p$, $K[\sigma \cap M]$ is split-$F$-regular.
		\item For all fields $K$ of characteristic $p$, $K[\sigma \cap M]$ is $F$-pure-regular. 
		\item There exists a field $k$ of characteristic $p$ such that $k[\sigma \cap M]$ is $F$-pure-regular. 
		\item There exists a field $k$ of characteristic $p$ such for all $\beta \in \sigma \cap M$, there exists an integer $e_\beta > 0$ such that the unique $k[\sigma \cap M]$-linear map 
		\begin{align*}
			 k[\sigma \cap M] &\to F^{e_\beta}_*k[\sigma \cap M]\\
			1 &\mapsto F^{e_\beta}_*X^\beta
		\end{align*}
		is pure.
		\item For all $\beta\in \sigma \cap M$, there exists an integer $e_\beta > 0$ such that the set
        \[\{\gamma \in M: p^{e_\beta}\gamma + \beta \in \sigma\}\]
		is contained in $\sigma \cap M$.
	\end{enumerate}
\end{theoremA*}

\subsection{An application to graded rings associated with valuations} 
Suppose $\nu$ is a  real-valued valuation of $K/k$ that is centered on a Noetherian local ring $(R,\mathfrak m)$ that is essentially of finite type over $k$ with fraction field $K$.  In this context, Teissier introduced in \cite{TeissierValuationsToric}, the graded ring of $R$ along $\nu$, $\gr_\nu(R) \coloneqq \bigoplus_{m \in \mathbb{R}_{\geq 0}} \mathfrak{a}_{\geq m}(R)/\mathfrak{a}_{> m}(R)$ where $\fra_{\geq m}(R)$ is the ideal of $r \in R$ such that $\nu(r) \geq m$ (and similarly with $\mathfrak{a}_{> m}(R)$).  These rings have seen recent applications in algebraic $K$-stability, see \cite{LiuXuZhuangFiniteGeneration}.  
One can show that the associated graded rings are monoid algebras in the maximal rational rank case, and hence are Noetherian if one assumes they are split-$F$-regular, see \autoref{thm:finite-generation-SFR-valuation-algebra-max-rank}.


\subsection*{Acknowledgements}  The authors would like to thank Harold Blum, Holger Brenner, Izzet Coskun, Dale Cutkosky, Lawrence Ein, Han Huang, Jonathan Monta\~no,  Mircea \mustata{}, Zsolt Patakfalvi, Anurag Singh, Karen Smith, Joe Waldron, and Jakub Witaszek for valuable discussions. We especially thank Harold for insightful conversations pertaining to \autoref{thm:finite-generation-SFR-valuation-algebra-max-rank}. The authors thank the hospitality of the Centre International de Rencontres Math{\'e}matiques, the Mathematisches Forschungsinstitut Oberwolfach, the American Institute of
Mathematics, and the Simons-Laufer Mathematical Sciences Institute (formerly the Mathematical Sciences Research Institute, and which is supported by National Science Foundation Grant No. DMS-1928930) where parts of the research were carried out.
{The pictures above were generated with the help of Inkscape.}

\section{The geometry of convex cones}

\subsection*{Conventions}
All rings are commutative with a multiplicative identity.
All vector spaces are finite dimensional over $\mathbb{Q}$ or $\mathbb R$. A \textbf{lattice} $\Lambda$ in an $\mathbb{R}$-vector space $V$ is a free Abelian subgroup of $V$ of finite rank such that $\Lambda \otimes_{\mathbb Z} \mathbb{R} = V$.

In this preliminary section we recall definitions and facts about convex cones. Most sources in the commutative algebra and algebraic geometry literature focus mainly on (rational) polyhedral cones. However, we will need to work with arbitrary convex cones in finite dimensional vector spaces. Thus, we collect the material we need for ease of reference. Our primary sources are \cite[Appendix A]{deFernexEinMustataBook} (an unpublished manuscript in birational geometry) and the classic sources \cite{FenchelConvexConesSetsFunctionsBook,RockafellarConvexAnalysis}.

Throughout this section, we fix a finite dimensional $\mathbb R$-vector space $V$. 

\subsection{Convex cones}
 A \textbf{cone} $\sigma$ in $V$ is any subset such that if $x \in \sigma$, then $t x \in \sigma$ for all $t \in \mathbb{R}_{> 0}$. We say $\sigma$ is a \textbf{convex cone} if it is additionally convex, that is, 
if for all $x, y \in \sigma$, $t \in [0,1]$,
$$(1-t)x + ty \in \sigma.$$

{}

\begin{lemma}
\label{convex-cone-semigroup}
Let $\sigma$ be a subset in $V$. Then we
have the following:
\begin{enumerate}[label=(\arabic*)]
	\item $\sigma$ is a convex cone if and only if whenever $x_1, \dots, x_n 	\in \sigma$, then $t_1x_1 + \dots + t_nx_n \in \sigma$ for all $t_1, 	\dots, t_n \in \mathbb{R}_{> 0}$.  \label{convex-cone-semigroup.1}

	\item Let $\sigma$ be a convex cone. Then the linear
	 span of $\sigma$ is the set $\{x - y: x, y \in \sigma\}$.
     \label{convex-cone-semigroup.2}

	\item If $\sigma$ is a convex cone, then so is its closure 
	$\overline{\sigma}$.
    \label{convex-cone-semigroup.3}
\end{enumerate}
\end{lemma}

\begin{proof}
\ref{convex-cone-semigroup.1} follows by \cite[Cor.\ 2.6.1]{RockafellarConvexAnalysis}, \ref{convex-cone-semigroup.2} follows by \cite[Thm.\ 2.7]{RockafellarConvexAnalysis} (it is assumed that $0 \in \sigma$ in this result, but the hypothesis is not needed for \ref{convex-cone-semigroup.2}), and \ref{convex-cone-semigroup.3} follows by \cite[Thm.\ 6.2]{RockafellarConvexAnalysis}.
\end{proof}

{}

\begin{remark}
\noindent The lemma implies that a convex cone is closed under addition, but it may not contain $0$. For example, the ample cone of a projective variety is a convex cone that does not contain $0$.
On the other hand, a closed convex cone always contains $0$, and hence 
is a sub-monoid of $V$.
\end{remark}

{}

Given a subset $S \subseteq V$, the \textbf{convex cone generated
 by $S$}, denoted $\pos(S)$, is the intersection of all convex 
cones of $V$ that contain $S$. One can verify that
$$\pos(S) = \bigg{\{}\sum_{i=1}^r t_is_i : s_i \in S, t_i \in \mathbb{R}_{> 0}\bigg{\}}.$$
Note that $\pos(S)$ need not contain $0$. We say a convex cone is
 $\textbf{polyhedral}$ if it is generated by a finite set.

{}

As shown by the next result, convex polyhedral cones are usually closed.

{}

\begin{lemma}
If $\sigma$ is a convex polyhedral cone in $V$ and $0 \in \sigma$, then $\sigma$ is closed in $V$.
\end{lemma}

\begin{proof}
This follows from the proof of \cite[Lem.\ A.1.3]{deFernexEinMustataBook}. Note that
 according to our (and also \cite{deFernexEinMustataBook}'s) conventions, a convex 
polyhedral cone need not contain $0$. This is why we have to
 assume that $0 \in \sigma$ in the hypothesis of this lemma.
\end{proof}

{}

The following result will be useful in the sequel.

{}

\begin{theorem}[Carath\'eodory]
\label{Caratheodory}
If $\sigma$ is a convex cone generated by a subset $S$ of $V$, then $\sigma$ is the union of the convex cones generated by
 subsets of $S$ that are $\mathbb{R}$-linearly independent.
\end{theorem}

\begin{proof}
See \cite[Prop.\ A.5.7]{deFernexEinMustataBook}
or \cite[Cor.\ 17.1.2]{RockafellarConvexAnalysis}.
\end{proof}

{}

\subsection{Dual cones and duality}
The $\mathbb{R}$-dual space of $V$ will be denoted by $V^*$.
Evaluation gives a bilinear pairing
\[\langle \_ ,\_ \rangle : V^* \times V \rightarrow \mathbb{R}\]
The \textbf{support} of a linear form $\varphi \in V^*$ is the set
\[\Supp(\varphi) \coloneqq \{x \in V: \langle \varphi, x \rangle \geq 0\} = \{x \in V : \varphi(x) \geq 0\}.\]
It follows readily that $\Supp(\varphi)$ is a closed convex cone. We say that $\varphi \in V^*$ \textbf{supports} a cone $\sigma$ if 
$\sigma \subseteq \Supp(\varphi).$
We define
\[\varphi^{\perp} \coloneqq \{x \in V: \langle \varphi, x \rangle = 0\},\]
that is, $\varphi^{\perp}$ is the kernel of $\varphi$.
Thus, if $\varphi \neq 0$, then $\varphi^{\perp}$ is a hyperplane of $V$. Similarly, for a subset $S \subseteq V$, we define
\[
S^\perp \coloneqq \{\varphi \in V^* \colon S \subseteq \varphi^\perp\}.	
\]
If $S = \{x\}$ is a singleton, then we denote $S^\perp$ by just $x^\perp$.
Note that if $\mathbb{R}S$ is the linear span of $S$, then $S^\perp = (\mathbb{R}S)^\perp$.

{}

Given a cone $\sigma \subseteq V$, we define the \textbf{dual cone}, denoted $\sigma^{\vee}$, to be
$$\sigma^{\vee} \coloneqq \{\varphi \in V^* : \sigma \subseteq \Supp(\varphi)\}.$$

{}

\begin{lemma}
If $\sigma$ is a cone in $V$, then $\sigma^\vee = \overline{\sigma}^\vee.$
\end{lemma}

\begin{proof}
Indeed, since $\sigma \subseteq \overline{\sigma}$, we have $\overline{\sigma}^\vee \subseteq \sigma^\vee$. On the other hand, if $\sigma \subseteq \Supp(\varphi)$, then $\overline{\sigma} \subseteq \Supp(\varphi)$ because $\Supp(\varphi)$ is closed. Thus $\sigma^\vee \subseteq \overline{\sigma}^\vee$, and we obtain the desired equality.
\end{proof}

{}

As a consequence, one can deduce that the 
dual
of a cone is always a closed cone.
{}

\begin{corollary}
\label{dual-cone-closed}
If $\sigma$ is a convex cone in $V$ (not necessarily closed), then $\sigma^\vee$ is a closed convex cone in $V^*$.
\end{corollary}

\begin{proof}
One can verify that the dual of a closed convex cone is a closed 
convex cone. Since $\sigma^\vee = \overline{\sigma}^\vee$ by the 
previous lemma, and $\overline\sigma$ is a closed convex cone, we 
are done.\end{proof} 

{}

The basic duality result for convex cones is the following.

{}

\begin{proposition}
\label{duality-cones}
Let $\sigma$ be a convex cone in $V$. Then
\[\overline{\sigma} = \bigcap_{\varphi \in \sigma^\vee} \Supp(\varphi).\]
Said differently, after canonically identifying $(V^*)^*$ with $V$, we have $(\sigma^\vee)^\vee = \overline{\sigma}$.
\end{proposition}

\begin{proof}
The proof of this proposition uses the hyperplane separation theorem for convex sets. See \cite[Prop.\ A.2.1]{deFernexEinMustataBook} for the case of a closed convex cone. The statement for non-closed convex cones then follows because $\sigma^\vee = \overline{\sigma}^\vee$. Alternatively, see \cite[Thm.\ 3]{FenchelConvexConesSetsFunctionsBook}.
\end{proof}

We also get a nice characterization of the dual cone of the image of a convex cone under a surjective linear map.

\begin{proposition}
	\label{prop:dual-cones-quotient}
	Let $\pi \colon V \twoheadrightarrow W$ be a surjective linear map of finite dimensional $\mathbb R$-vector spaces. Let $\sigma$ be a convex cone in $V$ and consider the image cone $\pi(\sigma)$. Then under the injection
	$
	\pi^* \colon W^* \hookrightarrow V^*	
	$
	of dual spaces, we have $\pi^*(\pi(\sigma)^\vee) = \sigma^\vee \cap \ker(\pi)^\perp$.
\end{proposition}



\begin{proof}
If $\phi \in \ker(\pi)^\perp \subseteq V^*$, then $\phi$ induces $\varphi \in W^*$ so that $\phi = \varphi \circ \pi = \pi^*(\varphi)$. When $\phi \in \sigma^\vee \cap \ker(\pi)^\perp$, we have additionally that $\varphi(\pi(\sigma)) = \phi(\sigma) \in [0, \infty)$, so that $\varphi \in \pi(\sigma)^\vee$. It follows that $\sigma^\vee \cap \ker(\pi)^\perp \subseteq \pi^*(\pi(\sigma)^\vee)$. The reverse inclusion follows similarly: given $\psi \in \pi(\sigma)^\vee \subseteq W^*$, note that $\pi^*(\psi)(\sigma) = \psi (\pi(\sigma)) \geq 0$ and $\pi^*(\psi)(\ker(\pi)) = \psi(\pi(\ker(\pi))) = 0$. Thus, $\pi^*(\psi) \in \sigma^\vee \cap \ker(\pi)^\perp$ and the proposition follows.
\end{proof}

{}

\subsection{Strong convexity}
A convex cone is \textbf{strongly convex} if it does not contain a non-trivial linear subspace. Equivalently, a convex cone $\sigma$ is strongly convex if
whenever $x, -x \in \sigma$, then $x = 0$.

{}

\begin{lemma}
Let $\sigma$ be a convex cone in $V$. If $W$ is a subspace of $V$ contained in $\sigma$, then for all $\varphi \in \sigma^\vee$,
$W \subseteq \varphi^\perp.$
\end{lemma}

\begin{proof}
Let $\varphi \in \sigma^\vee$. We need to show that if $x \in W$, then $\varphi(x) = 0$. If not, then $\varphi(x) > 0$ since $W \subseteq \sigma \subseteq \Supp(\varphi)$. This implies $\varphi(-x) = -\varphi(x) < 0$,
 which is impossible since $W \subseteq \Supp(\varphi)$.
\end{proof}

{}

\begin{remark}
\label{strong-conv-not-closed}
The closure of a strongly convex cone need not be strongly convex. For example, consider the cone $\sigma$ in $\mathbb{R}^2$ given by 
\[
\sigma \coloneqq \{(x,y) \in \mathbb{R}^2: y > 0\}.
\]
In other words, $\sigma$ is the upper half-plane. Then $\sigma$ is strongly convex but $\overline{\sigma}$ is not since $\overline{\sigma}$ contains the line $y = 0$.
\end{remark}

{}

\begin{lemma}
    \label{lem-characterize-strong-conv}
    Let $\sigma$ be a convex cone in $V$. The following are equivalent:
    \begin{enumerate}[label=(\arabic*)]
        \item The linear span of $\sigma$ equals $V$. \label{lem-characterize-strong-conv.LinSpanIsV}
        \item The dual cone $\sigma^\vee$ is strongly convex.  \label{lem-characterize-strong-conv.DualIsStronglyConvex}
    \end{enumerate}
\end{lemma}

\begin{proof}
\ref{lem-characterize-strong-conv.LinSpanIsV} $\Rightarrow$ \ref{lem-characterize-strong-conv.DualIsStronglyConvex}: Suppose $\varphi, -\varphi \in \sigma^\vee$. We have to show $\varphi = 0$. By assumption, 
$$\sigma \subseteq \Supp(\varphi) \cap \Supp(-\varphi) = \varphi^\perp.$$
Therefore the linear span of $\sigma$ is contained in $\varphi^\perp$. But the linear span of $\sigma$ equals $V$, which means $\varphi = 0$.

\noindent \ref{lem-characterize-strong-conv.DualIsStronglyConvex} $\Rightarrow$ \ref{lem-characterize-strong-conv.LinSpanIsV}: Let $\mathbb{R}\sigma$ denote the linear span of $\sigma$. If $\mathbb{R}\sigma \neq V$, there exists a hyperplane $H$ of $V$ such that $\mathbb{R}\sigma \subseteq H$. Let $\varphi \in V^*$ such that $\ker(\varphi) = H$. Then $\varphi \neq 0$, $\ker(-\varphi) = H$ and consequently $\varphi, -\varphi \in \sigma^\vee$ since $\sigma \subseteq H$. This is impossible by strong convexity of 
$\sigma^\vee$.
\end{proof}

{}

Using duality and \autoref{lem-characterize-strong-conv}, we then have the following result.

{}

\begin{corollary}
\label{closed-strongly-convex-cone}
Let $\sigma$ be a closed convex cone in $V$. The following are equivalent:
\begin{enumerate}[label=(\arabic*)]
\item $\sigma$ is strongly convex.
\item $\sigma^\vee$ spans $V^*$.
\end{enumerate}
\end{corollary}

\begin{proof}
By \autoref{duality-cones}, we have $\sigma = (\sigma^\vee)^\vee$ since $\sigma$ is closed (after canonically identifying $V$ and $V^{**}$). The corollary now follows by \autoref{lem-characterize-strong-conv}.
\end{proof}

{}



\subsection{Faces of cones}
Let $\sigma$ be a \emph{closed} convex cone. A subset $\tau$ of $\sigma$ is called a \textbf{face of $\sigma$} if there exists $\varphi \in \sigma^\vee$ such that
$\tau  = \sigma \cap \varphi^{\perp}.$
Taking $\varphi = 0$, it follows that $\sigma$ is a face of itself. A \textbf{proper face of $\sigma$} is a face that is different from $\sigma$. 

Since $\varphi^{\perp}$ is a closed convex cone, any face of a closed convex cone is also a closed convex cone. In addition, we have the following properties of faces.

{}

\begin{lemma}
\label{extremal-nature-faces}
Let $\sigma$ be a closed convex cone in $V$. We have the following.
\begin{enumerate}[label=(\arabic*)]
	\item If $\sigma$ contains a linear subspace $L$ of 	$V$, then $L$ is contained in every face of $\sigma$.
        \label{extremal-nature-faces.1}

	\item Let $\tau$ be a face of $\sigma$. Then $	\tau$ has the property that if $x_1 + x_2 \in \tau$ for $x_1, x_2 \in 	\sigma$, then $x_1, x_2 \in \tau$. \label{extremal-nature-faces.2} 

	\item An arbitrary intersection of faces of $\sigma$
	 is also a face of $\sigma$. In fact, if $\{\tau_\alpha\}_{\alpha \in A}$ 	is a collection of faces of $\sigma$, then there exist $\alpha_1, \dots, 	\alpha_n \in A$ such that
	\[\bigcap_{\alpha \in A} \tau_\alpha = \tau_{\alpha_1} \cap \dots \cap 	\tau_{\alpha_n}.\] \label{extremal-nature-faces.3}
\end{enumerate}
\end{lemma}

\begin{proof}
\ref{extremal-nature-faces.1}  If $L = 0$ this is obvious. If not, choose a non-zero $x \in L$. Then $x$ and $-x \in \sigma$. Thus, for any $\varphi \in \sigma^\vee$, we have both 
$\varphi(x), \varphi(-x) \geq 0.$
Since $\varphi$ is linear, this can only happen if $\varphi(x) = 0$, that is, if $x \in \varphi^{\perp}$. This shows that
for any 
$\varphi \in \sigma^\vee$, $L \subseteq \varphi^\perp$, and so, $L$ is contained in any face of $\sigma$.

\noindent \ref{extremal-nature-faces.2} If $\varphi \in \sigma^\vee$ such that 
$\tau = \sigma \cap \varphi^\perp$, then $\varphi(x_1), \varphi(x_2) \geq 0$, and so, $\varphi(x_1 + x_2) = 0$ precisely
when $\varphi(x_1) = 0 = \varphi(x_2)$, that is, if 
$x_1, x_2 \in \tau$.

\noindent \ref{extremal-nature-faces.3} See \cite[Lem.\ A.3.2]{deFernexEinMustataBook}. The key point is that the intersection of a collection $\mathcal{C}$ (possibly infinite) of subspaces of a finite dimensional vector space $V$ is always equal to the intersection of finitely many subspaces from $\mathcal C$ by dimension considerations.
\end{proof}

{}

\begin{remark}
For an arbitrary closed convex cone $\sigma$, it is not true that if $\sigma'$ is a face of $\sigma$ and $\tau$ is a face of $\sigma'$, then $\tau$ is a face of $\sigma$. See \cite[E.g.\ A.4.5]{deFernexEinMustataBook} for an example. However, when $\sigma$ is polyhedral, then the face of a face of $\sigma$ is a face of $\sigma$ \cite[Prop.\ A.5.1]{deFernexEinMustataBook}.
\end{remark}

\subsection{Relative interiors of cones and their duals} 
Let $\sigma$ be a convex cone in $V$. The \textbf{relative interior} of $\sigma$, denoted $\Int(\sigma)$, is the topological interior of $\sigma$ as a subset of its affine hull. 

{}

\begin{lemma}
\label{rel-int-open-span}
Let $\sigma$ be a convex cone in $V$. Then the affine hull of $\sigma$ coincides with the linear span of $\sigma$ in $V$. Consequently, $\Int(\sigma)$ is the topological interior of $\sigma$ as a subset of its linear span. 
\end{lemma}

\begin{proof}
Recall that the linear span, $\mathbb{R}\sigma$, of $\sigma$ equals $\{x - y: x, y \in \sigma\}$ (\autoref{convex-cone-semigroup}(2)). The affine hull of $\sigma$, denoted $\Aff(\sigma)$, is the set
\[\Aff(\sigma) \coloneqq \bigg{\{}\sum_{i = 1}^n t_ix_i: x_i \in \sigma, t_i \in \mathbb{R}, \sum_{i=1}^n t_i = 1\bigg{\}}.\]
Clearly $\Aff(\sigma) \subseteq \mathbb{R}\sigma$. On the other hand, if $x, y \in \sigma$, then 
$x - y = 2(x/2) - y \in \Aff(\sigma).$
Therefore $\Aff(\sigma) = \mathbb{R}\sigma$, and consequently,
 $\Int(\sigma)$ is the topological interior of $\sigma$ in $\mathbb{R}\sigma$.
\end{proof}

{}

\begin{lemma}
\label{relative-interior-nonempty}
Let $\sigma$ be a convex cone in $V$ with a nonzero element. Then
we have the following: 
\begin{enumerate}[label=(\arabic*)]
	\item $\Int(\sigma)$ is non-empty. \label{relative-interior-nonempty.1} 
	\item If $x \in \Int(\sigma)$ and $y \in \overline{\sigma}$, 
	then $[x, y) \coloneqq \{(1-t)x + ty \colon t \in [0,1)\}$
	is contained in $\Int(\sigma)$.\label{relative-interior-nonempty.2}
	\item $\Int(\sigma) = \Int(\overline{\sigma})$. \label{relative-interior-nonempty.3}
	\item $\overline{\Int(\sigma)} = \overline{\sigma}$. Here all
	closures are in $V$. In particular, $\Int(\sigma)$ is dense in
	$\sigma$.\label{relative-interior-nonempty.4}
\end{enumerate}
\end{lemma}

\begin{proof}
\ref{relative-interior-nonempty.1}  The linear span $\mathbb{R}\sigma$ of $\sigma$ is a non-trivial 
subspace by hypothesis. Choosing a basis $\{x_1, \dots, x_r\}$ of $\mathbb{R}\sigma$ which is contained in $\sigma$, we see that
\[\{t_1x_1 + \dots + t_rx_r: t_i \in \mathbb{R}_{> 0}\}\]
is a non-empty open subset of $\mathbb{R}\sigma$ which is contained in $\sigma$.

For \ref{relative-interior-nonempty.2} see \cite[Theorem 6.1]{RockafellarConvexAnalysis}, and for \ref{relative-interior-nonempty.3} and \ref{relative-interior-nonempty.4} see \cite[Theorem 6.3]{RockafellarConvexAnalysis}.
\end{proof}

As a consequence we obtain the following result.

\begin{corollary}
	\label{cor:dense-convex-cone}
	Let $\sigma$ be a convex cone in $V$. Then $\sigma^\vee = \{0\}$ if and only if $\sigma =  V$.
\end{corollary}

\begin{proof}
	It is clear that if $\sigma = V$ that $\sigma^{\vee} = \{0\}$. Conversely, suppose $\sigma^\vee = \{0\}$. By duality \autoref{duality-cones}, $(\sigma^{\vee})^\vee = \overline{\sigma}$. Since $\sigma^\vee = \{0\}$ we get $\overline{\sigma} = V$. Then 
	\[
	\Int(\sigma) = \Int(\overline{\sigma}) = \Int(V) = V,
	\]
	where the first equality follows by \autoref{relative-interior-nonempty}\ref{relative-interior-nonempty.3}. Since $\Int(\sigma) \subseteq \sigma$, we have $\sigma = V$.
\end{proof}

{}

The relative interior of a convex cone consists of those points of the cone that are not contained in any proper face of the closure of the cone. 

{}

\begin{proposition}
\label{rel-int-as-comp}
Let $\sigma$ be a convex cone in $V$ with a nonzero element. Then
\[
	\Int(\sigma) = \overline{\sigma} - \bigcup_{\tau \subsetneq \overline{\sigma}} \tau,
\]
where the union is over all proper faces of $\overline{\sigma}$.
\end{proposition}

\begin{proof}
When $\sigma$ is a closed convex cone, that is, $\sigma = \overline{\sigma}$, the Proposition is \cite[Proposition A.3.4]{deFernexEinMustataBook}. If $\sigma$ is not necessarily closed, the
result follows by that for closed cones and \autoref{relative-interior-nonempty}\ref{relative-interior-nonempty.3} because $\Int(\sigma) = \Int(\overline{\sigma})$. The result can also be deduced from \cite[Thm.\ 11.6]{RockafellarConvexAnalysis}.
\end{proof}

{}



\begin{remark}
	\label{rem:Relint-full-dim-cones}
{\*}
	Re-interpreting \autoref{rel-int-as-comp} in 
	terms of the dual cone, we see that if $\sigma$ is a convex 
	cone of $V$ which is not contained in any hyperplane of $V$ 
	(equivalently, the linear span $\sigma$ is all of $V$), then
	\begin{equation}
		\label{eq:relint-full-dimensional}
	\Int(\sigma) = \big{\{}x \in \sigma: \forall \varphi \in 
	\sigma^\vee - \{0\}, \varphi(x) > 0 \big{\}}.
	\end{equation}
	Indeed, the hypothesis that $\sigma$ is not contained in a 
	hyperplane of $V$ implies that for any nonzero linear form 
	$\varphi \in \sigma^\vee = \overline{\sigma}^\vee$, that
	$\overline{\sigma} \cap \varphi^{\perp}$
	is a proper face of $\overline{\sigma}$. Thus the union of the
	 proper faces of $\overline{\sigma}$ is precisely the set
	$$\big{\{}x \in \overline{\sigma}: \exists \varphi \in 	\sigma^{\vee} - \{0\} \hspace{1mm} \textrm{such that} 	\hspace{1mm} \varphi(x) = 0 \big{\}}.$$
	The complement of this set, which equals $\Int(\sigma)$ by 
	\autoref{rel-int-as-comp}, is then
	$$\big{\{}x \in \sigma: \varphi(x) > 0, \forall \varphi \in \sigma^\vee - 	\{0\} \big{\}}.$$

\end{remark}

\begin{corollary}
	\label{cor:Relint-dual-cone}
	Let $\sigma$ be a closed strongly convex cone in $V$. Then
	\[
		\Int(\sigma^\vee) = \big{\{}\varphi \in \sigma^\vee \colon \varphi(x) > 0, \forall x \in \sigma - \{0\}\big{\}}.
	\]
\end{corollary}

\begin{proof}
	Since $\sigma$ is closed and strongly convex, by \autoref{closed-strongly-convex-cone}, $\sigma^\vee$ spans $V^*$. Thus, by \autoref{eq:relint-full-dimensional} in \autoref{rem:Relint-full-dim-cones} applied to $\sigma^\vee$, we have
	\[
	\Int(\sigma^\vee) = \big{\{}\varphi \in \sigma^\vee \colon \xi(\varphi) > 0, \forall \xi \in (\sigma^{\vee})^\vee - \{0\}\big{\}}.
	\]
	Since $\sigma$ is closed, by duality one can identify $(\sigma^\vee)^\vee$ with $\sigma$. Under this identification, 
	\[
		\big{\{}\varphi \in \sigma^\vee \colon \xi(\varphi) > 0, \forall \xi \in (\sigma^{\vee})^\vee - \{0\}\big{\}} = \big{\{}\varphi \in \sigma^\vee \colon \varphi(x) > 0, \forall x \in \sigma - \{0\}\big{\}}.	
	\]
	This proves the desired result.
\end{proof}

{}

\begin{corollary}
\label{rel-int-addition}
Let $\sigma$ be a convex cone in $V$. If $x \in \overline{\sigma}$ and $y \in \Int(\sigma)$, then for all $s, t \in \mathbb{R}_{> 0}$, $sx + ty \in \Int(\sigma)$. As a consequence, $\Int(\sigma)$ is also a convex cone.
\end{corollary}

\begin{proof}
If $sx+ ty \notin \Int(\sigma)$, then by \autoref{rel-int-as-comp}, there exists a proper face $\tau$ of $\overline{\sigma}$ such that $sx + ty \in \tau$. But then 
$sx, ty \in \tau$ 
by \autoref{extremal-nature-faces}, and hence
$x, y \in \tau$ 
because $\tau$ is a cone. In particular, $y \notin \Int(\sigma)$ by \autoref{rel-int-as-comp} again, contradicting our hypothesis.

Applying the above assertion to $x \in \Int(\sigma) \subseteq \overline{\sigma}$ we get that $\Int(\sigma)$ is a convex cone.
\end{proof}






{}

\subsection{Extremal subcones and lineality spaces} Let $\sigma$ be a convex cone in $V$ such that $0 \in \sigma$ and such that $\sigma$ has a nonzero element. A \textbf{ray} of $\sigma$ is a subcone of the form $\mathbb{R}_{\geq 0} \cdot x$, for some $x \in \sigma - \{0\}$. Thus, by our convention $\{0\}$ is not a ray of $\sigma$.

{}

\begin{definition}
	\label{def:extremal-subcone-ray}
	Suppose $\sigma$ is a convex cone such that $0 \in \sigma$. An \textbf{extremal subcone of $\sigma$}\footnote{We differ from \cite{deFernexEinMustataBook} who define the notion of an extremal subcone only for \emph{closed} convex cones \cite[Def.\ A.4.1]{deFernexEinMustataBook}. One reason for our choice is that the notion of the lineality space of a convex cone makes sense even for non-closed cones, and the lineality space satisfies the defining property of being extremal. Hence it is convenient to call the lineality space of a convex cone an extremal subcone even when the original cone is not closed.} is a \emph{non-empty} convex cone $\tau \subseteq \sigma$ such that for all $x, y \in \sigma$, if $x + y \in \tau$ then $x, y \in \tau$.
	An \textbf{extremal ray} of $\sigma$ is an extremal subcone which is a ray. If $\tau$ is an extremal ray of $\sigma$, then a nonzero element $x \in \tau$ is called an \textbf{extremal point} of $\sigma$.
\end{definition}

\begin{lemma}
	\label{lem:extremal-subcone-closed}
	Suppose $\sigma$ is a convex cone such that $0 \in \sigma$. Let $\tau$ be an extremal subcone of $\sigma$. Then we have the following:
	\begin{enumerate}[label=(\arabic*)]
		\item $0 \in \tau$.\label{lem:extremal-subcone-closed.1}
		\item If $\sigma$ is closed, then $\tau$ is closed.\label{lem:extremal-subcone-closed.2}
		\item If $L$ is a linear subspace of $V$ that is contained in $\sigma$, then $L \subseteq \tau$.\label{lem:extremal-subcone-closed.3}
		\item If $\tau$ is an extremal ray, then $\sigma$ is strongly convex.\label{lem:extremal-subcone-closed.4}
	\end{enumerate}
\end{lemma}

\begin{proof}
	\ref{lem:extremal-subcone-closed.1}: By definition of an extremal subcone, $\tau$ is non-empty. So let $x \in \tau$. Since $0 \in \sigma$ and $0 + x \in \tau$, it follows that $0 \in \tau$.

	\ref{lem:extremal-subcone-closed.2}: If $\Int(\tau) = \emptyset$, then $\tau = \{0\}$ by \autoref{relative-interior-nonempty}, and hence $\tau$ is closed. If $\Int(\tau) \neq \emptyset$, fix \emph{any} $x \in \Int(\tau)$. Then for all $y \in \overline{\tau}$, $y + x \in \Int(\tau) \subseteq \tau$ by \autoref{rel-int-addition}. Note that $\overline{\tau} \subseteq \sigma$ because $\sigma$ is a closed cone. Since $\tau$ is extremal, we then get $y \in \tau$ for all $y \in \overline{\tau}$. Thus, $\tau$ is closed.

	\ref{lem:extremal-subcone-closed.3}: Let $L$ be a linear subspace of $V$ that is contained in $\sigma$. Since $0 \in \tau$, it follows that $L \subseteq \tau$ since $x + (-x) \in \tau$ for all $x \in L$. 
	
	\ref{lem:extremal-subcone-closed.4}: By \ref{lem:extremal-subcone-closed.3} any linear subspace of $V$ that is contained in $\sigma$ must be contained in $\tau$. Since $\tau$ is a ray, it can only contain a $0$-dimensional linear space, and so, $\sigma$ is strongly convex.
\end{proof}

\begin{remark}
	By \autoref{extremal-nature-faces}, any face of a closed convex cone is an extremal subcone. However, the converse does not hold in general; see \cite[Example A.4.5]{deFernexEinMustataBook} which gives an example of an extremal ray of a cone which is not a face of the cone. Nevertheless, if the closed cone is polyhedral, that is, if the cone has a finite generating set, then the extremal subcones are precisely the faces of the cone \cite[Proposition A.5.1]{deFernexEinMustataBook}.
	\end{remark}

We have seen that the existence of an extremal ray implies strong convexity of the cone. For closed cones that contain a nonzero element, the converse also holds.

\begin{proposition}
\label{extremal-strongly-convex}
Let $\sigma$ be a closed convex cone in $V$ with a nonzero element. Then the following are equivalent:
\begin{enumerate}[label=(\arabic*)]
	\item $\{0\}$ is an extremal subcone of $\sigma$.
	\label{extremal-strongly-convex.1}
	\item $\{0\}$ is a face of $\sigma$.\label{extremal-strongly-convex.2}
	\item $\sigma$ is strongly convex.\label{extremal-strongly-convex.3}
	\item $\sigma$ has an extremal ray.\label{extremal-strongly-convex.4}
	\item $\sigma$ is generated as a convex cone by its extremal rays.\label{extremal-strongly-convex.5}
\end{enumerate}
\end{proposition}

\begin{proof}
\ref{extremal-nature-faces.1} $\Rightarrow$ \ref{extremal-strongly-convex.3} is a consequence of \autoref{lem:extremal-subcone-closed}\ref{lem:extremal-subcone-closed.3} and $\ref{extremal-strongly-convex.4} \Rightarrow \ref{extremal-strongly-convex.3}$ follows by \autoref{lem:extremal-subcone-closed}\ref{lem:extremal-subcone-closed.4}.

We also clearly have \ref{extremal-strongly-convex.5} $\Rightarrow$ \ref{extremal-strongly-convex.4} since we are assuming $\sigma$ has a nonzero element, while $\ref{extremal-strongly-convex.2} \Rightarrow \ref{extremal-strongly-convex.1}$ follows by \autoref{extremal-nature-faces}\ref{extremal-nature-faces.2}.

$\ref{extremal-strongly-convex.3} \Rightarrow \ref{extremal-strongly-convex.5}$ follows by \cite[Prop.\ A.4.6]{deFernexEinMustataBook} or, alternatively, by \cite[Cor.\ 18.5.2]{RockafellarConvexAnalysis}.

It remains to show \ref{extremal-strongly-convex.3} $\Rightarrow$ \ref{extremal-strongly-convex.2}. Since $\sigma$ is strongly convex, $\sigma^\vee$ spans $V^*$ by \autoref{lem-characterize-strong-conv}. Note $\dim_{\mathbb R}(V^*) \geq 1$ because $\sigma$ (hence $V$) has a nonzero element. Thus, $\Int(\sigma^\vee)$ is a non-empty open subset of $V^*$ by \autoref{relative-interior-nonempty}\ref{relative-interior-nonempty.1}. Choose $\varphi \in \Int(\sigma^\vee)$. Then for all $x \in \sigma \setminus \{0\}$, 
$
\varphi(x) > 0
$
by \autoref{cor:Relint-dual-cone}, that is, $\sigma \cap \varphi^{\perp} = \{0\}$. Thus, $\{0\}$ is a face of $\sigma$. 
\end{proof}



Using the notion of an extremal subcone, we next examine when the image of a convex cone under a surjective linear map of vector spaces stays strongly convex.

\begin{proposition}
	\label{prop:image-strongly-convex}
	Let $\sigma$ be a convex cone in $V$ such that $0 \in \sigma$. Let $\tau$ be an extremal subcone of $\sigma$ and let $W \coloneqq \mathbb{R}\tau$ be the linear span of $\tau$. If
	$
	\pi \colon V \twoheadrightarrow V/W	
	$
	is the projection map, then $\pi(\sigma)$ is a strongly convex cone.
\end{proposition}

\begin{proof}
	Since $\tau$ is an extremal subcone of $\sigma$, we first claim that
	\begin{equation}
		\label{eq:intersection-span-extremal}
	\sigma \cap W = \tau.	
	\end{equation}
	We already know that $\tau \subseteq \sigma \cap W$ because both $\sigma$ and $W = \mathbb{R}\tau$ contain $\tau$. Now let $z \in \sigma \cap W$. By the fact that $z \in W$, we can write $z = z_1 - z_2$ for $z_1, z_2 \in \tau$, and so, $z + z_2 = z_1 \in \tau$. Since $z, z_2 \in \sigma$ and $\tau$ is extremal, we get $z \in \tau$. This shows $\sigma \cap W \subseteq \tau$, and so, $\sigma \cap W = \tau$.

	That $\pi(\sigma)$ is a convex cone is a consequence of $\pi$ being a linear map. Suppose for contradiction that $\pi(\sigma)$ is not strongly convex. Then there exists $x \in \sigma \setminus W$ such that 
	\[
	x + W, -x + W \in \pi(\sigma).	
	\]
	Since $-x + W \in \pi(\sigma)$, there exists $y \in \sigma$ such that 
	$
	-x + W = y + W.	
	$
	Thus, $x + y \in W$. Moreover, $x + y \in \sigma$ because $x, y \in \sigma$, and so, $x + y \in \tau$ by \autoref{eq:intersection-span-extremal}. Since $\tau$ is extremal, we get $x \in \tau \subseteq W$, which contradicts our choice of $x$.
\end{proof}

\begin{definition}
	\label{def:lineality-space}
	The \textbf{lineality space} of a convex cone $\sigma$ that contains $0$, denoted $\Lin(\sigma)$, is the largest linear space that is contained in $\sigma$.
\end{definition}

\begin{lemma}
	\label{lem:properties-lineality-space}
	Let $\sigma$ be a convex cone such that $0 \in \sigma$. Then we have the following:
	\begin{enumerate}[label=(\arabic*)]
		\item $\Lin(\sigma) = \sigma \cap -\sigma$.\label{lem:properties-lineality-space.1}
		\item $\Lin(\sigma)$ is an extremal subcone of $\sigma$. Moreover, $\Lin(\sigma)$ is the intersection of all extremal subcones of $\sigma$.\label{lem:properties-lineality-space.2}
	\end{enumerate}
\end{lemma}

\begin{proof}
	\ref{lem:properties-lineality-space.1} is clear and we omit the proof.

	\ref{lem:properties-lineality-space.2}: Suppose $x, y \in \sigma$ such that $x + y \in \sigma \cap -\sigma$. Then $-x - y \in \sigma$ as well, and so, $-x = (-x - y) + y \in \sigma$. Thus, $-x \in \sigma \cap -\sigma$, or equivalently, $x \in \sigma \cap -\sigma$. We similarly get $y \in \sigma \cap -\sigma$. Thus, $\Lin(\sigma)$ is an extremal subcone of $\sigma$.

	It remains to show that any extremal subcone of $\sigma$ contains $\Lin(\sigma)$. So let $\tau$ be any extremal subcone of $\sigma$. Since $0 \in \tau$ by \autoref{lem:extremal-subcone-closed}\ref{lem:extremal-subcone-closed.1}, it follows that $\tau$ contains any linear subspace of $V$ that is contained in $\sigma$. In particular, then $\Lin(\sigma) \subseteq \tau$. 
\end{proof}

As a special case of \autoref{prop:image-strongly-convex}, we obtain the following:

\begin{proposition}
	\label{prop:lineality-quotient-strong-convex}
	Let $\sigma$ be a convex cone in $V$ such that $0 \in \sigma$. Let $\Lin(\sigma)$ be the lineality space of $\sigma$ and $\pi \colon V \twoheadrightarrow V/\Lin(\sigma)$ be the projection map. 
	Then we have the following:
	\begin{enumerate}[label=(\arabic*)]
		\item $\pi(\sigma)$ is a strongly convex cone.\label{prop:lineality-quotient-strong-convex.1}
		\item If $\sigma$ is closed, then $\pi(\sigma)$ is a closed strongly convex cone. \label{prop:lineality-quotient-strong-convex.2}
	\end{enumerate}
\end{proposition}

\begin{proof}
	\ref{prop:lineality-quotient-strong-convex.1} follows by \autoref{prop:image-strongly-convex} since $\Lin(\sigma)$ is an extremal subcone of $\sigma$ by \autoref{lem:properties-lineality-space}\ref{lem:properties-lineality-space.2}.
	
	\ref{prop:lineality-quotient-strong-convex.2}: By \ref{prop:lineality-quotient-strong-convex.1} it remains to check that $\pi(\sigma)$ is closed. The topology on $V/\Lin(\sigma)$ is the quotient topology induced by $\pi$. Thus, it suffices to check that 
	\[
	\pi^{-1}(\pi(\sigma)) = \sigma + \Lin(\sigma)	
	\]
	is closed in $V$. However, $\Lin(\sigma )\subseteq \sigma$, and so, $\sigma + \Lin(\sigma) = \sigma$, and the latter is closed by assumption.
\end{proof}

In the proof of \autoref{thm:dual-cone-rational-SFR} we will also need the following fact about lineality spaces and extremal subcones.

\begin{lemma}
	\label{lem:lineality-sum-extremal-subcone}
	Let $\sigma$ be a convex cone such that $0 \in \sigma$. Let $\tau \subseteq \sigma$ be an extremal subcone and let $W \coloneqq \mathbb{R}\tau$ be the linear span of $\tau$. Then 
	$
	\Lin(\sigma +  W) = W.	
	$
\end{lemma}

\begin{proof}
	Let $x \in \Lin(\sigma + W)$, that is, $x \in (\sigma + W) \cap (-\sigma + W)$. Let $y_1, y_2 \in \sigma$ and $z_1, z_2 \in W$ such that
	\[
	y_1 + z_1 = x = -y_2 + z_2.	
	\]
	Then $y_1 + y_2 = z_2 - z_1 \in W$ and $y_1 + y_2 \in \sigma$. Thus, $y_1 + y_2 \in \sigma \cap W \stackrel{\autoref{eq:intersection-span-extremal}}{=} \tau$, where the equality holds because $\tau$ is extremal. Again, by the extremal nature of $\tau$, we get $y_1, y_2 \in \tau$, and so, $x = y_1 + z_1 \in W$. This shows that $\Lin(\sigma + W) \subseteq W$. The containment $W \subseteq \Lin(\sigma + W)$ is clear since $\Lin(\sigma + W)$ is the largest subspace of $V$ that is contained in $\sigma + W$ and clearly $W \subseteq \sigma + W$.
\end{proof}

Even though an extremal ray of a closed convex cone is not a face of the cone in general, such a ray is always contained in a proper of face of the cone except in the case where the cone is a single ray. This and some other properties of extremal rays are highlighted in the next proposition.

{}

\begin{proposition}
\label{extremal-ray-properties}
Let $\sigma$ be a closed convex cone in $V$ and let $\tau$ be an extremal ray of $\sigma$. Let $\mathcal{L}$ be the linear span of $\sigma$. Then we have the following:
\begin{enumerate}[label=(\arabic*)]
\item $\dim_{\mathbb R} \mathcal{L} = 1$ if and only if 
$\sigma = \tau$.\label{extremal-ray-properties.dimlineality1}

\item If $\dim_{\mathbb R} \mathcal{L} \geq 2$, then $\tau$ is contained in a proper face of $\sigma$.\label{extremal-ray-properties.dimlinealitygeq2}

\item If $\dim_{\mathbb R} V \geq 2$, then there exists a nonzero linear functional $\varphi \in \sigma^{\vee}$ such that 
$\tau \subseteq \sigma \cap \varphi^\perp$.\label{extremal-ray-properties.dimambientspacegeq2}
\end{enumerate}
\end{proposition}

\begin{proof}
Since $\sigma$ has an extremal ray, we know that it is strongly convex by \autoref{extremal-strongly-convex}. The proof 
of the forward implication in (1) then follows by the strong
convexity of $\sigma$, while the backward implication follows 
because a ray spans a one-dimensional space.

The assertion in (3) follows from (1) and (2). Indeed, if $\tau$ is contained in a proper face $F$ of $\sigma$, then expressing $F$ as $\sigma \cap \varphi^\perp$, for some $\varphi \in \sigma^\vee$, we note that $\varphi \neq 0$ because $F$ is a proper face. On the other hand, if $\sigma = \tau$, then using the fact that $\dim_{\mathbb R} V \geq 2$, we note that $\sigma$ is contained in some hyperplane of $V$. This hyperplane is realized as the kernel of a nonzero linear functional $\varphi$, which consequently contains $\sigma$ in its support.

Thus, it suffices to verify (2). Without loss of generality assume $\mathcal{L} = \mathbb{R}^{\oplus d}$, where $d \geq 2$ by hypothesis. Since $\sigma$ is a closed cone, the assertion follows via \autoref{rel-int-as-comp} if we can show that 
$\Int(\sigma) \cap \tau = \emptyset$. 
Recall that $\Int(\sigma)$ is the topological interior of $\sigma$ in $\mathcal L$ by \autoref{rel-int-open-span}, that is, $\Int(\sigma)$ is an open subset of $\mathcal L$. Therefore for any $x \in \Int(\sigma)$, by enclosing $x$ in an open ball $B$ of $\mathcal{L}$ centered at $x$ such that $B \subseteq \Int(\sigma)$, we see that $x$ lies on a line segment between two points $y, y' \in B$ such that $y$ and $y'$ are not on the ray spanned by $x$ (here we are using the fact that $d \geq 2$). Thus, $x \notin \tau$ since $x \neq 0$ and $\tau$ is extremal.
\end{proof}

{}




{}

\begin{remark}
It is important that the dimension of $V$ is $\geq 2$ in \autoref{extremal-ray-properties}\ref{extremal-ray-properties.dimambientspacegeq2}. Otherwise, one may take $V = \mathbb{R}$ and $\sigma = R_{\geq 0}$. Then $\sigma$ is an extremal ray of itself, and there is no nonzero $\varphi \in \sigma^\vee$ such that $\sigma \subseteq \varphi^{\perp}$.
\end{remark}

{}


{}







{}


{}

\subsection{Diophantine approximation}

Throughout this subsection we fix a finite dimensional 
$\mathbb Q$-vector space $V_{\mathbb Q}$ and let 
\[
V \coloneqq V_{\mathbb Q} \otimes_{\mathbb Q} \mathbb R.
\]
We view  $V_{\mathbb{Q}}$ as a subset of $V$. 
A vector $v \in V$ is \textbf{rational} if $v \in V_{\mathbb Q}$. An $\mathbb R$-linear subspace $W$ of $V$ is \textbf{defined over the rationals} if it is spanned by rational vectors of $V$. By a \textbf{rational hyperplane of $V$} we mean a hyperplane of $V$ which is defined over the rationals.

More generally, an affine subspace $A \subseteq V$ is \textbf{defined over the rationals} if $A$ is the affine hull of a set of rational vectors of $V$. In other words, $A$ is defined over the rationals if there exists $S \subseteq V_{\mathbb Q}$ such that
$$A = \bigg{\{}\sum^n_{i = 1}\lambda_is_i : n \in \mathbb{Z}_{> 0}, s_i \in S, \lambda_i \in \mathbb{R}, \sum_{i=1}^n \lambda_i = 1\bigg{\}}.$$

If $\sigma$ is a convex cone in $V$ containing $0$, then a ray $\tau$ of $\sigma$ is \textbf{rational} if there exists $x \in V_{\mathbb Q}$ such that $\tau = \mathbb{R}_{\geq 0}\cdot x$.


\begin{lemma}
\label{irrational-lines}
Let $V_{\mathbb Q}$ be a finite dimensional $\mathbb Q$-vector space of dimension $n$ and $V = V_{\mathbb Q} \otimes_{\mathbb Q} \mathbb R$. 
\begin{enumerate}[label=(\arabic*)]
	\item If $W$ is a subspace of $V$ and $W_{\mathbb Q} \coloneqq W \cap V_{\mathbb Q}$, then $W$ is defined over the rationals if and only if $W = W_{\mathbb Q} \otimes_{\mathbb Q} \mathbb R$. \label{irrational-lines.1}
	\item A one-dimensional vector subspace $L$ of $V$ is not defined over the rationals if and only if $L \cap V_{\mathbb Q} = \{0\}$.\label{irrational-lines.2}
	\item Suppose $Q$ is a subspace of $V$. Then there exists a smallest subspace $W$ of $V$ defined over the rationals such that $Q \subseteq W$.\label{irrational-lines.3}
	\item Suppose $e_1, \dots, e_n \in V_{\mathbb Q}$ is a basis of rational vectors of $V$ over $\mathbb R$. A vector $v \in V$ is not contained in any proper $\mathbb R$-vector subspace of $V$ defined over the rationals if and only if when we express $v$ in terms of the basis $\{e_1,\dots,e_n\}$ as
	\[v = \sum_{i=n}^n x_ie_i,\]
	then the real numbers $x_1, \dots, x_n$ are linearly independent over $\mathbb Q$.\label{irrational-lines.4}
	\item Let $Q$ be a vector subspace of $V$ such that $Q$ is not contained in any proper subspace of $V$ that is defined over the rationals. Then there exists $v \in Q$ such that $v$ is not contained in any proper \emph{affine} subspace of $V$ that is defined over the rationals.\label{irrational-lines.5}
\end{enumerate}
\end{lemma}

\begin{proof}
\ref{irrational-lines.1} and \ref{irrational-lines.2} follow from the definition of a subspace defined over the rationals, so we omit their proofs. 

  \ref{irrational-lines.3}: Let $\Sigma$ be the collection of subspaces $W$ of $V$ defined over the rationals such that $Q \subseteq W$. Then 
  \begin{equation}
	\label{eq:IF}
\bigcap_{W \in \Sigma} W = \bigcap_{W \in \Sigma} (W_{\mathbb Q} \otimes_{\mathbb Q} \mathbb{R}) = \bigg{(}\bigcap_{W \in \Sigma} W_{\mathbb Q}\bigg{)} \otimes_{\mathbb Q} \mathbb{R}	
  \end{equation}
  is defined over the rationals as well. Note here that by dimension considerations, the possibly infinite intersection $\bigcap_{W \in \Sigma} W$ is actually a finite intersection (alternatively, since $V$ is an Artinian $\mathbb R$-module, an arbitrary intersection of $\mathbb R$-submodules of $V$ always equals an intersection of a finite collection of the submodules). Thus, the second equality  in \autoref{eq:IF} follows because (flat) base change to $\mathbb{R}$ commutes with finite intersections. Then $\bigcap_{W \in \Sigma} W$ is the smallest subspace of $V$ defined over the rationals that contains $Q$.

  \ref{irrational-lines.4}: For the forward implication of \ref{irrational-lines.4}, we prove the contrapositive. If $x_1, \dots, x_n$ are linearly
dependent over $\mathbb Q$, then one can, without loss of 
generality, express $x_1$ as a $\mathbb Q$-linear combination of
$x_2,\dots, x_n$, say 
$$x_1 = q_2x_2 + \dots + q_nx_n,$$
where $q_i \in \mathbb{Q}$. Then
$v = x_2(q_2e_1 + e_2) + \dots + x_n(q_ne_1 + e_n)$ is a vector
in the proper $\mathbb{R}$-subspace of $V$ spanned by the $n-1$ rational vectors 
$q_2e_1 + e_2, \dots, q_ne_1 + e_n$.  

Conversely, suppose $x_1, \dots, x_n$ are linearly independent over $\mathbb Q$. If $v$ is  contained in a proper $\mathbb R$-subspace of $V$ defined over the rationals, then there exists $w_1, \dots, w_m \in V_{\mathbb Q}$ and real numbers $r_1, \dots, r_m$, with $m < n$, such that $v = r_1w_1 + \dots + r_mw_m$. Expressing each $w_i$ as a $\mathbb{Q}$-linear combination of the basis $\{e_1, \dots, e_n\}$ of $V_{\mathbb Q}$, it follows by linear independence of $\{e_1, \dots, e_n\}$ that $x_1, \dots, x_n$ is in the $\mathbb{Q}$-linear span of $r_1, \dots, r_m$. This is impossible since $m < n$, completing the proof of \ref{irrational-lines.4}.

\ref{irrational-lines.5}: Assume that no such $v \in Q$ exists, that is, assume every $v \in Q$ is contained in some proper affine subspace of $V$ defined over the rationals. Since $V_{\mathbb Q}$ is countable, the set $\Sigma$ of proper affine subspaces of $V$ defined over the rationals is countable. Note that for any $W \in \Sigma$, $Q \cap W \neq Q$. Otherwise, $W$ would be subspace of $V$ since it would contain $0$, which contradicts the assumption that $Q$ is not contained in any proper subspace of $V$ that is defined over the rationals. However, by our assumption,
\begin{equation}
	\label{eq:Baire-cat-contradiction}
	Q = \bigcup_{W \in \Sigma} Q \cap W.
\end{equation}
Note that $Q \cap W$ is a proper affine subspace of $Q$, so its complement is dense in $Q$. Then \autoref{eq:Baire-cat-contradiction} is impossible by the Baire category theorem since in a finite dimensional vector space over $\mathbb R$, a countable intersection of open dense sets is still dense.
\end{proof}

{}

The main result of this subsection is the following Diophantine 
approximation result.

{}

\begin{proposition}{\cite[Lem.\ 3.7.6]{BirkarCasciniHaconMcKernan}}
\label{BCHM-density}
Let $V$ be a finite-dimensional $\mathbb{R}$-vector space with the
property that there exists a $\mathbb{Q}$-vector space $V_{\mathbb Q}$ such that
$V = V_{\mathbb Q} \otimes_{\mathbb Q} \mathbb R.$
Let $\Lambda$ be a lattice in $V$ generated by rational vectors. If $v \in V$ is not contained in any proper affine subspace of $V$ which is defined over the rationals, then the set
$$\{mv + \lambda: m \in \mathbb{Z}_{\geq 0}, \lambda \in \Lambda\}$$
is dense in $V$.
\end{proposition}


{}

We draw the following conclusion from \autoref{BCHM-density}.

{}

\begin{corollary}
\label{BCHM-density-consequence}
Let $V$ be a finite-dimensional $\mathbb{R}$-vector space with the
property that there exists a $\mathbb{Q}$-vector space $V_{\mathbb Q}$ such that $V = V_{\mathbb Q} \otimes_{\mathbb Q} \mathbb R.$ Let $\Lambda$ be a lattice in $V$ generated by rational vectors. If 
$$\tau \coloneqq R_{\geq 0}\cdot x$$ 
is a ray in $V$ which is not contained in any proper vector subspace of $V$ that is defined over the rationals, then the set
$$\tau + \Lambda \coloneqq \{v + \lambda: v \in \tau, \lambda \in \Lambda\}$$
is dense in $V$. 
\end{corollary}

\begin{proof}
By \autoref{BCHM-density}, it suffices to show that there exists $v \in \tau$ such that $v$ is not contained in any proper affine subspace of $V$ that is defined over the rationals. Let $Q = \mathbb{R}\tau$ be the $\mathbb R$-span of $\tau$. Then $Q$ is not contained in any proper vector subspace of $V$ that is defined over the rationals. By \autoref{irrational-lines}\ref{irrational-lines.5}, there exists $v \in Q$ such that $v$ is not contained in any proper affine subspace of $V$ that is defined over the rationals. Then replacing $v$ by $-v$ if necessary, we may assume $v \in \tau$.\hfill \qedhere
\end{proof}

\section{Finite generation of split-{$F$}-regular monoid algebras}
Let $R$ be a commutative ring of prime characteristic $p > 0$ and $e > 0$ be an integer. The \textbf{$e$-th iterate of the (absolute) Frobenius on $R$} is the ring homomorphism
\begin{align*}
	F^e \colon R &\to R\\
	r &\mapsto r^{p^e}.
\end{align*}
The target copy of $R$ with $R$-module structure given by restriction of scalars along $F^e$ is denoted $F^e_*R$. Furthermore, for $x \in R$, the corresponding element of $F^e_*R$ is denoted $F^e_*x$. Thus, in this notation the $R$-module structure on $F^e_*R$ is given, more explicitly, as follows: for $r \in R$ and $F^e_*x \in F^e_*R$, we have
\[
r \cdot F^e_*x = F^e_*r^{p^e}x.	
\]
When $e = 1$, $F^1_*R$ is just denoted $F_*R$.

We say $R$ is \textbf{$F$-finite} if $F_*R$ is a finite $R$-module. We say that $R$ is \textbf{Frobenius split} if $F \colon R \to F_*R$ splits as a map of $R$-modules. Note that if $R$ is Frobenius split, then for all integers $e > 0$, $F^e \colon R \to F^e_*R$ also splits.

\subsection{Background on split-{$F$}-regularity and {$F$}-pure-regularity} 
We now introduce the notions of split-$F$-regularity and $F$-pure-regularity for integral domains. 

\emph{We will not assume that our rings are Noetherian or $F$-finite.}

\begin{definition}
	\label{def:split-F-regular}
	Let $R$ be an integral domain of characteristic $p > 0$. Then $R$ is \textbf{split-$F$-regular} if for all $c \in R \setminus \{0\}$ there exists an integer $e > 0$ such that the unique $R$-linear map 
	\begin{align}
		\label{split-F-regular}
	\lambda_{c,e} \colon R &\to F^e_*R \nonumber\\
	1 &\mapsto F^e_*c	
	\end{align}
	splits as a map of $R$-modules.
\end{definition}

\begin{remark}
	When $R$ is Noetherian and $F$-finite, then \autoref{def:split-F-regular} is due to Hochster and Huneke \cite{HochsterHunekeTC1,HochsterHunekeTightClosureAndStrongFRegularity}. In this setting the notion is called \emph{strongly $F$-regular}. 
\end{remark}

We collect some basic facts about this notion.

\begin{lemma}
	\label{lem:basic-facts-SFR}
Let $R$ be a split-$F$-regular domain of characteristic $p > 0$. Then we have the following:
\begin{enumerate}[label=(\arabic*)]
	\item $R$ is Frobenius split.\label{lem:basic-facts-SFR.1}
	\item Any localization of $R$ is split-$F$-regular. \label{lem:basic-facts-SFR.2}
	\item Let $c \in R \setminus \{0\}$ and suppose $\lambda_{c,e}$ is as in \autoref{split-F-regular} above. If $\lambda_{c,e}$ splits, then for all integers $f > 0$, there exists an $R$-linear map $F^{e+f}_*R \to R$ that sends $F^{e+f}_*c \mapsto 1$.\label{lem:basic-facts-SFR.3}
	\item $R$ is a direct summand of any finite extension of $R$. Hence, $R$ is a splinter.\label{lem:basic-facts-SFR.4}
	\item $R$ is integrally closed in its fraction field. \label{lem:basic-facts-SFR.5}
	\item Suppose $S$ is a direct summand of $R$. Then $S$ is split-$F$-regular. \label{lem:basic-facts-SFR.6}
\end{enumerate}
\end{lemma}

\begin{proof}
	\ref{lem:basic-facts-SFR.1} Choose an integer $e \gg 0$ and an $\lambda_{1,e} \colon R \to F^e_*R$ that sends $1 \mapsto F^e_*1$ splits. But $\lambda_{1,e}$ is precisely $F^e$, the $e$-iterate of Frobenius on $R$. Since $F^e \colon R \to F^e_*R$ factors through $F \colon R \to F_*R$, it follows that $F$ splits as well. Hence $R$ is Frobenius split.

	\ref{lem:basic-facts-SFR.2} Let $W \subset R$ be a multiplicative set and $c/w \in W^{-1}R$ be a nonzero element. Choose 
	\begin{align*}
		\phi_{c,e} \colon F^e_*R &\to R\\
		F^e_*c &\mapsto 1.	
	\end{align*}
	Then $W^{-1}\phi_{c,e}$ is a $W^{-1}R$-linear map $F^e_*W^{-1}R \to W^{-1}R$ that sends $F^e_*(c/1) \mapsto 1$ (we are using that $F^e_*$ commutes with localization). The composition
	\[
	F^e_*W^{-1}R \xrightarrow{F^e_*(s/1) \cdot} F^e_*W^{-1}R \xrightarrow{W^{-1}\phi_{c,e}} W^{-1}R 	
	\]
	is then a $W^{-1}R$-linear map such that $F^e_*(c/s) \mapsto 1$.

	\ref{lem:basic-facts-SFR.3} Let $\phi_{c,e} \colon F^e_*R \to R$ be a splitting of $\lambda_{c,e}$, that is $\phi_{c,e}(F^{e}_*c) = 1$. Then for all integers $f > 0$, we get an $R$-linear map 
	\begin{align*}
	F^f_*(\phi_{c,e}) \colon F^{e + f}_*R &\to F^f_*R\\
	F^{e+f}_*c &\mapsto F^f_*1.	
	\end{align*}
	Since $R$ is Frobenius split, the $f$-th iterate of Frobenius on $R$ also splits. So let $\varphi \colon F^f_*R \to R$ be such a Frobenius splitting. Then $\varphi \circ F^f_*(\phi_{c,e}) \colon F^{e+f}_*R \to R$ maps $F^{e+f}_*c \mapsto 1$.

	\ref{lem:basic-facts-SFR.4} Suppose $R \hookrightarrow T$ is a module finite ring extension. Then it is well-known that $\Hom_R(T,R)\neq 0$ (see, for instance, \cite[Ex.\ 5.1.3]{DattaMurayamaF-Solid}). Consequently, there exists an $R$-linear map $\varphi \colon T \to R$ such that $\varphi(1) \neq 0$.
	Let $c \coloneqq \varphi(1)$. Since $R$ is a domain and split-$F$-regular, there exists an integer $e > 0$ and an $R$-linear map 
	$\phi_{c,e} : F^e_*R \rightarrow R$
	such that $\phi_{c,e}(F^e_*c) = 1$. Then the composition
	\[T \xrightarrow{F^e_T} F^e_*T \xrightarrow{F^e_*\varphi} F^e_*R \xrightarrow{\phi_{c,e}} R\]
	is an $R$-linear map which sends $1$ to $1$, that is, $R \hookrightarrow T$ splits.

	\ref{lem:basic-facts-SFR.5} Let $K$ be the fraction field of $R$. Suppose $a/b \in K$ be integral over $R$. Then the extension $R \hookrightarrow R[a/b]$ is finite, and hence it splits by \ref{lem:basic-facts-SFR.4}. Let $\varphi$ be a splitting. Then 
	\[
	a = a\varphi(1) = \varphi(a) = \varphi(b\cdot a/b) = b\varphi(a/b).	
	\] 
	Since $\varphi(a/b) \in R$, we get $b|a$, that is, $a/b \in R$. Thus $R$ is integrally closed in $K$.

	\ref{lem:basic-facts-SFR.6} We may assume $S \subseteq R$. Thus, $S$ is also a domain. Let $\varphi \colon R \to S$ be a splitting of $S \subseteq R$, and let $c \in S$ be a nonzero element. Choose an $R$-linear map $\phi_{c,e} \colon F^e_*R \to R$ that maps $F^e_*c \mapsto 1$. Then the composition
	\[
	F^e_*S \hookrightarrow F^e_*R \xrightarrow{\phi_{c,e}} R \xrightarrow{\varphi} S	
	\]
	is an $S$-linear map that sends $F^e_*c \to 1$. Thus, $S$ is split-$F$-regular.
\end{proof}

We will now introduced a weaker variant of split-$F$-regularity that is a more natural notion to study in the non-Noetherian or non-$F$-finite setting. The definition (\autoref{def:F-pure-regularity}) will be obtained by replacing the splitting of the map $\lambda_{c,e}$ in \autoref{split-F-regular} by its purity. 

Recall that a map of modules $M \to N$ over a ring $A$ is \textbf{pure} or \textbf{universally injective} if for all $A$-modules $P$, the induced map $M \otimes_A P \to N \otimes_A P$ is injective. Split maps are pure, and we have the following partial converse.

\begin{lemma}\cite[\href{https://stacks.math.columbia.edu/tag/058L}{Tag 058L}]{stacks-project}
	\label{lem:pure-implies-split}
	Let $A$ be a ring and $\varphi \colon M \to N$ be an $A$-linear map. Suppose that $\varphi$ is pure and $\coker(\varphi)$ is a finitely presented $A$-module. Then $\varphi$ splits.
\end{lemma}

We now introduce the pure variant of split-$F$-regularity. We will again defined the notion only for integral domains. For a definition in the non-domain case, see \cite{DattaSmithFrobeniusAndValuationRings}.

\begin{definition}
	\label{def:F-pure-regularity}
	Let $R$ be an integral domain of characteristic $p > 0$. Then $R$ is \textbf{$F$-pure-regular} if for all $c \in R \setminus \{0\}$ there exists an integer $e > 0$ such that the unique $R$-linear map
	\begin{align*}
		\lambda_{c,e} \colon R &\to F^e_*R\\
		1 &\mapsto F^e_*c
	\end{align*}
	is pure/universally injective.
\end{definition}

Since split maps are pure, it is clear that a split-$F$-regular domain is also $F$-pure-regular. However, the two notions do not coincide even in nice settings; see \autoref{rem:difference-splitF-vs-Fpure}.

\begin{remark}
	$F$-pure-regularity was first studied in the Noetherian local setting in Smith's thesis \cite{SmithThesis}. The notion was then considered for arbitrary Noetherian rings in \cite{HashimotoF-pure-homomorphisms}, where it is called \emph{very strongly $F$-regular}. The reason behind Hashimoto's terminology is that for Noetherian rings, $F$-pure-regularity is a strengthening of a variant of strong $F$-regularity which was defined for arbitrary Noetherian rings by Hochster using tight closure theory \cite{HochsterFoundations}. The $F$-pure-regular terminology was first used in \cite{DattaSmithFrobeniusAndValuationRings}, where the first author and Smith studied this notion for non-Noetherian rings (in particular, valuation rings). More recently, this notion has been studied in depth by Hochster and Yao \cite{HochsterYaoGenericLocal} for excellent rings. Since tight closure theory has been primarily explored in the Noetherian setting, we prefer the $F$-pure-regular terminology over the very strongly $F$-regular terminology.  It is important not to confuse this with \emph{purely $F$-regular} singularities, which are the positive characteristic analog of purely log terminal singularities, see for example \cite[Definition 3.1]{TakagiInversion} or \cite{TakagiAdjointIdealsAndACorrespondence}.
\end{remark}

\autoref{lem:basic-facts-SFR} has an analog for $F$-pure-regularity. With the exception of the splinter property, the analogs of \autoref{lem:basic-facts-SFR}\ref{lem:basic-facts-SFR.1}-\ref{lem:basic-facts-SFR.6} may be found in \cite[Prop.\ 6.1.3, Rem.\ 6.1.5, Prop.\ 6.2.2]{DattaSmithFrobeniusAndValuationRings}. We include the proof of the splinter property for $F$-pure-regularity here as it does not seem to appear in the literature.

\begin{lemma}
	\label{lem:F-pure-regular-splinter}
	Let $R$ be an $F$-pure-regular ring of characteristic $p > 0$. Then any finite ring extension $R \hookrightarrow T$ is pure.
\end{lemma}

\begin{proof}
	Suppose $R \hookrightarrow T$ is a module finite ring extension. Then it is well-known that $\Hom_R(T,R)\neq 0$ (see, for instance, \cite[Ex.\ 5.1.3]{DattaMurayamaF-Solid}). Consequently, there exists an $R$-linear map $\varphi \colon T \to R$ such that $\varphi(1) \neq 0$.
	Let $c \coloneqq \varphi(1)$. Choose an integer $e > 0$ such that the map 
	\begin{align*}
	\lambda_{c,e} \colon R &\to F^e_*R\\	
	1 &\mapsto F^e_*c
	\end{align*}
	is $R$-pure. Then the composition
	\[
	R \hookrightarrow T	\xrightarrow{F^e} F^e_*T \xrightarrow{F^e_*\varphi} F^e_*R
	\]
	equals $\lambda^e_c$, which was chosen to be $R$-pure. Thus, $R\hookrightarrow T$ is $R$-pure as well.
\end{proof}

\begin{remarks}
	\label{rem:difference-splitF-vs-Fpure}
	We record some differences between split-$F$-regularity and $F$-pure-regularity.
	\begin{enumerate}
		\item If $R$ is a split-$F$-regular domain, then $R$ is a direct summand of every finite ring extension. On the other hand, if $R$ is $F$-pure-regular, then we only get purity of finite extensions of $R$ and not their splitting. Of course this distinction ceases to exist by \autoref{lem:pure-implies-split} if the extension is finite and finitely presented.
		\item Every regular local ring is $F$-pure-regular \cite[Thm.\ 6.2.1]{DattaSmithFrobeniusAndValuationRings}. In particular, a DVR is always $F$-pure-regular. However, in \cite[Cor.\ 7.4.3]{DattaMurayamaF-Solid} it shown that a split-$F$-regular DVR is necessarily excellent. Note that non-excellent DVRs abound in positive characteristic, even in the function field of $\mathbb{P}^2_{\overline{\mathbb{F}_p}}$ \cite{DattaSmithExcellence}.
		\item Split-$F$-regularity and $F$-pure-regularity do not coincide in the class of excellent Henselian regular local rings by recent examples in \cite{DattaMurayamaTateNotF-split}. 
	\end{enumerate}
\end{remarks}

\autoref{rem:difference-splitF-vs-Fpure} indicates that it is more natural to study $F$-pure-regularity instead of split-$F$-regularity for non-Noetherian or non-$F$-finite rings. However, we will show that for the monoid algebras we care about in this paper, the distinction between split-$F$-regularity and $F$-pure-regularity is superfluous (see \autoref{split-F-regular-cones}, \autoref{cor:splitF-regular-equiv-Fpure-regular}).

\subsection{From monoids to convex cones}\label{subsec:monoids-to-cones}
Let $L$ be a free Abelian group of finite rank and $S$ be a submonoid of $L$. We let $\mathbb{Z}S$ denote the $\mathbb Z$-submodule of $L$ generated by $S$. Note $\mathbb{Z}S$ is also a free Abelian group of finite rank. If 
$$\mathbb{Z}S \cong \mathbb{Z}^{\oplus d}$$
we call $d$ the \textbf{rank} of the monoid $S$.  We let $\mathbb{R}S$ denote the $\mathbb R$-vector space $\mathbb{Z}S \otimes_{\mathbb Z} \mathbb R$ and similarly $\mathbb{Q}S \coloneqq \mathbb{Z}S \otimes_{\mathbb Z} \mathbb Q$. Both vector spaces have dimension $d$. Let $\sigma_S$ be the convex cone of $\mathbb{R}S$ generated by $S$ and $\overline{\sigma_S}$ the closure of $\sigma_S$ in $\mathbb{R}S$. Thus,
\[
\sigma_S = \bigg{\{}\sum^n_{i = 1}a_is_i: n \in \mathbb{Z}_{> 0}, a_i \in \mathbb{R}_{\geq 0}, s_i \in S\bigg{\}}.
\]
Note that the $\mathbb{R}$-span of $\sigma_S$ is $\mathbb{R}S$, that is, $\sigma_S$ is a full-dimensional cone in $\mathbb{R}S$.

{}

\begin{lemma}
	\label{cone-rational-intersection}
	With notation as above, 
	\[
	\mathbb{Q}S \cap \sigma_S = \bigg{\{}\sum^n_{i = 1}a_is_i: n \in \mathbb{Z}_{> 0}, a_i \in \mathbb{Q}_{\geq 0}, s_i \in S\bigg{\}}.
	\]
	\end{lemma}
	
	\begin{proof}
	The inclusion 
	\[
	\bigg{\{}\sum^n_{i = 1}a_is_i: n \in \mathbb{Z}_{> 0}, a_i \in \mathbb{Q}_{\geq 0}, s_i \in S\bigg{\}} \subseteq \mathbb{Q}S \cap \sigma_S
	\]
	is clear. Let $x \in \mathbb{Q}S \cap \sigma_S$. We want to show that $x$ is in the $\mathbb{Q}_{\geq 0}$-span of $S$. Since $\sigma_S$ is the convex cone generated by $S$, by Carath\'eodory's theorem (\autoref{Caratheodory}), there exists an $\mathbb{R}$-linearly independent subset $\{s_1, \dots, s_k\} \subseteq S$ and $b_1, \dots b_k \in \mathbb{R}_{\geq 0}$
	such that
	\begin{equation}
		\label{eq:*}
	x = b_1s_1 + \dots + b_ks_k.
	\end{equation}
	Since $x$ is also an element of $\mathbb{Q}S$, extending $\{s_1,\dots,s_k\}$ to a $\mathbb{Q}$-basis $\{s_1, \dots, s_k, t_1, \dots, t_{d-k}\}$ of $\mathbb{Q}S$, there exist $a_1, \dots, a_d \in \mathbb{Q}$
	such that
	\begin{equation}
		\label{eq:**}
	x = a_1s_1 + \dots + a_ks_k + a_{k+1}t_1 + \dots + a_dt_{d-k}.
	\end{equation}
	But $\{s_1, \dots, s_k, t_1, \dots, t_{d-k}\}$ is also an $\mathbb{R}$-linearly independent subset of $\mathbb{R}S = \mathbb{Q}S \otimes_{\mathbb Q} \mathbb R$. Thus, comparing the coefficients in \autoref{eq:*} and \autoref{eq:**}, we get
	$
	a_1 = b_1, \dots, a_k = b_k
	$ 
	and 
	$
		a_{k+1} = \dots = a_{d} = 0.		
	$ 
	Thus, $b_1, \dots, b_k \in \mathbb{Q}_{\geq 0}$ and $x$ is in the $\mathbb{Q}_{\geq 0}$-span of $s_1, \dots, s_k$.
	\end{proof}

{}

\begin{center}
\noindent\fbox{%
    \parbox{\textwidth-75pt}{
From now until the end of the paper, $k$ will denote a  field of characteristic $p > 0$.
	}}
\end{center}

Let $k[S]$ be the monoid algebra on $S$ over $k$. Since $S \subseteq \mathbb{Z}S$, $k[S]$ is a $k$-subalgebra of the Noetherian $k$-algebra $k[\mathbb{Z}S]$. Moreover, $k[S]$ and $k[\mathbb{Z}S]$ have the same fraction field, hence 
\[\td{\Frac(k[S])/k} = d.\]

The main question that we will investigate in this paper is:

\begin{question}
\label{ques-split-monoid-algebra}
If $k[S]$ is split-$F$-regular, is $k[S]$ a finitely generated $k$-algebra, or equivalently, is $S$ a finitely generated monoid?
\end{question}

We will see in \autoref{thm:split-F-regular-monoid-finite-generation} that \autoref{ques-split-monoid-algebra} has an affirmative answer. But this will require some work.

We first draw some preliminary conclusions from split-$F$-regularity of $k[S]$.

\begin{proposition}
\label{normal-semigroup}
Suppose $k[S]$ is split-$F$-regular. Then we have the following: 
\begin{enumerate}[label=(\arabic*)]
\item $k[S]$ is a splinter. \label{normal-semigroup.1}
\item $k[S]$ is integrally closed in its fraction field. \label{normal-semigroup.2}
\item For all $\alpha \in \mathbb{Z}S$, if $n\alpha \in S$ for some integer $n > 0$, then $\alpha \in S$. In other words, $S$ is a normal\footnote{Normal submonoids are also called \emph{saturated} submonoids in the literature. For instance, see \cite[A.6]{deFernexEinMustataBook}.} submonoid of $\mathbb{Z}S$.\label{normal-semigroup.3}
\end{enumerate}
\end{proposition}

\begin{proof}
\ref{normal-semigroup.1} follows by \autoref{lem:basic-facts-SFR}\ref{lem:basic-facts-SFR.4} and \ref{normal-semigroup.2} by \autoref{lem:basic-facts-SFR}\ref{lem:basic-facts-SFR.5}.

\ref{normal-semigroup.3} Suppose $\alpha \in \mathbb{Z}S$ such that $n\alpha \in S$ for some integer $n > 0$. Since $\Frac(k[S]) = \Frac(k[\mathbb{Z}S])$, we see that $X^{\alpha} \in \Frac(k[S])$ and $(X^\alpha)^n \in k[S]$. Consequently $X^\alpha \in k[S]$ since $k[S]$ is integrally closed in $\Frac(k[S])$ by \ref{normal-semigroup.2}. Thus, $\alpha \in S$.
\end{proof}

{}

\begin{proposition}
\label{integrally-closed-cones}
If $k[S]$ is integrally closed in its fraction field, then $\mathbb{Z}S \cap \sigma_S = S$.
\end{proposition}

\begin{proof}
It suffices to show $\mathbb{Z}S \cap \sigma_S \subseteq S$. We have $\mathbb{Z}S \cap \sigma_S = \mathbb{Z}S \cap (\mathbb{Q}S \cap \sigma_S)$, and by \autoref{cone-rational-intersection},
\[
\mathbb{Q}S \cap \sigma_S = \bigg{\{}\sum^n_{i = 1}a_is_i: n \in \mathbb{Z}_{\geq 0}, a_i \in \mathbb{Q}_{\geq 0}, s_i \in S\bigg{\}}.
\]
Hence if $\alpha \in \mathbb{Z}S \cap \sigma_S$ and $\sum_i a_is_i$ is a representation of $\alpha$ as above with $a_i \in \mathbb{Q}_{\geq 0}$, upon clearing the denominators of all the nonzero $a_i$, there exists a positive integer $m$ such that
$m\alpha \in S.$
But the proof of \autoref{normal-semigroup}\ref{normal-semigroup.3} then implies $\alpha \in S$ because $k[S]$ is integrally closed in its fraction field. 
\end{proof}

Using the aforementioned results, in order to examine \autoref{ques-split-monoid-algebra} we can reduce to the following setup from convex geometry.

\begin{center}
\noindent\fbox{%
    \parbox{\textwidth-65pt}{
\begin{setup}
\label{setup}
	\parbox[t]{335pt}{
Let $M$ be a free Abelian group of finite rank $d \geq 1$ (the lattice) and let 
\begin{center}
$\displaystyle{
	\textrm{$M_{\mathbb Q} \coloneqq M \otimes_{\mathbb Z} \mathbb Q$ and $M_\mathbb{R} \coloneqq M \otimes_{\mathbb Z} \mathbb R = M_{\mathbb Q} \otimes_{\mathbb Q} \mathbb R$.}}
	$
\end{center}
\noindent Fix a convex cone $\sigma$ in $M_\mathbb{R}$, \emph{not necessarily closed}, such that 

	\begin{enumerate}
\item $0 \in \sigma,$ and
\item the $\mathbb R$-linear span of $\sigma$ is all of $M_{\mathbb R}$.
\end{enumerate}
}
\end{setup}
	}}
\end{center}

{}

Suppose $\sigma$ satisfies the hypotheses of \autoref{setup}. We are interested in examining how split-$F$-regularity of the monoid algebra $k[\sigma \cap M]$ is related to the geometry of $\sigma$. In particular, we will examine the following variant of \autoref{ques-split-monoid-algebra}:

{}

\begin{question}
\label{revised-question}
If $k[\sigma \cap M]$ is split-$F$-regular, is $k[\sigma \cap M]$ Noetherian?
\end{question}

{}

We note that the question will have an affirmative answer if $\sigma$ is a closed cone generated by \emph{finitely many rational rays}. 

{}

\begin{remarks}
\label{some-remarks}
{\*}
\begin{enumerate}
\item Suppose $S$ is a submonoid of a lattice $L$. If $M = \mathbb{Z}S$ and $\sigma_S$ is the convex cone of $\mathbb{R}S = M_\mathbb{R}$ generated by $S$, then the pair $(\mathbb{Z}S, \sigma_S)$ satisfies the hypotheses of \autoref{setup}. In addition, if the monoid algebra $k[S]$ is normal (which it is when $k[S]$ is split-$F$-regular), then $k[S] = k[\sigma_S \cap \mathbb{Z}S]$ by \autoref{integrally-closed-cones}. Thus, an affirmative answer to \autoref{revised-question} under the hypotheses of Setup \ref{setup} also provides an affirmative answer to \autoref{ques-split-monoid-algebra}. Actually \autoref{setup} is more general because the cone $\sigma_S$ has the additional property that it is generated by rational rays by construction, whereas no rationality assumptions are a priori made on the cone in \autoref{setup}. \label{some-remarks.a}
\smallskip

\item \autoref{setup} implies that the dual cone $\sigma^\vee = \overline{\sigma}^\vee$ is strongly convex in the dual space $M_{\mathbb R}^*$ by \autoref{duality-cones} and \autoref{closed-strongly-convex-cone}. If in addition we assume that $\overline{\sigma}$ is strongly convex, then $\overline{\sigma}^\vee$ will span $M_{\mathbb R}^*$.\label{some-remarks.b}
\smallskip

\item Since $\sigma$ is closed under scaling by elements of $\mathbb{R}_{> 0}$, if $\alpha \in M$ and $n > 0$ is a positive integer such that $n\alpha \in \sigma \cap M$, then $\alpha \in \sigma \cap M$. In other words, the monoid $\sigma \cap M$ is automatically saturated in $M$ (c.f. \autoref{normal-semigroup}\ref{normal-semigroup.3}).\label{some-remarks.c}
\smallskip

\item By the hypotheses of \autoref{setup}, the linear span of $\sigma$ is $M_{\mathbb R}$. Since $M_{\mathbb R}$ is an $\mathbb R$-vector space of dimension $d \geq 1$, $\sigma$ must contain a nonzero element of $M_{\mathbb R}$. Furthermore, by \autoref{rel-int-open-span}
$\Int(\sigma)$ 
coincides with the topological interior of $\sigma$ as a subset of $M_{\mathbb R}$. Thus, $\Int(\sigma)$ is open in $M_{\mathbb R}$ and it non-empty by \autoref{relative-interior-nonempty}\ref{relative-interior-nonempty.1}. \label{some-remarks.d}

\smallskip

\item Let $\mathbb{Z}(\sigma \cap M)$ be the subgroup of $M$ generated by $\sigma \cap M$. Then the following are equivalent:
\begin{enumerate}[label = (\arabic*)]
	\item  $\Int(\sigma)$ is a non-empty open set in $M_{\mathbb R}$.\label{some-remarks.e.i}
	\item $\mathbb{Z}(\Int(\sigma) \cap M) = M.$\label{some-remarks.e.ii}
	\item $\sigma$ is full-dimensional in $M_{\mathbb R}$.\label{some-remarks.e.iii}
\end{enumerate}
\ref{some-remarks.e.i} $\Rightarrow$ \ref{some-remarks.e.ii}: $M_\mathbb{Q}$ is dense in $M_{\mathbb R}$ and $\Int(\sigma)$ is a non-empty open subset of $M_{\mathbb R}$. Scaling an element of $\Int(\sigma) \cap M_{\mathbb Q}$ by a large enough integer, it follows that there exists  $\alpha \in \Int(\sigma) \cap M$. As $\Int(\sigma)$ is open in the metric space $M_{\mathbb R}$, for any $e \in M$, there exists an integer $m \gg 0$, 
\[\alpha + \frac{1}{m}e \in \Int(\sigma).\]
Scaling by $m$, we get $m\alpha + e \in \Int(\sigma) \cap M$ (see \autoref{rel-int-addition}), and so, $e = (m\alpha + e) - m\alpha \in \mathbb{Z}(\Int(\sigma) \cap M)$. 

\noindent \ref{some-remarks.e.ii} $\Rightarrow$ \ref{some-remarks.e.iii}: $M_{\mathbb R}$ is spanned as a $\mathbb R$-vector space by $M$. Since every element of $M = \mathbb{Z}(\Int(\sigma) \cap M)$ can be expressed as a difference of two elements of $\Int(\sigma) \cap M$, it follows that $M_{\mathbb R}$ is also spanned as a $\mathbb R$-vector space by $\Int(\sigma) \cap M \subseteq \sigma$.

\noindent \ref{some-remarks.e.iii} $\Rightarrow$ \ref{some-remarks.e.i}: $\Int(\sigma)$ is non-empty by \autoref{relative-interior-nonempty}\ref{relative-interior-nonempty.1} since $\dim_{\mathbb R}(M_{\mathbb R}) \geq 1$ and open in $M_{\mathbb R}$ by \autoref{rel-int-open-span}.
\label{some-remarks.e}

\smallskip

\item Suppose $M$, $M_{\mathbb Q}$ and $M_{\mathbb R}$ are as in \autoref{setup}. Let $\sigma$ be a cone in $M_{\mathbb R}$ such that $\sigma \cap M$ has a nonzero lattice point and such that $0 \in \sigma$. Our main object of interest will be the $k$-algebra, $k[\sigma \cap M]$, determined by the monoid $\sigma \cap M$. For the study of $k[\sigma \cap M]$, we can always reduce to \autoref{setup}. Consider the subgroup 
\[
N \coloneqq \mathbb{Z}(\sigma \cap M)
\] 
of $M$ generated by $\sigma \cap M$. This is also a finite free Abelian group of rank $\geq 1$. Let 
\[
N_{\mathbb R} \coloneqq \mathbb{R}(\sigma \cap M)
\] 
be the subspace of $M_{\mathbb R}$ that is spanned by $\sigma \cap M$. Note that $N_{\mathbb R}$ can be identified with $N \otimes_{\mathbb Z} \mathbb{R}$. Then 
\[
\tau \coloneqq \sigma \cap N_{\mathbb R}
\] 
is a cone in $N_{\mathbb R}$ and one can check that
$
\tau \cap N = \sigma \cap M.
$  
The $\mathbb{R}$-span of $\tau$ is $N_{\mathbb R}$ because $N_{\mathbb R}$ is also the $\mathbb R$-span of $\sigma \cap M$. 
Replacing $M$ by $N$, $M_{\mathbb R}$ by $N_{\mathbb R}$ and $\sigma$ by $\tau$ we see that the monoid algebra $k[\sigma \cap M]$ remains unchanged and that we are now in the situation of \autoref{setup}.\label{some-remarks.f}
\end{enumerate}
\end{remarks}

{} 

\subsection{Splitting versus purity and independence of ground field}
\label{sec:Splitting versus purity and independence of ground field} 
{Our first observation is that the $k$-algebra $k[\sigma \cap M]$ is always Frobenius split.

\begin{proposition}
    \label{prop.SemiGroupAlgebrasAreFSplit}
Let $\sigma$ be a convex cone in $M_{\mathbb R}$ such that $0 \in \sigma$. Then the monoid algebra $k[\sigma \cap M]$ is always Frobenius split.
\end{proposition}

\begin{proof}
Fix a Frobenius splitting 
$
\phi: F_*k \rightarrow k$
of the ground field. Then $\phi$ extends to a Frobenius splitting 
$
	\varphi: F_*k[M] \rightarrow k[M]
	$ 
as follows: for $s \in k$ and $\alpha \in M$ we send
\[   
 F_*(s X^{\alpha}) \mapsto 
     \begin{cases}
       \phi(F_*s) X^{\alpha/p} &\quad\text{if $\alpha/p \in M$,}\\
       0 &\quad\text{otherwise.} \\ 
     \end{cases}
\]
Every element of $k[\sigma \cap M]$ is a sum of terms of the form $s X^{\beta}$, where $\beta \in \sigma \cap M$ and $s \in k$. Since $\sigma$ is closed under scaling, the restriction of $\varphi$ to $F_*k[\sigma \cap M]$ maps each $F_*(s X^{\beta})$ into $k[\sigma \cap M]$. Consequently, ${\varphi}$ maps $F_*k[\sigma \cap M]$ into $k[\sigma \cap M]$, and so, $k[\sigma \cap M]$ is Frobenius split.
\end{proof}

We will now try to understand what restrictions are placed on the geometry of $\sigma$ if $k[\sigma \cap M]$ is split-$F$-regular. A key observation for the rest of the paper is:

\begin{proposition}
\label{split-F-regular-cones}
Let $\sigma$ be a convex cone in $M_{\mathbb R}$ such that $0 \in \sigma$. Fix $\alpha \in \sigma \cap M$ and an integer $e_\alpha > 0$. The following are equivalent: 
\begin{enumerate}[label=(\arabic*)]
\item The unique $k[\sigma \cap M]$-linear map 
\begin{align*}
\lambda_{X^\alpha, e_\alpha} \colon k[\sigma \cap M] &\to F^{e_\alpha}_*k[\sigma \cap M]\\
1 &\mapsto F^{e_\alpha}_*X^\alpha
\end{align*}
splits.\label{split-F-regular-cones.split}

\item The unique $k[\sigma \cap M]$-linear map 
\begin{align*}
\lambda_{X^\alpha, e_\alpha} \colon k[\sigma \cap M] &\to F^{e_\alpha}_*k[\sigma \cap M]\\
1 &\mapsto F^{e_\alpha}_*X^\alpha
\end{align*}
is pure. \label{split-F-regular-cones.pure}

\item For all integers $e \geq e_\alpha$,
\[M_{\alpha, e} \coloneqq \{\beta \in M: p^e\beta + \alpha \in \sigma\}\]
is contained in $\sigma \cap M$. \label{split-F-regular-cones.monoid}

\item The set $M_{\alpha, e_\alpha} \coloneqq \{\beta \in M: p^{e_\alpha}\beta + \alpha \in \sigma\}$ is contained in $\sigma \cap M$. \label{split-F-regular-cones.monoid.one-e}
\end{enumerate}
\end{proposition}

Note that the sets $M_{\alpha, e}$ have appeared in several contexts when studying splitting of monoid algebras, see for example \cite{DeStefaniMontanoNunezBetancourtPurity}.

We will use the following lemma in the proof of \autoref{split-F-regular-cones}.

\begin{lemma}
	\label{lem:torsion-free-quotient}
	Let $\sigma$ be a convex cone in $M_{\mathbb R}$ such that $0 \in \sigma$. Then $M/\mathbb{Z}(\sigma \cap M)$ is a torsion-free Abelian group.
\end{lemma}

\begin{proof}[Proof of \autoref{lem:torsion-free-quotient}]
	Let $\gamma \in M$ and $n \in \mathbb{Z}$ be a nonzero integer such that $n\gamma \in \mathbb{Z}(\sigma \cap M)$. We want to show $\gamma  \in \mathbb{Z}(\sigma \cap M)$. Replacing $n$ by $-n$ if necessary, we may assume $n > 0$. Thus, there exist $\beta_1, \beta_2 \in \sigma \cap M$ such that 
	$
	n\gamma = \beta_1 - \beta_2.	
	$
	Then $n\gamma + \beta_2 = \beta_1 \in \sigma$. Consequently, $(n\gamma + \beta_2) + (n-1)\beta_2 \in \sigma$ as well since $\sigma$ is a convex cone (here we use that $n - 1 \geq 0$). In other words, $n(\gamma + \beta_2) \in \sigma$. Since $\sigma$ is a cone, this gives us $\gamma + \beta_2 \in \sigma$. Clearly, $\gamma + \beta_2 \in M$ as well since $\gamma, \beta_2 \in M$. Thus, $\gamma + \beta_2 \in \sigma \cap M$, and so, $\gamma  = (\gamma + \beta_2) - \beta_2 \in \mathbb{Z}(\sigma \cap M)$, as desired.
\end{proof}

\begin{proof}[Proof of \autoref{split-F-regular-cones}]
	\ref{split-F-regular-cones.split} $\Rightarrow$ \ref{split-F-regular-cones.pure} because split maps are pure.

	\ref{split-F-regular-cones.pure} $\Rightarrow$ \ref{split-F-regular-cones.monoid}:  It is well-known that for all integers $e \geq e_\alpha$ the $k[\sigma \cap M]$-linear map 
	\begin{align*}
	\lambda_{X^\alpha,e} \colon k[\sigma \cap M] &\to F^e_*k[\sigma \cap M]\\
	1 &\mapsto F^e_*X^\alpha	
	\end{align*}
	is pure if $\lambda_{X^\alpha, e_\alpha}$ is pure. See for instance \cite[Rem.\ 6.1.5]{DattaSmithFrobeniusAndValuationRings}.

	So fix an integer $e \geq e_\alpha$ and let $\beta \in M_{\alpha,e}$. We want to show that $\beta \in \sigma \cap M$.

	\begin{claim}\label{claim:first}
		$\beta \in \mathbb{Z}(\sigma \cap M)$.
	\end{claim}
	
	\begin{proof}[Proof of \autoref{claim:first}]
		By the defining property of $M_{\alpha, e}$, we have $p^e\beta + \alpha \in \sigma$. Then we also have $p^e\beta + \alpha \in \sigma \cap M$ because $\beta, \alpha \in M$. This shows
		\[
		p^e\beta = (p^e\beta + \alpha) - \alpha	\in \mathbb{Z}(\sigma \cap M).
		\]
		Since $M/\mathbb{Z}(\sigma \cap M)$ is a torsion-free Abelian group by \autoref{lem:torsion-free-quotient}, we get $\beta \in \mathbb{Z}(\sigma \cap M)$. This proves \autoref{claim:first}. \phantom\qedhere
	\end{proof}

	By \autoref{claim:first} there exist $\beta_1, \beta_2 \in \sigma \cap M$ such that
	\begin{equation}
		\label{eq:beta-difference}
	\beta = \beta_1 - \beta_2.
	\end{equation}
	We claim that the Frobenius map $F^e \colon k[\sigma \cap M] \to F^e_*k[\sigma \cap M]$ factors via the $k[\sigma \cap M]$-algebra
	\[
	S \coloneqq \frac{k[\sigma \cap M][T_1,T_2]}{(T_1^{p^e} - X^\alpha, T_2^{p^e} - X^{p^e\beta +\alpha}, X^{\beta_1}T_1 - X^{\beta_2}T_2)},	
	\]
	where $T_1, T_2$ are indeterminates. Indeed, let $\overline{T_1}$ (resp. $\overline{T_2}$) denote the class of $T_1$ (resp. $T_2$) in $S$. Then the $k[\sigma \cap M]$-algebra factorization map $\pi \colon S \to F^e_*k[\sigma \cap M]$ sends
	\begin{align*}
	\overline{T_1} &\mapsto F^e_*X^\alpha\\
	\overline{T_2} &\mapsto F^e_*X^{p^e\beta +\alpha},
	\end{align*}
	and $\pi$ acts as the $e$-th iterate of Frobenius on the elements of $k[\sigma \cap M]$.
	
	Consider the $k[\sigma \cap M]$-linear map
	\begin{align*}
	\widetilde{\lambda_e} \colon k[\sigma \cap M] &\to S\\
	1 \mapsto \overline{T_1},
	\end{align*}
	Then $\pi \circ \widetilde{\lambda_e} = \lambda_{X^\alpha, e}$ by construction, and since $\lambda_{X^\alpha, e}$ is $k[\sigma \cap M]$-pure, we get $\widetilde{\lambda_e}$ is also $k[\sigma \cap M]$-pure. In particular, $\widetilde{\lambda_e}$ is injective. Note that $S$ is a finite and finitely presented $k[\sigma \cap M]$-algebra (finiteness follows since $S$ is integral over $k[\sigma \cap M]$). Therefore, $S$ is a finitely presented $k[\sigma \cap M]$-module by \cite[\href{https://stacks.math.columbia.edu/tag/0564}{Tag 0564}]{stacks-project}.
	Then by \cite[\href{https://stacks.math.columbia.edu/tag/0519}{Tag 0519}(4)]{stacks-project}, $\coker(\widetilde{\lambda_e})$ is a finitely presented $k[\sigma \cap M]$-module. Consequently, $\widetilde{\lambda_e}$ splits as a map of $k[\sigma \cap M]$-modules by \autoref{lem:pure-implies-split}. Let $\phi \colon S \to k[\sigma \cap M]$ be a splitting of $\widetilde{\lambda_e}$, that is, $\phi$ is a $k[\sigma \cap M]$-linear map such that $\phi(\overline{T_1}) = 1$.

	Note that in $S$ we have
	$
	X^{\beta_1}\overline{T_1} = X^{\beta_2}\overline{T_2}.	
	$ 
	Applying $\phi$ to both sides of this equation and using $k[\sigma \cap M]$-linearity gives us 
	$
	X^{\beta_1} = X^{\beta_2}\phi(\overline{T_2}).	
	$ 
	Since $\phi(\overline{T_2}) \in k[\sigma \cap M]$, then in the ring $k[\mathbb{Z}(\sigma \cap M)]$, we have
	\[
	X^{\beta} \stackrel{\autoref{eq:beta-difference}}{=} X^{\beta_1 - \beta_2} = \phi(\overline{T_2}) \in k[\sigma \cap M].	
	\]
	Thus, we must have $\beta \in \sigma \cap M$, as desired.

	\ref{split-F-regular-cones.monoid} $\Rightarrow$ \ref{split-F-regular-cones.monoid.one-e} is clear.

	\ref{split-F-regular-cones.monoid.one-e} $\Rightarrow$ \ref{split-F-regular-cones.split}: Suppose $\phi \colon F^{e_\alpha}_*k \to k$ is a Frobenius splitting of the $e_\alpha$-th iterate of the Frobenius on $k$. Then as in the proof of \autoref{prop.SemiGroupAlgebrasAreFSplit}, we get a Frobenius splitting
	$
	\pi \colon F^{e_\alpha}_* k[M] \to k[M]	
	$
	of the $e_\alpha$-th iterate of Frobenius on $k[M]$ given by the following rule: for $s \in k$ and $\beta \in M$,
	\[   
 	\pi(F^{e_\alpha}_*(s X^\beta)) = 
     \begin{cases}
       \phi(F^{e_\alpha}_*s) X^{\beta/p^{e_\alpha}} &\quad\text{if $\beta/p^{e_\alpha} \in M$},\\
       0 &\quad\text{otherwise.} \\ 
     \end{cases}
\]
Note that since $\alpha \in \sigma \cap M$, we have $-\alpha \in M$. Then we get a $k[M]$ (hence also $k[\sigma \cap M]$)-linear map
\[
\pi_{-\alpha} \coloneqq F^{e_\alpha}_*k[M] \xrightarrow{F^{e_\alpha}_*X^{-\alpha}\cdot} F^{e_\alpha}_*k[M] \xrightarrow{\pi} k[M]	
\]
that sends $F^{e_\alpha}_* X^\alpha \mapsto 1$. It suffices to show that $\pi_{-\alpha}$ restricted to $F^{e_\alpha}_*k[\sigma \cap M]$ maps into $k[\sigma \cap M]$. So suppose $\gamma \in \sigma \cap M$. Then for any $s \in k$, we get
\[
	\pi_{-\alpha}(F^{e_\alpha}_*s X^{\gamma}) =
	\begin{cases}
	  \phi(F^{e_\alpha}_*s) X^{(\gamma - \alpha)/p^{e_\alpha}} &\quad\text{if $(\gamma - \alpha)/p^{e_\alpha} \in M$,}\\
	  0 &\quad\text{otherwise.} 
	\end{cases}
\]
Suppose $\gamma \in \sigma \cap M$ such that $(\gamma - \alpha)/p^{e_\alpha} \in M$. Then
$(\gamma - \alpha)/p^{e_\alpha} \in M_{\alpha,e_\alpha} \stackrel{\ref{split-F-regular-cones.monoid.one-e}}{\subseteq} \sigma \cap M$.	
Thus, $\pi_{-\alpha}(F^{e_\alpha}_*s X^{\gamma}) \in k[\sigma \cap M]$ for all $s \in k$ and $\gamma \in \sigma \cap M$, or equivalently, $\pi_{-\alpha}$ maps $F^{e_\alpha}_*k[\sigma \cap M]$ into $k[\sigma \cap M]$. Therefore, upon restricting the domain and codomain of $\pi_{-\alpha}$ we obtain a $k[\sigma \cap M]$-linear map $F^{e_\alpha}_*k[\sigma \cap M] \to k[\sigma \cap M]$ that splits $\lambda_{X^\alpha, e_\alpha}$.
\end{proof}

\begin{remark}
\label{rem:independence-field}
    Notice that in \autoref{split-F-regular-cones}, conditions \ref{split-F-regular-cones.monoid} and \ref{split-F-regular-cones.monoid.one-e} depend only on $\sigma$, $M$ and the characteristic $p$ and not on the field $k$ of characteristic $p$ that we form the monoid algebra over. Said differently, if there exists a field $k$ of characteristic $p$ such that the map $k[\sigma \cap M] \xrightarrow{1 \mapsto F^e_*X^\alpha} F^e_*k[\sigma \cap M]$ splits for some $\alpha \in \sigma \cap M$, then for \emph{any} field $K$ of characteristic $p$, the map $K[\sigma \cap M] \xrightarrow{1 \mapsto F^e_*X^\alpha} F^e_*K[\sigma \cap M]$ splits as well.
\end{remark}

The proof of \autoref{split-F-regular-cones}\ref{split-F-regular-cones.monoid.one-e} $\Rightarrow$ \autoref{split-F-regular-cones}\ref{split-F-regular-cones.split} shows that if 
\begin{align*}
	\lambda_{X_\alpha, e_\alpha} \colon k[\sigma \cap M] &\to F^{e_\alpha}_*k[\sigma \cap M]\\
	1 &\mapsto F^{e_\alpha}_*X^\alpha
\end{align*}
splits, then up to fixing a splitting $F^{e_\alpha}_*k \to k$ of the $e_\alpha$-th iterate of the ground field $k$, there is a canonical $k[M]$-linear map
\begin{align*}
F^{e_\alpha}_*k[M] &\to k[M]\\
F^{e_\alpha}_*X^\alpha &\mapsto 1
\end{align*}
that restricts to a $k[\sigma \cap M]$-linear map $F^{e_\alpha}_*k[\sigma \cap M] \to k[\sigma \cap M]$. We single out this observation in the next result.

\begin{lemma}
	\label{lem:split-implies-torus-compatible-split}
	Let $\sigma$ be a convex cone in $M_\mathbb{R}$ such that $0 \in \sigma$. Let $\alpha \in \sigma \cap M$ and $e_\alpha > 0$ be an integer such that the $k[\sigma \cap M]$-linear map 
	\begin{align*}
		\lambda_{X^\alpha, e_\alpha} \colon k[\sigma \cap M] &\to F^{e_\alpha}_*k[\sigma \cap M]\\
		1 &\mapsto F^{e_\alpha}_*X^\alpha
	\end{align*}
	splits. Let $\phi \colon F^{e_\alpha}_*k \to k$ be a splitting of the $e_\alpha$-th iterate of Frobenius on $k$. Define a group homomorphism
	$
	\pi_{-\alpha} \colon F^{e_\alpha}_*k[M] \to k[M]	
	$
	by the rule that for $s \in k$ and $\gamma \in M$, 
	\[
		\pi_{-\alpha}(F^{e_\alpha}_*s X^\gamma) =
		\begin{cases}
		  \phi(F^{e_\alpha}_*s) X^{(\gamma - \alpha)/p^{e_\alpha}} &\quad\text{if $(\gamma - \alpha)/p^{e_\alpha} \in M$,}\\
		  0 &\quad\text{otherwise.} 
		\end{cases}
	\]
	Then $\pi_{-\alpha}$ satisfies the following properties:
	\begin{enumerate}
		\item $\pi_{-\alpha}$ is $k[M]$-linear.
		\item $\pi_{-\alpha}(F^{e_\alpha}_*X^\alpha) = 1$.
		\item $\pi_{-\alpha}$ maps $F^{e_\alpha}_*k[\sigma \cap M]$ into $k[\sigma \cap M]$. Thus, $\pi_{-\alpha}$ induces a splitting of $\lambda_{X_\alpha, e_\alpha}$.
	\end{enumerate}
\end{lemma}

\begin{proof}
	Since $\lambda_{X_\alpha, e_\alpha}$ splits, the set 
	\[
	M_{\alpha, e_\alpha} = \{\beta \in M \colon p^{e_\alpha}\beta + \alpha	\in \sigma\}
	\]
	is contained in $\sigma \cap M$ by \autoref{split-F-regular-cones}. Then the Lemma was verified in the proof of \autoref{split-F-regular-cones}\ref{split-F-regular-cones.monoid.one-e} $\Rightarrow$ \autoref{split-F-regular-cones}\ref{split-F-regular-cones.split}.
\end{proof}

\subsection{Splitting a relative interior algebra element} 
We will now show that the splitting of a single $X^\alpha \in k[\sigma \cap M]$, for $\alpha \in \Int(\sigma) \cap M$, implies that $k[\sigma \cap M]$ is split-$F$-regular.

\begin{theorem}
	\label{thm:splitting-rel-int-implies-SFR}
	Let $\sigma$ be a convex cone in $M_{\mathbb R}$ such that $0 \in \sigma$. Suppose there exists an element $\alpha \in \Int(\sigma) \cap M$ and an integer $e_\alpha > 0$ such that the $k[\sigma \cap M]$-linear map 
	\begin{align*}	
		\lambda_{\alpha, e_\alpha} \colon k[\sigma \cap M] &\to F^{e_\alpha}_*k[\sigma \cap M]\\
		1 &\mapsto F^{e_\alpha}X^\alpha
	\end{align*}
	splits. Then $k[\sigma \cap M]$ is split-$F$-regular.
\end{theorem}

\begin{proof}
	We first want to show that we can assume without loss of generality that $\sigma$ is full-dimensional in $M_{\mathbb R}$, that is, we are in \autoref{setup}.

	Let $N$ be the subgroup of $M$ generated by $\sigma \cap M$. Let $N_{\mathbb R}$ be the subspace of $M_{\mathbb R}$ generated by $\sigma \cap M$. Then $N_{\mathbb R} =  N \otimes_{\mathbb Z} \mathbb{R}$, and, as shown in \autoref{some-remarks}\autoref{some-remarks.f}, the cone 
	\[
	\tau \coloneqq \sigma \cap N_{\mathbb R}	
	\]
	satisfies the properties that the linear span of $\tau$ equals $N_{\mathbb R}$ and 
	$
	\tau \cap N = \sigma \cap M.	
	$ 
	We claim that $\alpha  \in \Int(\tau) \cap N$ as well. We already have $\alpha \in \sigma \cap M \subseteq N$. Thus, it remains to show $\alpha \in \Int(\tau)$.
	Since $\Int(\sigma)$ is an open subset of $\mathbb{R}\sigma$ containing $\alpha$ and since $N_{\mathbb R} \subseteq \mathbb{R}\sigma$, it follows that $\Int(\sigma) \cap N_{\mathbb R}$ is an open subset of $N_{\mathbb R}$ that contains $\alpha$. Moreover, $\Int(\sigma) \cap N_{\mathbb R} \subseteq \sigma \cap N_{\mathbb R} = \tau$. But this means $\Int(\sigma) \cap N_{\mathbb R} \subseteq \Int(\tau)$, and so, $\alpha \in \Int(\tau)$. Thus, replacing $\sigma$ by $\tau$, $M$ by $N$ and $M_{\mathbb R}$ by $N_{\mathbb R}$, we may assume that we are in \autoref{setup}.
	
	

	We will next show that for any $\beta \in \sigma \cap M$, there exists an integer $e_\beta > 0$ such that the $k[\sigma \cap M]$-linear map 
	\begin{align*}
	\lambda_{\beta,{e_\beta}} \colon k[\sigma \cap M] &\to F^{e_\beta}_*k[\sigma \cap M]\\
	1 &\mapsto F^{e_\beta}_*X^\beta
	\end{align*}
	splits.
	
	Since we are in \autoref{setup}, $\Int(\sigma)$ is a non-empty open subset of $M_{\mathbb R}$ by \autoref{some-remarks}\autoref{some-remarks.d}. Upon fixing a metric on $M_{\mathbb R}$, it follows that there exists an integer $n \gg 0$ such that 
	\[
		\alpha - \beta/n \in \Int(\sigma), 
	\]
	or equivalently, that $n\alpha - \beta \in \Int(\sigma)$.
	Thus, $X^{\beta}, X^{n\alpha - \beta} \in k[\sigma \cap M]$. Fix a splitting 
	\begin{align*}
		\varphi \colon F^{e_\alpha}_*k[\sigma \cap M] &\to k[\sigma \cap M]\\
		F^{e_\alpha}_*X^\alpha &\mapsto 1
	\end{align*}
	of $\lambda_{\alpha,e_\alpha}$. Then by 
    composing $\varphi$ with itself and pre-multiplication by a monomial, 
    there exists an integer $e' > 0$ and a $k[\sigma \cap M]$-linear map $F^{e'}_*k[\sigma \cap M] \to k[\sigma \cap M]$ such that $F^{e'}_*X^{n\alpha} = F^{e'}_*(X^\alpha)^n =  1$. Since $X^{n\alpha} = X^\beta \cdot X^{n\alpha - \beta}$, for the same $e'$, the $k[\sigma \cap M]$-linear map 
	\begin{align*}
		\lambda_{\beta, e'} \colon k[\sigma \cap M] &\to F^{e'}_*k[\sigma \cap M]\\
		1 &\mapsto F^{e'}_*X^\beta
	\end{align*}
	splits as well by pre-multiplying the splitting $F^{e'}_*X^{n\alpha} \mapsto 1$ by $F^e_* X^{n\alpha - \beta}$. 

	Now suppose that one has an arbitrary nonzero element 
	\begin{equation}\label{eq:arbitrary-element}
	Y \coloneqq \sum_{i=1}^n s_iX^{\alpha_i}	
	\end{equation}
	of $k[\sigma \cap M]$, where for all $i = 1, \dots, n$, the $\alpha_i \in \sigma \cap M$ are distinct and $s_i \in k \setminus \{0\}$. Since $M$ is free of finite rank, we can choose an integer $e'' \gg 0$ such that for all $i \neq 1$,
	\begin{equation}\label{eq:different-classes}
	\alpha_i - \alpha_1 \notin p^{e''}M,
	\end{equation}
	and also such that the $k[\sigma \cap M]$-linear map 
	\begin{align*}
	\lambda_{\alpha_1, e''} \colon k[\sigma \cap M] &\to F^{e''}_*k[\sigma \cap M]\\ 
	 1 &\mapsto F^{e''}_*X^{\alpha_1}.	
	\end{align*}
	splits.	

	In fact, by \autoref{lem:split-implies-torus-compatible-split}, upon fixing a splitting $\phi \colon F^{e''}_*k \to k$ of the $e$-th iterate of Frobenius on $k$, a splitting of $\lambda_{\alpha_1, e''}$ is given by the map
	\[
	\pi_{-\alpha_1} \colon F^{e''}_*k[\sigma \cap M] \to k[\sigma \cap M]
	\]
	with the property that for $s \in k$ and $\gamma \in \sigma \cap M$, 
	\[
	\pi_{-\alpha_1}(F^{e''}_*s X^{\gamma}) =
	\begin{cases}
	  \phi(F^{e''}_*s) X^{(\gamma - \alpha_1)/p^{e''}} &\quad\text{if $(\gamma - \alpha_1)/p^{e''} \in M$,}\\
	  0 &\quad\text{otherwise.} 
	\end{cases}
\]
Then by \autoref{eq:arbitrary-element} and \autoref{eq:different-classes}, the composition 
\[
F^{e''}_*k[\sigma \cap M] \xrightarrow{F^{e''}_*s^{-1}_1\cdot} F^{e''}_*k[\sigma \cap M] \xrightarrow{\pi_{-\alpha_1}} k[\sigma \cap M]	
\]
sends
\[
F^{e''}_*Y \mapsto F^{e''}_*s^{-1}_1Y \mapsto \sum_{i=1}^n\pi_{-\alpha_1}(F^{e''}_*s_is^{-1}_1X^{\alpha_i}) = \pi_{-\alpha_1}(F^{e''}_*s_1s^{-1}_1X^{\alpha_1}) = 1.
\]
Since $Y$ is an arbitrary nonzero element of $k[\sigma \cap M]$, we get $k[\sigma \cap M]$ is split-$F$-regular.
\end{proof}

\autoref{thm:splitting-rel-int-implies-SFR} along with results from \autoref{sec:Splitting versus purity and independence of ground field} allow us to show that there is no difference between split-$F$-regularity and $F$-pure-regularity for monoid algebras determined by convex cones. Moreover, the aforementioned properties are independent of the ground field as long as we are working in a fixed positive characteristic. The restriction on characteristic will later be removed in \autoref{cor:split-F-regular-independence-field-characteristic} as a consequence of finite generation.

\begin{corollary}
	\label{cor:splitF-regular-equiv-Fpure-regular}
	Let $\sigma$ be a convex cone in $M_{\mathbb R}$ such that $0 \in \sigma$. 
	Let $p > 0$ be a prime number. Then the following are equivalent:
	\begin{enumerate}[label=(\arabic*)]
		\item There exists a field $k$ of characteristic $p$ such that $k[\sigma \cap M]$ is split-$F$-regular. \label{cor:splitF-regular-equiv-Fpure-regular.split-one-field}
		\item For all fields $K$ of characteristic $p$, $K[\sigma \cap M]$ is split-$F$-regular.\label{cor:splitF-regular-equiv-Fpure-regular.split-any-field}
		\item For all fields $K$ of characteristic $p$, $K[\sigma \cap M]$ is $F$-pure-regular. \label{cor:splitF-regular-equiv-Fpure-regular.pure-any-field}
		\item There exists a field $k$ of characteristic $p$ such that $k[\sigma \cap M]$ is $F$-pure-regular. \label{cor:splitF-regular-equiv-Fpure-regular.pure-one-field}
		\item There exists a field $k$ of characteristic $p$ such for all $\beta \in \sigma \cap M$, there exists an integer $e_\beta > 0$ such that the unique $k[\sigma \cap M]$-linear map 
		\begin{align*}
			\lambda_{\beta,e_\beta} \colon k[\sigma \cap M] &\to F^{e_\beta}_*k[\sigma \cap M]\\
			1 &\mapsto F^{e_\beta}_*X^\beta
		\end{align*}
		is pure.\label{cor:splitF-regular-equiv-Fpure-regular.pure-one-field-lattice-points}
        \item For all $\beta\in \sigma \cap M$, there exists an integer $e_\beta > 0$ such that 
        \[M_{\beta, e_\beta} \coloneqq \{\gamma \in M: p^{e_\beta}\gamma + \beta \in \sigma\}\]
		is contained in $\sigma \cap M$.\label{cor:splitF-regular-equiv-Fpure-regular.pure-lattice-points-just-monoid}
	\end{enumerate}
\end{corollary}

\begin{proof}
	\ref{cor:splitF-regular-equiv-Fpure-regular.split-one-field} $\Rightarrow$ \ref{cor:splitF-regular-equiv-Fpure-regular.split-any-field}: By assumption that exists an integer $e_\alpha > 0$ such that the unique $k[\sigma \cap M]$-linear map
	$
        k[\sigma \cap M] \xrightarrow{1 \mapsto F^{e_\alpha}_*X^\alpha} F^{e_\alpha}_*k[\sigma \cap M]
    $
	splits. Then by \autoref{rem:independence-field}, for any field $K$ of characteristic $p$, the unique $K[\sigma \cap M]$-linear map 
	\begin{align*}
		K[\sigma \cap M] &\to F^{e_\alpha}_*K[\sigma \cap M]\\
		1 &\mapsto F^{e_\alpha}_*X^\alpha	
	\end{align*}
	splits as well. Upon applying \autoref{thm:splitting-rel-int-implies-SFR} to $K[\sigma \cap M]$, we conclude that this ring is split-$F$-regular.

	The implications \ref{cor:splitF-regular-equiv-Fpure-regular.split-any-field} $\Rightarrow$ \ref{cor:splitF-regular-equiv-Fpure-regular.pure-any-field} $\Rightarrow$ \ref{cor:splitF-regular-equiv-Fpure-regular.pure-one-field} $\Rightarrow$ \ref{cor:splitF-regular-equiv-Fpure-regular.pure-one-field-lattice-points} are clear.

	\ref{cor:splitF-regular-equiv-Fpure-regular.pure-one-field-lattice-points} $\Rightarrow$ \ref{cor:splitF-regular-equiv-Fpure-regular.pure-lattice-points-just-monoid} follows by \autoref{split-F-regular-cones}.

	\ref{cor:splitF-regular-equiv-Fpure-regular.pure-lattice-points-just-monoid} $\Rightarrow$ \ref{cor:splitF-regular-equiv-Fpure-regular.split-one-field}: By \ref{cor:splitF-regular-equiv-Fpure-regular.pure-lattice-points-just-monoid} and \autoref{split-F-regular-cones}, we get that for any field $k$ of characteristic $p > 0$ and for all $\beta \in \sigma \cap M$, the $k[\sigma \cap M]$-linear map
	\begin{align*}
		\lambda_{\beta, e_\beta} \colon k[\sigma \cap M] &\to F^{e_\beta}_*k[\sigma \cap M]\\
		1 &\mapsto F^{e_\beta}_*X^\beta
	\end{align*}
	splits. Thus, we would be done by \autoref{thm:splitting-rel-int-implies-SFR} if $\Int(\sigma) \cap M \neq \emptyset$.
	
	We now show that we can always reduce to the case where $\Int(\sigma) \cap M \neq \emptyset$ without altering the monoid $\sigma \cap M$. Indeed, we may assume $\sigma \cap M$ has a nonzero element, as otherwise there is nothing to prove. Next, by \autoref{some-remarks}\autoref{some-remarks.f}, we may assume that $\sigma$ is full-dimensional in $M_{\mathbb R}$. Since $\sigma$ has a nonzero element, we have that $\Int(\sigma)$ is a non-empty open set in $M_{\mathbb R}$ by \autoref{relative-interior-nonempty}\ref{relative-interior-nonempty.1} and \autoref{rel-int-open-span}. Then by \autoref{some-remarks}\autoref{some-remarks.e}, $\Int(\sigma) \cap M \neq \emptyset$. 
\end{proof}

\subsection{Split-{$F$}-regularity implies rationality of lineality space and extremal rays} As before, we let $\overline{\sigma}$ denote the topological closure of $\sigma$ in $M_{\mathbb R}$. The next result says that in the analysis of split-$F$-regularity of $k[\sigma \cap M]$ we can reduce to the case where $\sigma$ is closed.

\begin{lemma}
\label{rational-rays-closure}
Let $\sigma$ be a convex cone in $M_{\mathbb R}$ such that $0 \in \sigma$ and suppose $\Int(\sigma) \cap M \neq \emptyset$.  If $k[\sigma \cap M]$ is split-$F$-regular, then every rational ray of $\overline\sigma$ is contained in $\sigma$. Hence, $k[\overline\sigma \cap M] = k[\sigma \cap M]$.
\end{lemma}

\begin{proof}
Let $\alpha \in \Int(\sigma) \cap M$ and pick an integer $e_\alpha > 0$ such that the $k[\sigma\cap M]$-linear map 
$	
    k[\sigma \cap M] \xrightarrow{1 \mapsto F^{e_\alpha}_*X^\alpha} F^{e_\alpha}_*k[\sigma \cap M]
$
splits. 
Since a rational ray passes through a nonzero point in the lattice $M$, it suffices to show that
$\overline{\sigma} \cap M \subseteq \sigma \cap M.$
Let $\beta \in \overline{\sigma} \cap M$ be a nonzero element. 
Since $\alpha \in \Int(\sigma)$, by \autoref{rel-int-addition}, for all integers $f > 0$, 
$$p^f\beta + \alpha \in \Int(\sigma) \cap M \subseteq \sigma \cap M,$$
and hence $\beta \in \sigma \cap M$ by \autoref{split-F-regular-cones}\ref{split-F-regular-cones.pure}$\Rightarrow$\ref{split-F-regular-cones.monoid.one-e} applied to $f = e_\alpha$.
\end{proof}

The next result is the first key technical analysis of how certain rays in the boundary of $\sigma$ behave when $k[\sigma \cap M]$ is split-$F$-regular.

\begin{proposition}
\label{prop:common-setup-rationality-boundary}
Under \autoref{setup} suppose $k[\sigma \cap M]$ is split-$F$-regular. Let $\tau = \mathbb{R}_{\geq 0}\cdot x$ be a nontrivial ray in $\overline{\sigma}$. Let $W \subseteq M_{\mathbb R}$ be the linear subspace such that all the following conditions are satisfied:
\begin{enumerate}
	\item $W$ is defined over $M_{\mathbb Q}$,\label{prop:common-setup-rationality-boundary.a}
	\item $\tau \subseteq W$,\label{prop:common-setup-rationality-boundary.b}
	\item no proper subspace of $W$ defined over $M_{\mathbb Q}$ contains $\tau$, and \label{prop:common-setup-rationality-boundary.c}
	\item there exists a nonzero $\varphi \in (\sigma \cap W)^\vee \subseteq W^*$, such that $\tau \subseteq \varphi^{\perp}$. \label{prop:common-setup-rationality-boundary.d}
\end{enumerate}	
Then $\tau$ is rational. Consequently, $\tau \in \sigma$.
\end{proposition}

\begin{proof}
	Since $\sigma$ is full dimensional by hypothesis, $\Int(\sigma) \cap M \neq \emptyset$.  Thus, if $\tau \in \overline{\sigma}$ is rational then $\tau \in \sigma$ by \autoref{rational-rays-closure}.

Fix $\alpha \in \Int(\sigma) \cap M$. Assuming that $\tau$ is not rational, we will show that for any $e \in \mathbb{Z}_{> 0}$, there exists $\beta_e \in M - \overline{\sigma}$ such that $p^e\beta_e + \alpha \in \Int(\sigma) \cap M$. This will contradict \autoref{split-F-regular-cones} because it shows that for all integers $e > 0$, the set 
\[
M_{\alpha, e} \coloneqq \{\beta \in M \colon p^e\beta + \alpha \in \sigma \cap M\}	
\]
is not contained in $\sigma$.

So suppose that $\tau$ is not rational, that is,
$\tau \cap M_{\mathbb Q} = \{0\}.$
Since $W$ is defined over $M_{\mathbb Q}$ it follows that $\mathbb{R}\tau \subsetneq W$. Furthermore, if $W_{\mathbb Q} \coloneqq W \cap M_{\mathbb Q}$ is the $\mathbb Q$-vector subspace of rational vectors in $W$, then $W = W_{\mathbb Q} \otimes_{\mathbb Q} \mathbb R$.

For any $e \in \mathbb{Z}_{> 0}$,
$\Lambda_e \coloneqq p^eM \cap W$
is a lattice in $W$ generated by rational vectors: indeed, $W_{\mathbb Q}/\Lambda_e$ is a torsion group since it is a subgroup of the torsion group $M_\mathbb{Q}/p^eM$. So, $W_{\mathbb{Q}}/\Lambda_e \otimes_{\mathbb Z} \mathbb{R} = (W_{\mathbb{Q}}/\Lambda_e \otimes_{\mathbb Z} \mathbb{Q}) \otimes_{\mathbb Q} \mathbb{R} = 0$, that is, $\Lambda_e \otimes_{\mathbb Z} \mathbb{R} = W_{\mathbb Q} \otimes_{\mathbb Z} \mathbb{R} = W_{\mathbb Q} \otimes_{\mathbb Q} \mathbb{R} = W$. 

Since the ray 
$
-\tau \coloneqq \{-x: x \in \tau\}
$
is also not contained in any proper $\mathbb R$-vector subspace of $W$ defined over the rationals, \autoref{BCHM-density-consequence} implies that for all $e \in \mathbb{Z}_{> 0}$,
\begin{equation}
	\label{eq:dense}
	-\tau + \Lambda_e = \{y + \lambda: y \in -\tau, \lambda \in \Lambda_e\}
\end{equation}
is dense in $W$.

Let
\begin{equation}
\label{key-open-set}
U_\alpha \coloneqq -\alpha + \Int(\sigma).
\end{equation}
Openness of $\Int(\sigma)$ in $M_{\mathbb R}$ (see \autoref{some-remarks}\autoref{some-remarks.d}) shows that $U_\alpha$ is also an open subset of $M_{\mathbb R}$. Moreover, 
$
\tau \subseteq U_\alpha 
$
because $\tau + \alpha \subseteq \Int(\sigma)$ by \autoref{rel-int-addition}.  Thus, $U_\alpha$ is an open subset of $M_\mathbb{R}$ containing $\tau$ such that $U_\alpha + \alpha = \Int(\sigma)$. 

Note $\dim_{\mathbb R} W \geq 2$ because $\mathbb{R}\tau \subsetneq W$. 

\begin{claim}
\label{key-claim}
$U_\alpha \cap (W - \Supp(\varphi))$ is a non-empty open subset of $W$.
\end{claim}

\begin{proof}[Proof of \autoref{key-claim}]
The complement $W - \Supp(\varphi)$ is non-empty because $\varphi \neq 0$ in $W^*$. The set $U_\alpha \cap (W - \Supp(\varphi)) = (U_\alpha \cap W) \cap (W - \Supp(\varphi))$ is open in $W$ as it is the intersection of two opens. By our choice of $\varphi$, $\tau$ is contained in the topological boundary, $\varphi^{\perp}$, of $\Supp(\varphi)$. Since $\tau \subseteq U_\alpha \cap W$, the open set $U_\alpha \cap W$ intersects the topological boundary of $\Supp(\varphi)$ non-trivially. Consequently, $U_\alpha \cap W$ contains points that are \emph{not} in $\Supp(\varphi)$, that is, $U_\alpha \cap (W - \Supp(\varphi)) = (U_\alpha \cap W) \cap (W - \Supp(\varphi)) \neq \emptyset$, proving \autoref{key-claim}.\phantom\qedhere
\end{proof}

A dense set always intersects a non-empty open set non-trivially and recall that for all integers $e > 0$, $-\tau + \Lambda_e$ is dense in $W$ (see \autoref{eq:dense}). Therefore, for all integers $e > 0$, there exists 
\begin{equation}
	\label{eq:z_e-equation}
z_e \in [U_\alpha \cap (W - \Supp(\varphi))] \cap [-\tau + \Lambda_e].
\end{equation}
Choose $x_e \in \tau$ and $\beta_e \in M$ such that 
$p^e\beta_e \in \Lambda_e \subset W$ and $z_e = -x_e + p^e\beta_e$.

\begin{claim}
\label{claim-two}
$\beta_e \notin \sigma$.
\end{claim}

\begin{proof}[Proof of \autoref{claim-two}] Recall that by hypothesis, $\varphi \in (\sigma \cap W)^\vee$ and $\tau \subseteq \varphi^\perp$. Thus, $-x_e \in -\tau \subseteq \varphi^\perp$ as well. If $\beta_e \in \sigma$, then $p^e\beta_e \in \sigma \cap W$. But this would imply
\begin{equation}
\label{chain-equalities}
\varphi(z_e) = \varphi(-x_e + p^e\beta_e) =  \varphi(p^e\beta_e) \geq 0,
\end{equation}
contradicting $z_e$ being an element of $W - \Supp(\varphi)$. This proves \autoref{claim-two}. \phantom\qedhere
\end{proof}

By definition of $U_\alpha$ (see \autoref{key-open-set}) and the fact that $z_e \in U_\alpha$ (see \autoref{eq:z_e-equation}), for all $e \in \mathbb{Z}_{>0}$, we see that
$
z_e + \alpha \in \Int(\sigma).
$ 
Now because $x_e \in \tau \subset \overline{\sigma}$, by \autoref{rel-int-addition} we have
\[
\alpha + p^e\beta_e = \alpha + (z_e + x_e) = x_e + (z_e + \alpha) \in \Int(\sigma) \subseteq \sigma.
\]
Furthermore, since $\alpha, p^e\beta_e \in M$, we have 
$
\alpha + p^e\beta_e \in \sigma \cap M.
$
Then \autoref{claim-two} shows that for all  $e \in \mathbb{Z}_{> 0}$, there exists $\beta_e \in M - \sigma$ such that $\alpha + p^e\beta_e \in \sigma$. In other words, for all $e \in \mathbb{Z}_{> 0}$, the set
\[
M_{\alpha, e} \coloneqq \{\beta \in M \colon p^e\beta + \alpha \in \sigma\}	
\]
is not contained in $\sigma$. This is impossible by \autoref{split-F-regular-cones} since $k[\sigma \cap M]$ is split-$F$-regular. Hence $\tau$ must be rational.
\end{proof}

\autoref{prop:common-setup-rationality-boundary} has the following important consequence.

\begin{theorem}
	\label{extremal-rays-are-rational}
	Under \autoref{setup} suppose $k[\sigma \cap M]$ is split-$F$-regular.
Then we have the following:
\begin{enumerate}[label=(\arabic*)]
	\item The lineality space $\sigma \cap -\sigma$ of $\sigma$ is defined over $M_{\mathbb Q}$.
    \label{extremal-rays-are-rational.lineality}
	\item The lineality space of $\sigma$ equals the lineality space of $\overline{\sigma}$.
	\label{extremal-rays-are-rational.lineality-closure}
	\item If $\overline{\sigma}$ is strongly convex, then every extremal ray of $\overline{\sigma}$ is rational and $\sigma  = \overline{\sigma}$.
    \label{extremal-rays-are-rational.extremalRay}
\end{enumerate}
\end{theorem}

\begin{proof}
	\ref{extremal-rays-are-rational.lineality}  Let $L \coloneqq \sigma \cap -\sigma$ and assume for contradiction that $L$ is not defined over $M_{\mathbb Q}$. Let $W$ be the smallest subspace of $M_{\mathbb R}$ such that $W$ is defined over $M_{\mathbb Q}$ and
	$
	L \subseteq W.	
	$
	 Since $L$ is not rational, we have $\dim_{\mathbb R} W \geq 2$ and $\dim_{\mathbb R} L \geq 1$. Furthermore, $L \neq W$ (as otherwise $L$ would be defined over $M_{\mathbb Q}$).  Thus, by \autoref{irrational-lines}\ref{irrational-lines.5}, there exists $x \in L$ such that $W$ is the smallest linear subspace of $M_{\mathbb R}$ defined over $M_{\mathbb Q}$ that contains $x$.

	Now let 
	$
	\tau \coloneqq \mathbb{R}_{\geq 0}x
	$ 
	be the ray of $\sigma$ generated by $x$. The rational subspace $W$ satisfies conditions \autoref{prop:common-setup-rationality-boundary.a}, \autoref{prop:common-setup-rationality-boundary.b} and \autoref{prop:common-setup-rationality-boundary.c} of \autoref{prop:common-setup-rationality-boundary} for the ray $\tau$. We claim that condition \autoref{prop:common-setup-rationality-boundary.d} of \autoref{prop:common-setup-rationality-boundary} is also satisfied, that is, there exists a nonzero $\varphi \in (\sigma \cap W)^\vee \subseteq W^*$ such that $\tau \subseteq \varphi^\perp$. Indeed, we first observe that 
	$
	(\sigma \cap W)^\vee \neq \{0\}.
	$ 
	Otherwise $\sigma \cap W = W$ by \autoref{cor:dense-convex-cone}, which would imply that the lineality space $L$ of $\sigma$ contains $W$. But $L \subsetneq W$, so we get a contradiction. 

	Note that $L$ is also the lineality space of $\sigma \cap W$. Thus, for all $\varphi \in (\sigma \cap W)^\vee$, $L \subseteq \varphi^\perp$. Since $\tau \subseteq L$, choosing \emph{any} nonzero $\varphi \in (\sigma \cap W)^\vee$, we have $\tau \subseteq \varphi^\perp$. Now we are in the situation of \autoref{prop:common-setup-rationality-boundary}, and we can conclude that $\tau$ is rational. But then $\mathbb{R}\tau$ is the smallest subspace of $M_{\mathbb R}$ defined over $M_{\mathbb Q}$ that contains $\tau$, which implies $\mathbb{R}\tau = W$ by our choice of $x$ (and hence $\tau$). This is impossible since $\dim_{\mathbb R}(\mathbb{R}\tau) = 1$ while $\dim_{\mathbb{R}}(W) \geq 2$. By contradiction, $L$ must be defined over $M_{\mathbb Q}$, which finishes the proof of \ref{extremal-rays-are-rational.lineality}.

	\ref{extremal-rays-are-rational.lineality-closure}  
    It suffices to show that $\overline{\sigma} \cap -\overline{\sigma} \subseteq \sigma$. 
      Fix $\alpha \in \Int \sigma \cap M$ and choose an integer $e_\alpha > 0$ such that $k[\sigma \cap M] \xrightarrow{1 \mapsto F^{e_\alpha}_*X^\alpha} F^{e_\alpha}_*k[\sigma\cap M]$ splits. By \autoref{rational-rays-closure}, $\sigma \cap M = \overline{\sigma} \cap M$, and hence, $k[\sigma \cap M] = k[\overline{\sigma} \cap M]$.
    Since $\Int(\sigma) = \Int(\overline{\sigma})$ by \autoref{relative-interior-nonempty}\ref{relative-interior-nonempty.3}, we get $\alpha \in \Int(\overline{\sigma}) \cap M$. Thus\begin{align*}
		k[\overline{\sigma} \cap M] &\to F^{e_\alpha}_*k[\overline{\sigma} \cap M]\\
		1 &\mapsto F^{e_\alpha}_*X^\alpha
	\end{align*}	
	splits as well, so by \ref{extremal-rays-are-rational.lineality} applied to $\overline{\sigma}$ we get that $\overline{\sigma} \cap - \overline{\sigma}$ is defined over $M_{\mathbb Q}$. Let $\mathcal{B} \subseteq M$ be a basis of $\overline{\sigma} \cap - \overline{\sigma}$. Then $\mathcal{B} \cup -\mathcal{B}$ generates $\overline{\sigma} \cap - \overline{\sigma}$ as a convex cone and $\mathcal{B} \cup -\mathcal{B} \subseteq \overline{\sigma} \cap M = \sigma \cap M$. Thus, $\overline{\sigma} \cap - \overline{\sigma} \subseteq \sigma$.

	\ref{extremal-rays-are-rational.extremalRay} Let $\tau$ be an extremal ray of $\overline{\sigma}$. Again, let $W$ be the smallest subspace of $M_\mathbb{R}$ such that $W$ is defined over $M_{\mathbb Q}$ and $\tau \subseteq W$. 
	
	If $\dim_{\mathbb R}(W) = 1$, then $W = \mathbb{R}\tau$. So, $\tau$ is rational since $W$ is.

	Suppose $\dim_{\mathbb R}(W) \geq 2$. Note that $\tau$ is also an extremal ray of the closed strongly convex cone $\overline{\sigma} \cap W$. If the $\mathbb{R}$-linear span of $\overline{\sigma} \cap W$ has dimension $\geq 2$, then $\tau$ is contained in a proper face of $\overline{\sigma} \cap W$ by \autoref{extremal-ray-properties}\ref{extremal-ray-properties.dimlinealitygeq2}, and so, we get a nonzero $\varphi \in (\overline{\sigma} \cap W)^\vee$ such that $\tau \subseteq \varphi^\perp$. If the $\mathbb{R}$-linear span of $\overline{\sigma} \cap W$ has dimension $1$, then this span, and hence also $\overline{\sigma} \cap W$ is clearly contained in some hyperplane of $W$ since $\dim_{\mathbb R}(W) \geq 2$. This again gives us a nonzero $\varphi \in W^*$ such that $\tau \subseteq \overline{\sigma} \cap W \subseteq \varphi^\perp$. Thus, in either situation there exists 
	\[
	0 \neq \varphi \in (\overline{\sigma} \cap W)^\vee 	
	\]
	such that $\tau \subseteq \varphi^\perp$.
	Since $\sigma \cap W \subseteq \overline{\sigma} \cap W$, we have 
	$
	(\overline{\sigma} \cap W)^\vee \subseteq (\sigma \cap W)^\vee	
	$
	in $W^*$. Thus, $\varphi \in (\sigma \cap W)^\vee$ and we are again in the situation of \autoref{prop:common-setup-rationality-boundary}, which allows us to conclude that $\tau$ is rational. Hence, every extremal ray of $\overline{\sigma}$ is rational. Since every rational ray of $\overline{\sigma}$ must be contained in $\sigma$ by \autoref{rational-rays-closure}, we get that every extremal ray of $\overline{\sigma}$ is also a ray of $\sigma$. But $\overline{\sigma}$ is generated as a convex cone by its extremal rays by \autoref{extremal-strongly-convex}, and so, $\sigma = \overline{\sigma}$.
\end{proof}

Our next goal is to show that if $k[\sigma \cap M]$ is split-$F$-regular, then $\sigma$ is closed and is generated as a convex cone by rational rays. For this we will need the following lemma.

\begin{lemma}
	\label{lem:cone-ses}
	Under \autoref{setup} let $L$ be the lineality space of $\sigma$ and assume that $L$ is defined over $M_{\mathbb Q}$. Consider the short exact sequence of vector spaces
	\[
	0 \to L \to M_{\mathbb R} \xrightarrow{\pi} M_{\mathbb R}/L \to 0.	
	\]
	and the cone $\pi(\sigma)$ of $M_{\mathbb R}/L$. We have the following:
	\begin{enumerate}[label=(\arabic*)]
		\item Suppose $\pi(\sigma)$ is generated as a convex cone by (finitely many) elements of $\pi(\sigma) \cap \pi(M)$. Then $\sigma$ is generated as a convex cone by (finitely many) elements of $\sigma \cap M$. \label{lem:cone-ses.generators}
		\smallskip
		\item $k[\pi(\sigma) \cap \pi(M)]$ is a direct summand of $k[\sigma \cap M]$. Thus, if $k[\sigma \cap M]$ is split-$F$-regular, then $k[\pi(\sigma) \cap \pi(M)]$ is split-$F$-regular.\label{lem:cone-ses.splitF-regular}
	\end{enumerate}
\end{lemma}

\begin{proof}
	(1) Since $L$ is defined over $M_{\mathbb Q}$, after scaling one can assume that $L$ has a basis $\{\ell_1,\dots,\ell_d\}$ such that each $\ell_i \in M$. Since $L$ is the lineality of $\sigma$ note that we have
	\[
	\{\ell_1,\dots,\ell_d\} \cup \{-\ell_1,\dots,-\ell_d\} \subseteq \sigma \cap M.	
	\]
	We now claim that
	\begin{equation}
		\label{eq:intersection}
	\pi(\sigma) \cap \pi(M) = \pi(\sigma \cap M).	
	\end{equation}
	The non-trivial inclusion is to show that $\pi(\sigma) \cap \pi(M) \subseteq \pi(\sigma \cap M)$. So let $y \in \pi(\sigma) \cap \pi(M)$. Then there exist $x \in \sigma$ and $\alpha \in M$ such that
	$
	\alpha + L = y = x + L .	
	$
	Then $\alpha = x + \ell$ for some $\ell \in L$. But $x, \ell \in \sigma$ since $L \subseteq \sigma$. Thus, $\alpha = x + \ell \in \sigma$ as well, and hence, $\alpha \in \sigma \cap M$. This shows that $y = \alpha + L \in \pi(\sigma \cap M)$, proving the desired inclusion.

	Thus, if $\{y_i\}_{i \in I} \subseteq \pi(\sigma) \cap \pi(M) = \pi(\sigma \cap M)$ is a generating set of $\pi(\sigma)$ as a convex cone, choosing $x_i \in \sigma \cap M$ such that 
	$
		y_i = x_i + L,
	$
	one can readily verify that $\sigma$ is generated as a convex cone by the set
	\[
	\{x_i\}_{i \in I} \cup \{\ell_1,\dots,\ell_d\} \cup \{-\ell_1,\dots,-\ell_d\} \subseteq \sigma \cap M.
	\]
	Note that if $I$ is a finite set then the above generating set for $\sigma$ is also finite.

	(2) The second assertion follows from the first by \autoref{lem:basic-facts-SFR}\ref{lem:basic-facts-SFR.6}.

	Since $L$ is defined over $M_{\mathbb Q}$, it follows that $(L \cap M) \otimes_{\mathbb Z} \mathbb{R} = L$. Thus, applying $\_ \otimes_{\mathbb Z} \mathbb{R}$ to the short exact sequence of Abelian groups
	\begin{equation}
		\label{eq:ses-abelian-gp}
	0 \to (L \cap M) \to M \xrightarrow{\pi} \pi(M) \to 0	
	\end{equation}
	yields the short exact sequence of $\mathbb R$-vector spaces 
	\[
	0 \to L \to M_{\mathbb R} \xrightarrow{\pi} M_{\mathbb R}/L \to 0.	
	\]
	Note that $\pi(M)$ is a finitely generated torsion-free Abelian group, and is hence free. Therefore \autoref{eq:ses-abelian-gp} splits. Let $i \colon \pi(M) \to M$ be a splitting of $\pi \colon M \twoheadrightarrow \pi(M)$, that is, $\pi \circ i = \id_{\pi(M)}$. Since $i$ is an additive map, we get an induced ring homomorphism
	\begin{align*}
		k[\pi(M)] &\xrightarrow{i^*} k[M]\\
		X^{\gamma} &\mapsto X^{i(\gamma)}	
	\end{align*}
	which is the identity on $k$.
	Similarly, the additive map $\pi \colon M \to \pi(M)$ induces a ring homomorphism
	\begin{align*}
		k[M] &\xrightarrow{\pi^*} k[\pi(M)]\\
		X^{\alpha} &\mapsto X^{\pi(\alpha)}.
	\end{align*}
	By construction, $\pi^* \circ i^* = \id_{k[\pi(M)]}$ (thus $i^*$ splits not only as a map of $k[\pi(M)]$-modules but also as a map of $k[\pi(M)]$-algebras). Additionally, $\pi^*$ restricted to $k[\sigma \cap M]$ maps into $k[\pi(\sigma \cap M)] \stackrel{\autoref{eq:intersection}}{=} k[\pi(\sigma) \cap \pi(M)]$.
	 Thus, in order to show that $k[\pi(\sigma) \cap \pi(M)]$ is a direct summand of $k[\sigma \cap M]$, it is enough to show that the restriction of $i^*$ to the subring $k[\pi(\sigma) \cap \pi(M)]$ of $k[\pi(M)]$ maps into $k[\sigma \cap M]$.
	So suppose $y \in \pi(\sigma) \cap \pi(M)$. It suffices to show that
	\[
	i^*(X^y) \in k[\sigma \cap M],	
	\]
	or equivalently, that $i(y) \in \sigma \cap M$.	
	We already have 
	$
	i(y) \in M
	$
	because $i^*$ maps into $k[M]$. Also,
	\[
	i(y) + L = \pi(i(y)) = y \in \pi(\sigma) \cap \pi(M) \stackrel{\autoref{eq:intersection}}{=} \pi(\sigma \cap M).
	\]
	Thus, there exists $x \in \sigma \cap M$ such that
	$
	i(y) + L = x + L.	
	$
	Then $i(y) = x + \ell$, for some $\ell \in L$. Since $L \subseteq \sigma$, we have $i(y) = x + \ell \in \sigma$. Hence, $i(y) \in \sigma \cap M$, as desired.
\end{proof}

{}

Combining the previous results we can prove the following characterization of split-$F$-regular cone algebras for arbitrary cones satisfying the hypotheses of \autoref{setup}.

\begin{theorem}
	\label{thm:split-F-regular-arbitrary-cone-closed}
	Under \autoref{setup} suppose $k[\sigma \cap M]$ is split-$F$-regular. Then $\overline{\sigma}$ is generated as a convex cone by rational rays and $\sigma = \overline{\sigma}$.
\end{theorem}

\begin{proof}
As before, if we can show that $\overline{\sigma}$ is generated as a convex cone by rational rays, then it will follow that $\overline{\sigma} = \sigma$ by \autoref{rational-rays-closure}. In fact, \autoref{rational-rays-closure} shows that
$
k[\sigma \cap M] = k[\overline{\sigma} \cap M].
$ 
Thus, applying \autoref{thm:splitting-rel-int-implies-SFR} we get that $k[\sigma \cap M] = k[\overline{\sigma} \cap M]$ is split-$F$-regular

Let $L$ be the lineality space $\overline{\sigma}$.
Then $L$ is defined over $M_{\mathbb Q}$ by \autoref{extremal-rays-are-rational}\ref{extremal-rays-are-rational.lineality} and \ref{extremal-rays-are-rational.lineality-closure}. We may assume $L \neq M_{\mathbb R}$ as otherwise $\overline{\sigma} = L$ and we are done since $L$ is defined over $M_{\mathbb Q}$.
Let
\[
\pi \colon M_{\mathbb R} \twoheadrightarrow M_{\mathbb R}/L	
\]
be the projection map. Note that $\pi(\overline{\sigma})$ is a closed strongly convex cone by \autoref{prop:lineality-quotient-strong-convex}. Furthermore, since $\overline{\sigma}$ is full-dimensional in $M_{\mathbb R}$, $\pi(\overline{\sigma})$ is full-dimensional in $M_{\mathbb R}/L$. In other words, $\pi(\overline{\sigma})$ satisfies the hypotheses of \autoref{setup}. We also have that $\pi(M)$ is a lattice of $M_{\mathbb R}/L$ and $\pi(M) \otimes_{\mathbb Z} \mathbb R = M_{\mathbb R}/L$. Since $k[\pi(\overline{\sigma}) \cap \pi(M)]$ is a direct summand of $k[\overline{\sigma} \cap M]$ by \autoref{lem:cone-ses}\ref{lem:cone-ses.splitF-regular}, we get $k[\pi(\overline{\sigma}) \cap \pi(M)]$ is split-$F$-regular. Hence we can apply \autoref{extremal-rays-are-rational}(2) to conclude that $\pi(\overline{\sigma})$ has rational extremal rays (here rationality is with respect to the lattice $\pi(M)$), that is, $\pi(\overline{\sigma})$ is generated as a convex cone by elements of $\pi(\overline{\sigma}) \cap \pi(M)$. Then $\overline{\sigma}$ is generated as a convex cone by $\overline{\sigma} \cap M$ (see \autoref{lem:cone-ses}\ref{lem:cone-ses.generators}). This completes the proof of \autoref{thm:split-F-regular-arbitrary-cone-closed}. 
\end{proof}

\subsection{Split-{$F$}-regularity and rationality of the dual cone} 
Our ultimate aim is to show that if $k[\sigma \cap M]$ is split-$F$-regular, then $\sigma$ is a closed cone generated by finitely many rational rays. This will imply that $\sigma^\vee$ is also generated by finitely many rational points/rays in the dual space by \cite[Thm.\ 16]{FenchelConvexConesSetsFunctionsBook}. The goal of this subsection is to establish the rationality of $\sigma^\vee$ without assuming a priori that $\sigma$ is generated by \emph{finitely many} rational rays. Our main result is:

\begin{theorem}
	\label{thm:dual-cone-rational-SFR}
	Under \autoref{setup} suppose that $\sigma$ is 
	strongly convex in $M_{\mathbb R}$. If $k[\sigma \cap M]$ is split-$F$-regular, then $\sigma^\vee$ is a strongly convex full-dimensional closed cone in $M_{\mathbb R}^*$ generated by rational extremal rays.
\end{theorem}

\begin{proof}
	Split-$F$-regularity of $k[\sigma \cap M]$ implies that $\sigma$ is a closed cone by \autoref{thm:split-F-regular-arbitrary-cone-closed}.
	We know that $\sigma^\vee$ is always a closed convex cone in $M_{\mathbb R}^*$ (\autoref{dual-cone-closed}). Since $\sigma$ is full-dimensional by the hypothesis of \autoref{setup}, it follows by \autoref{lem-characterize-strong-conv} that $\sigma^\vee$ is strongly convex. Thus, by \autoref{extremal-strongly-convex}, $\sigma^\vee$ is generated by its extremal rays.  By \autoref{closed-strongly-convex-cone}, $\sigma^\vee$ is also full-dimensional in $M_{\mathbb R}^*$  because $(\sigma^\vee)^\vee = \overline{\sigma} = \sigma$ is strongly convex. 
	
	Rationality of the extremal rays of $\sigma^\vee$ here means that if 
	 $
	 \{e_1, \dots, e_d\} \subset M
	 $
	 is a $\mathbb Z$-basis for $M$ and hence also an $\mathbb R$-basis of $M_{\mathbb R}$, and $M^* = \mathbb{Z}e_1^* \oplus \dots \oplus \mathbb{Z}e_d^*$ is the dual lattice in $M_{\mathbb R}^*$, then for any extremal ray $\tau$ of $\sigma^\vee$, $\tau \cap M^*$ has a nonzero element.

	 We proceed by induction on $d = \dim_{\mathbb{R}}(M_\mathbb{R})$. The $d = 1$ case is clear because $\sigma$ is then a rational ray and hence so is $\sigma^\vee$. So assume $d \geq 2$.

	 Observe that the extremal rays of $\sigma$ are rational by \autoref{extremal-rays-are-rational}(2) since $k[\sigma \cap M]$ is split-$F$-regular (here we are using that $\overline{\sigma} = \sigma$ and $\sigma$ is strongly convex). Now, let
	 $
		\tau \subseteq \sigma^\vee
	 $ 
	 be an extremal ray of $\sigma^\vee$. Since $\sigma^\vee$ is full-dimensional in a vector space of dimension $\geq 2$, it follows by \autoref{extremal-ray-properties}\ref{extremal-ray-properties.dimlinealitygeq2} that $\tau$ is contained in a proper face of $\sigma$, and hence, by \autoref{rel-int-as-comp} that 
	 \begin{equation}
		\label{eq:extremal-boundary}
		\tau \cap \Int(\sigma^\vee) = \emptyset.	
	 \end{equation}
	 Let $\{\mathbb{R}_{\geq 0}x_i \colon i \in I\}$ be the extremal rays of $\sigma$, and assume that for all $i \in I$, $x_i \in \sigma \cap M$. Then there exists $i \in I$ such that 
	 $ 
		\tau \subseteq x_i^\perp.
	 $ 
	 Otherwise, choosing \emph{any} nonzero $\varphi \in \tau$, we would get that $\varphi(x_i) > 0$ for all $i \in I$ (since $\sigma$ is contained in $\Supp(\varphi)$), which would consequently imply that for all nonzero $x \in \sigma$, $\varphi(x) > 0$ since $x$ is a $\mathbb{R}_{> 0}$-linear combination of some of the $x_i$'s. But then $\varphi$ would be an element of $\Int(\sigma^\vee)$ by \autoref{cor:Relint-dual-cone}, contradicting \autoref{eq:extremal-boundary}.

	 So choose $x_i \in \sigma \cap M$ to be an extremal point of $\sigma$ such that $\tau \subseteq x_i^\perp$ and set
	 $
		\rho_i = \mathbb{R}x_i
	 $
	 to be the span of $x_i$. Note that $\rho_i$ is defined over $M_{\mathbb Q}$ by construction. Let 
	 $
		\pi \colon M_{\mathbb R} \twoheadrightarrow M_{\mathbb R}/\rho_i
	 $
	 be the projection map. Since $\mathbb{R}_{\geq 0}x_i$ is extremal, \autoref{prop:image-strongly-convex} implies that
	 \(
		\pi(\sigma)
	 \)
	 is a strongly convex cone of $M_{\mathbb R}/\rho_i$. Since $\sigma$ is full-dimensional in $M_{\mathbb R}$ and contains $0$, we get that $\pi(\sigma)$ is also full-dimensional and contains $0$, that is, $\pi(\sigma)$ satisfies the hypotheses of \autoref{setup}. 
	 Furthermore, by \autoref{prop:dual-cones-quotient}, we know that $\pi(\sigma)^\vee$ can be identified with
	 \(
		\sigma^\vee \cap \rho_i^\perp = \sigma^\vee \cap x_i^\perp. 
	 \)
	 Note that by our choice of $x_i^\perp$, we have that $\tau$ is a ray of the cone
	 \(
		\sigma^\vee \cap x_i^\perp.
	 \)
	 Moreover, $\tau$ remains extremal in the subcone $\sigma^\vee \cap x_i^\perp$ of $\sigma^\vee$. Thus, if we can show that $k[\pi(\sigma) \cap \pi(M)]$ is split-$F$-regular then it will follow by induction (since $\dim_{\mathbb R}M_{\mathbb R}/\rho_i = d-1$) that all the extremal rays of $\sigma^\vee \cap x_i^\perp = \pi(\sigma)^\vee$ are rational. Then $\tau$ would be rational in $M_{\mathbb R}^*$ (here we are using that $x_i^\perp$ is rational because $x_i$ is a rational point), which is what we want to establish.

	Since $k[\sigma \cap M]$ is split-$F$-regular, the localization
	\(
	k[\sigma \cap M]\bigg{[}\frac{1}{X^{x_i}}\bigg{]}	
	\)
	is split-$F$-regular as well by \autoref{lem:basic-facts-SFR}\ref{lem:basic-facts-SFR.2}. Consider the cone 
	\(
	\sigma + \rho_i.
	\)
	\begin{claim}
		\label{claim:localization-cone}
	We claim that $k[\sigma \cap M]\big{[}\frac{1}{X^{x_i}}\big{]}	= k[(\sigma+\rho_i) \cap M]$.
	\end{claim}
	
	\begin{proof}[Proof of \autoref{claim:localization-cone}]
		The non-trivial inclusion is 
	\[
		k[(\sigma+\rho_i) \cap M] \subseteq k[\sigma \cap M]\bigg{[}\frac{1}{X^{x_i}}\bigg{]}.
	\]
	So let $\alpha \in (\sigma+\rho_i) \cap M$. The claim follows if we can show that $X^\alpha \in k[\sigma \cap M]\bigg{[}\frac{1}{X^{x_i}}\bigg{]}$. Let $\beta \in \sigma$ and $\lambda \in \mathbb{R}$ such that 
	\(
	\alpha = \beta + \lambda x_i.	
	\)
	Choose an integer $n \gg 0$ such that $n + \lambda > 0$. We then get 
	\(
	\alpha + nx_i = \beta + (n + \lambda)x_i \in \sigma \cap M,	
	\)
	since $\alpha, x_i \in M$ and $\beta, x_i \in \sigma$. Thus, 
	\[
	X^\alpha = {X^{\alpha + nx_i}}\frac{1}{(X^{x_i})^n} \in k[\sigma \cap M]\bigg{[}\frac{1}{X^{x_i}}\bigg{]},
	\]
	proving \autoref{claim:localization-cone}. \phantom\qedhere
	\end{proof}

	By our claim, we get that $k[(\sigma +\rho_i) \cap M]$ is split-$F$-regular. 
	Note that 
	\(
	\pi(\sigma + \rho_i) = \pi(\sigma).
	\)
	Furthermore, the lineality space of $\sigma + \rho_i$ equals $\rho_i$ by \autoref{lem:lineality-sum-extremal-subcone} because $\rho_i$ is the linear span of the extremal ray $\mathbb{R}_{\geq 0}\cdot x_i$. Then $k[\pi(\sigma)\cap \pi(M)] = k[\pi(\sigma + \rho_i) \cap \pi(M)]$ is also split-$F$-regular by \autoref{lem:cone-ses}\ref{lem:cone-ses.splitF-regular}. This completes the proof of the Theorem.
\end{proof}

\subsection{Split-{$F$}-regularity and affine hyperplane sections of cones} We continue to work in \autoref{setup}. Suppose $H$ is a hyperplane in $M_{\mathbb R}$ that is defined over $M_{\mathbb Q}$. Then 
$
H \cap M	
$
is a lattice of rank $d-1$. Indeed, $H \cap M$ is a subgroup of the free Abelian group $M$, and the latter has rank $d$. Thus, $H \cap M$ is also a free Abelian group of finite rank, say $d'$. Furthermore, any basis of $H$ consisting of elements in $M$ is contained in $H \cap M$. Thus, $H \cap M$ contains a free Abelian group of rank $d-1$, and so, $d' = d-1$. 

\begin{proposition}
	\label{prop:hyperplane-section-split-F-regular}
	Under \autoref{setup}, suppose $k[\sigma \cap M]$ is split-$F$-regular. Let $H$ be a hyperplane of $M_{\mathbb R}$ that is defined over $M_{\mathbb Q}$. Let $\sigma_H \coloneqq \sigma \cap H$ and $M_H$ be the lattice $H \cap M$. Then $k[\sigma_H \cap M_H]$ is also split-$F$-regular. 
\end{proposition}

\begin{proof}
	Split-$F$-regularity of $k[\sigma_H \cap M_H]$ will follow if we can show that $k[\sigma_H \cap M_H]$ is a direct summand of $k[\sigma \cap M]$ (\autoref{lem:basic-facts-SFR}\ref{lem:basic-facts-SFR.6}). We first claim that
	\begin{equation}
		\label{eq:hyperplane-cone-equality}
	\sigma_H \cap M_H = (\sigma \cap M) \cap \mathbb{Z}(\sigma_H \cap M_H).	
	\end{equation}
	It is clear that $\sigma_H \cap M_H \subseteq (\sigma \cap M) \cap \mathbb{Z}(\sigma_H \cap M_H)$. Note that $\sigma_H \cap M_H \subseteq H$, and so, the group $\mathbb{Z}(\sigma_H \cap M_H)$ generated by $\sigma_H \cap M_H$ is also contained in $H$. Thus,
	\[
	(\sigma \cap M) \cap \mathbb{Z}(\sigma_H \cap M_H) \subseteq (\sigma \cap M) \cap H = \sigma_H \cap M_H,
	\]
	proving \autoref{eq:hyperplane-cone-equality}. One can now define a $k$-linear map
	\[
	\varphi \colon k[\sigma \cap M] \to k[\sigma_H \cap M_H]	
	\]
	as follows: for $\alpha \in \sigma \cap M$ and $s \in k$,
	\[
	\varphi(sX^\alpha) =
	\begin{cases}
		\textrm{$sX^\alpha$ \hspace{4mm} if  $\alpha \in \sigma_H \cap M_H$},\\
		\textrm{$0$ \hspace{7.5mm}  if $\alpha \in (\sigma \cap M) \setminus (\sigma_H \cap M_H)$}. 
	\end{cases}	
	\]
	By definition, $\varphi(1) = 1$. We claim that $\varphi$ is $k[\sigma_H \cap M_H]$-linear. Suppose $\alpha \in \sigma \cap M$ and $\beta \in \sigma_H \cap M_H$. If $\alpha \in \sigma_H \cap M_H$, then $\beta + \alpha \in \sigma_H \cap M_H$. On the other hand, if $\alpha \in (\sigma \cap M) \setminus (\sigma_H \cap M_H)$, then $\beta + \alpha \in (\sigma \cap M) \setminus (\sigma_H \cap M_H)$ as well. Otherwise, $\alpha = (\beta + \alpha) - \beta \in \mathbb{Z}(\sigma_H \cap M_H)$, and so, $\alpha \in (\sigma \cap M) \cap \mathbb{Z}(\sigma_H \cap M_H) = \sigma_H \cap M_H$ by \autoref{eq:hyperplane-cone-equality}, which is a contradiction. The upshot of this discussion is that
	\[
	\varphi(X^\beta \cdot X^\alpha) = \varphi(X^{\beta + \alpha}) = 
	\begin{cases}
		\textrm{$X^{\beta + \alpha} = X^\beta \varphi(X^\alpha)$ \hspace{5mm} if $\alpha \in \sigma_H \cap M_H$},\\
		\textrm{\hspace{8mm}$0 = X^\beta \varphi(X^{\alpha})$ \hspace{4.5mm} if $\alpha \in (\sigma \cap M) \setminus (\sigma_H \cap M_H)$}.
	\end{cases}	
	\]
	The $k[\sigma_H \cap M_H]$-linearity of $\varphi$ now easily follows. Thus, $\varphi$ splits $k[\sigma_H \cap M_H] \hookrightarrow k[\sigma \cap M]$.
\end{proof}

\autoref{prop:hyperplane-section-split-F-regular} can be generalized to affine hyperplanes of $M_{\mathbb R}$. However, the proof is more involved.

\begin{theorem}
	\label{thm:affine-hyperplane-split-F-regular}
	Suppose $\sigma$ is a cone in \autoref{setup}. Let $H$ be an affine hyperplane of $M_{\mathbb R}$ that is defined over $M_{\mathbb Q}$. 
	Furthermore, assume that there exists an $\alpha \in \sigma \cap M \cap H$.
	Consider the hyperplane
\[
H' \coloneqq H - \alpha	= \{\textbf{r} - \alpha \colon \textbf{r} \in H\}.
\]
	Let $\sigma_{H'}$ denote the cone of the hyperplane $H'$ that is generated by $(\sigma \cap H) - \alpha$. If $k[\sigma \cap M]$ is split-$F$-regular, then $k[\sigma_{H'} \cap M]$ is split-$F$-regular.
\end{theorem}

\begin{proof}[Proof of \autoref{thm:affine-hyperplane-split-F-regular}]
	We may assume $\sigma_{H'} \cap M$ has a nonzero element as otherwise there is nothing to prove. 
	We begin the proof with a claim.

	\begin{claim}
		\label{claim:hyperplane-original-cone-points}
		Let $\beta \in H' \cap M$. Then the following are equivalent:
		\begin{enumerate}
			\item $\beta \in \sigma_{H'} \cap M$.\label{claim:hyperplane-original-cone-points.a}
			
			\item For all integers $\ell \gg 0$, $\beta + \ell\alpha \in \sigma$.\label{claim:hyperplane-original-cone-points.b}

			\item There exists an integer $\ell  > 0$ such that $\beta + \ell\alpha \in \sigma$.\label{claim:hyperplane-original-cone-points.c}
		\end{enumerate}
	\end{claim}

	\begin{proof}[Proof of \autoref{claim:hyperplane-original-cone-points}]
	\autoref{claim:hyperplane-original-cone-points.a} $\Rightarrow$ \autoref{claim:hyperplane-original-cone-points.b}: Choose $a_1,\dots,a_n \in \sigma \cap H$ and $r_1,\dots,r_n \in \mathbb{R}_{\geq 0}$ such that
	\[
	\beta =  r_1(a_1 - \alpha) + \dots + r_n(a_n - \alpha).	
	\]
	Then for any integer $\ell$ such that $\ell > r_1 + \dots + r_n$, we have
	\[
	\beta + \ell\alpha = r_1a_1 + \dots + r_na_n + \left(\ell - \sum_{i=1}^nr_i\right)\alpha \in \sigma
	\]
	because $a_1,\dots,a_n,\alpha \in \sigma$. 

	\autoref{claim:hyperplane-original-cone-points.b} $\Rightarrow$ \autoref{claim:hyperplane-original-cone-points.c} is clear.

	\autoref{claim:hyperplane-original-cone-points.c} $\Rightarrow$ \autoref{claim:hyperplane-original-cone-points.a}: Since $\sigma$ is closed under scaling, we have
	$
	(1/\ell)\beta + \alpha \in \sigma,	
	$
	and, since $H'$ is a hyperplane, we have
	$
	(1/\ell)\beta \in H'.
	$ 
	Thus, $(1/\ell)\beta + \alpha \in \sigma \cap H$, and hence,  $(1/\ell)\beta \in (\sigma \cap H) - \alpha$. Since $\sigma_{H'}$ is the cone in the hyperplane $H'$ over $(\sigma \cap H) - \alpha$, upon scaling by $\ell$ we get $\beta \in \sigma_{H'}$. Additionally, $\beta \in M$ by assumption, and so, $\beta \in \sigma_{H'} \cap M$. This completes the proof of \autoref{claim:hyperplane-original-cone-points}. \phantom\qedhere
	\end{proof}

	So let $\beta \in \sigma_{H'} \cap M$ be an \emph{arbitrary} element. Fix an integer $\ell \gg 0$ as above such that $\beta + \ell\alpha \in \sigma \cap M$. 
	Since $k[\sigma \cap M]$ is split-$F$-regular, there exists an integer $e > 0$ such the $k[\sigma\cap M]$-linear map 
	\begin{align*}
	\lambda_{\beta + \ell\alpha, e} \colon k[\sigma \cap M] &\to F^e_*k[\sigma \cap M]\\
	1 &\mapsto F^e_*X^{\beta + \ell\alpha}
	\end{align*}
	splits. Fix a Frobenius splitting
	$
	\phi \colon F^e_*k \to k	
	$
	of the $e$-th iterate of the Frobenius on $k$. Then by \autoref{lem:split-implies-torus-compatible-split}, the map
	\[
		\pi_{-(\beta +\ell\alpha)} \colon F^e_*k[M] \to k[M]
	\]
	given by the rule
	\begin{equation}\label{eq:cases-key}
		\pi_{-(\beta +\ell\alpha)}(F^e_*s X^\gamma) =
		\begin{cases}
		  \phi(F^e_*s) X^{(\gamma - (\beta + \ell\alpha))/p^{e}} &\quad\text{if $(\gamma - (\beta + \ell\alpha))/p^{e} \in M$,}\\
		  0 &\quad\text{otherwise,} 
		\end{cases}
	\end{equation}
	for all $s  \in k$ and $\gamma \in M$, restricts to yield a splitting of $\lambda_{\beta +\ell\alpha,e}$.

	Consider the composition
	\[
	\pi_{-\beta} \coloneqq 	F^e_*k[M] \xrightarrow{F^e_*X^{\ell\alpha}\cdot} F^e_*k[M] \xrightarrow{\pi_{-(\beta + \ell\alpha)}} k.
	\]
	For all $s \in K$ and $\gamma \in M$,
	\begin{equation*}
		\pi_{-\beta}(F^e_*s X^\gamma) =
		\begin{cases}
		  \phi(F^e_*s) X^{(\gamma - \beta)/p^{e}} &\quad\text{if $(\gamma - \beta)/p^{e} \in M$,}\\
		  0 &\quad\text{otherwise.} 
		\end{cases}	
	\end{equation*}
	In particular, $\pi_{-\beta}(F^e_*X^\beta) = \phi(F^e_*1)X^0 = 1$, since $\phi \colon F^e_*k \to k$ is a Frobenius splitting.

	\begin{claim}
		\label{claim:induced-splitting-hyperplane}
		The restriction of $\pi_{-\beta}$ to $F^e_*k[\sigma_{H'} \cap M]$ maps into $k[\sigma_{H'} \cap M]$.
	\end{claim}

	\begin{proof}[Proof of \autoref{claim:induced-splitting-hyperplane}]
		It suffices to show that for all $s \in k$ and $\beta' \in \sigma_{H'} \cap M$, $\pi_{-\beta}(F^e_*s X^{\beta'}) \in k[\sigma_{H'}\cap M]$. By the construction of $\pi_{-\beta}$, we have
		\[
			\pi_{-\beta}(F^e_*s X^{\beta'}) = 	
			\begin{cases}
				\phi(F^e_*s) X^{(\beta' - \beta)/p^{e}} &\quad\text{if $(\beta' - \beta)/p^{e} \in M$,}\\
				0 &\quad\text{otherwise.} 
			\end{cases}
		\]
		Thus, if $\beta' - \beta \notin p^eM$, we are already done. So suppose $\beta' - \beta \in p^eM$ and choose $\delta \in M$ such that 
		\begin{equation}
			\label{eq:situation-we-are-in}
		\beta' - \beta = p^e\delta.	
		\end{equation}
		Since $\beta, \beta' \in H'$ and $H'$ is a linear space, it follows that $\delta \in H'$ as well. Consequently, 
		$
		\delta \in H' \cap M.	
		$
		In order to finish the proof, we need to show that $\delta \in \sigma_{H'} \cap M$. Equivalently, it suffices to show by \autoref{claim:hyperplane-original-cone-points} that there exists an integer $i > 0$ such that 
		\begin{equation}\label{eq:need-to-show}
			\delta + i\alpha \in \sigma.	
		\end{equation}
		Since $\beta' \in \sigma_{H'} \cap M$, by \autoref{claim:hyperplane-original-cone-points} again, there exists an integer $\ell' \gg \ell$ such that 
		$
		\beta' + \ell'\alpha \in \sigma 
		$
		and such that $\ell' - \ell = ip^e$ for some integer $i > 0$. Then 
		\[
		(\beta' + \ell'\alpha) - (\beta + \ell\alpha) \stackrel{\autoref{eq:situation-we-are-in}}{=} p^e(\delta + i\alpha),
		\]
		and so,
		\[
		\pi_{-(\beta + \ell\alpha)}(F^e_*X^{\beta' + \ell'\alpha}) \stackrel{\autoref{eq:cases-key}}{=} X^{\delta + i\alpha}.	
		\]
 		Since $\beta' + \ell'\alpha \in \sigma \cap M$ by our choice of the integer $\ell'$ (note, $\beta', \alpha \in M$ so $\beta' + \ell'\alpha \in M$) and since $\pi_{-(\beta + \ell\alpha)}$ restricted to $F^e_*k[\sigma \cap M]$ maps into $k[\sigma \cap M]$, we get $X^{\delta + i \alpha} \in k[\sigma \cap M]$, or equivalently, that 
		$
		\delta + i\alpha \in \sigma \cap M.	
		$
		This establishes \autoref{eq:need-to-show}, thereby completing the proof of \autoref{claim:induced-splitting-hyperplane}. \phantom\qedhere
	\end{proof}

\autoref{claim:induced-splitting-hyperplane} shows that for any $\beta \in \sigma_{H'} \cap M$ there exists an integer $e > 0$ such that 
\begin{align*}
\lambda_{\beta,e} \colon k[\sigma_{H'}\cap M] &\to F^e_*k[\sigma_{H'} \cap M]\\
1 &\mapsto F^e_*X^\beta
\end{align*}
splits (hence is pure). Then $k[\sigma_{H'} \cap M]$ is split-$F$-regular by \autoref{cor:splitF-regular-equiv-Fpure-regular}.
\end{proof}

\subsection{Split-{$F$}-regularity and finite generation} We finally have all the tools to address \autoref{revised-question}. The main theorem of this subsection is:

\begin{theorem}
	\label{thm:SFR-cone-finite-generation}
	Under \autoref{setup} suppose that $k[\sigma \cap M]$ is split-$F$-regular. Then $\sigma$ is a closed cone generated by finitely many elements of $\sigma \cap M$. Consequently, $k[\sigma \cap M]$ is a finite type $k$-algebra.
\end{theorem}

In the proof of \autoref{thm:SFR-cone-finite-generation} we will reduce to the case where $\sigma$ is closed and strongly convex. The proof of finite generation in this setting will rely on the following technical result.

\begin{proposition}
	\label{prop:extremal-rays-cannot-accumulate}
	Under \autoref{setup} suppose $\sigma$ is closed and strongly convex. Let $d$ be the rank of the lattice $M$. Assume $k[\sigma \cap M]$ is split-$F$-regular. Also assume that for any subgroup $M'$ of $M$ of rank $d-1$ and for any submonoid $S$ of $M'$, if $k[S]$ is split-$F$-regular, then $S$ is finitely generated as a monoid.
	Let $\rho$ be a ray of $\sigma$ and suppose there exists $\varphi \in \sigma^\vee$ such that $\varphi^\perp$ is defined over $M_{\mathbb Q}$ and 
	$
	\rho \subseteq \varphi^{\perp}.	
	$
	Then $\rho$ cannot be an accumulation ray of infinitely many extremal rays of $\sigma$.
\end{proposition}

\begin{proof}
	We begin by fixing a reference hyperplane. Since $\sigma$ is strongly convex, the dual $\sigma^\vee$ is full-dimensional, and so we may choose a rational ray $\nu$ in the interior of $\sigma^\vee$. This gives a rational hyperplane $\nu^\perp$ with $\nu^\perp \cap \sigma = \{ 0 \}$ by \autoref{cor:Relint-dual-cone}. Since $\sigma$ is full-dimensional, $\Int(\sigma)$ is a non-empty open subset of $M_{\mathbb R}$ by \autoref{some-remarks}\autoref{some-remarks.d}. Thus, there exists $w \in \Int(\sigma) \cap M$ by \autoref{some-remarks}\autoref{some-remarks.e}. Let $c := \nu(w) > 0$, and let $H$ be the rational affine hyperplane $H:= \nu^\perp + w = \{ v \in V :  \nu(v) = c \}$. Note that $\Int(\sigma) \cap H$ is an open subset of $H$ containing $w$. Thus, $\sigma \cap H$ is a full-dimensional closed convex set in $H$. Moreover, $\sigma \cap H$ is bounded by \cite[Corollary 8.4.1]{RockafellarConvexAnalysis}.
	If  $\mu$ is a (non-trivial) ray of $\sigma$, we will denote by $\hat\mu := \mu \cap H$ its unique point of intersection with $H$. To prove the proposition, it suffices to exhibit an open neighborhood of $\hat\rho := \rho \cap H$ in $H$ that intersects at most finitely many extremal rays of $\sigma$.

	Consider first the hyperplane $\varphi^\perp$ and the resulting sub-monoid $\sigma \cap \varphi^\perp \cap M$ of $M \cap \varphi^\perp \cong \mathbb{Z}^{d-1}$. By \autoref{prop:hyperplane-section-split-F-regular}, we have that the corresponding monoid algebra $k[\sigma \cap \varphi^\perp \cap M]$ is split-$F$-regular, and so by assumption $\sigma \cap \varphi^\perp \cap M$ is a finitely generated monoid. Since $\sigma \cap \varphi^\perp$ is an extremal subcone of $\sigma$, it follows that the extremal rays of the closed strongly convex cone $\sigma \cap \varphi^\perp$ are precisely the extremal rays of $\sigma$ that are contained in $\varphi^\perp$. Thus, $\sigma \cap \varphi^\perp$ is generated by a subset of the extremal rays of $\sigma$. The latter are all rational by \autoref{extremal-rays-are-rational}\ref{extremal-rays-are-rational.extremalRay} since $k[\sigma \cap M]$ is split-$F$-regular. As a result, it follows that $\sigma \cap \varphi^\perp$ is rational polyhedral since $\sigma \cap \varphi^{\perp}$ must be generated by the finitely generated monoid $\sigma \cap \varphi^\perp \cap M$. 
	As such, 
	\begin{equation}
		\label{eq:definition-P}
		P := \sigma \cap \varphi^\perp \cap H	
	\end{equation}
	is a rational polytope because $\nu^\perp$ is also a reference hyperplane of the rational polyhedral cone $\sigma \cap \varphi^\perp$. Note that $\hat\rho \in P$, and let $p_0, \ldots, p_\ell \in P$ be the extremal points of $P$. We have $p_0, \ldots, p_\ell \in M_{\mathbb{Q}}$, and can choose $N \in \mathbb{Z}_{> 0}$ sufficiently large so that $Np_0, \ldots, Np_\ell \in M$.  Replacing $H$ by its parallel translate $\{ v \in V : \nu(v) = Nc \} = \nu^\perp + Nw$ exchanges $P$ with $NP$, and so we may assume without loss of generality that $p_0, \ldots, p_\ell \in M$. By \cite[Thm.\ 18.2]{RockafellarConvexAnalysis}, there exists a unique face $g$ of $P$ such that $\hat\rho$ is in the relative interior of $g$. Moreover, by \cite[Thm.\ 18.1]{RockafellarConvexAnalysis}, $g$ is the intersection of all faces of $P$ that contain $\hat\rho$, and hence is the smallest face of $P$ that contains $\hat\rho$.
	Being a face of a polytope, $g$ is the convex hull of those extremal points $p_i$ of $P$ such that $p_i \in g$ \cite[Prop.\ 2.2(i), Prop.\ 2.3(ii)]{ZieglerLecturesonPolytopes}. Thus, after reordering the extremal points of $P$ if necessary, we may assume that $p_0 \in g$. By our choice of $g$, any face of $P$ that contains $\hat\rho$ must also contain $p_0$.

	Let $H' := H - p_0$ and consider $\sigma_{H'} := \mathbb{R}_{\geq 0}((\sigma \cap H) - p_0)$ the cone that is generated by $(\sigma \cap H) - p_0$. Note that $H' = \nu^\perp$ since $H$ is a translate of $\nu^\perp$ and $p_0 \in H$. However, we are using notation to match that of \autoref{thm:affine-hyperplane-split-F-regular}, which we will soon apply.  The cone $\sigma_{H'}$ is a full-dimensional in $H' = \nu^\perp$. Indeed, $\sigma \cap H$ contains a non-empty open subset of $H$, and so, $(\sigma \cap H) - p_0$ (and hence $\sigma_{H'}$) contains a non-empty open subset of $H'$. Thus, the cone $\sigma_{H'}$ in the rational vector space $\nu^\perp = H'$ satisfies the hypotheses of \autoref{setup}. By \autoref{thm:affine-hyperplane-split-F-regular}, we have that the monoid algebra of $\sigma_{H'} \cap M \subseteq M \cap \nu^\perp \cong \mathbb{Z}^{d-1}$ is split-$F$-regular, and so again by our assumption $\sigma_{H'} \cap M$ is finitely generated. In particular, $\sigma_{H'}$ is rational polyhedral since it is generated by $\sigma_{H'} \cap M$ by \autoref{thm:split-F-regular-arbitrary-cone-closed}; let $p_{\ell+1}, \ldots, p_m \in \sigma \cap H \cap M_{\mathbb{Q}}$ be such that $\sigma_{H'} = \sum_{i=1}^m \mathbb{R}_{\geq 0}(p_i - p_0)$. Since $\sigma \cap \varphi^\perp \cap H$ is the convex hull of $p_0, \ldots, p_{\ell}$, we may assume without loss of generality that $p_{\ell+1}, \ldots, p_m \not\in \varphi^\perp$. Set $Q$ to be the convex hull of $p_0, \ldots, p_m$, so that also $\sigma_{H'} = \mathbb{R}_{\geq 0}(Q-p_0)$. Note that $p_0 \in Q$ because $p_0 \in P \subseteq Q$. Since $\sigma_{H'}$ is full-dimensional in $H'$, $Q - p_0$ is a full-dimensional convex set in $H'= \nu^\perp$. In addition, we have 
	\[P = Q \cap \varphi^\perp,\] 
	where $P$ is the polytope defined in \autoref{eq:definition-P}. Thus, $P$ is a face of $Q$ (here we use that $\varphi^\perp$ is a supporting hyperplane of $\sigma$, and hence, $Q$). Moreover, $P$ is a proper face of $Q$ because $Q$ is a full-dimensional convex set in $H$.

	Suppose that $f_1, \ldots, f_n$ are the facets of the rational polytope $Q - p_0$, where 
	\[
	f_j = (Q - p_0)\cap \{ v \in \nu^\perp : \psi_j(v) = \lambda_j \}
	\] 
	for $\psi_1, \ldots, \psi_n \in (\nu^\perp)^*$ and $\lambda_1, \ldots, \lambda_n \in \mathbb{Q}$ with $\psi_j(Q - p_0) \in (\lambda_j,\infty)$ for $j = 1, \ldots, n$. In other words, we have that $Q - p_0 = \bigcap_{j=1}^n \{ v \in \nu^\perp : \psi_j(v)\geq \lambda_j \}$ is an irreducible representation of $Q - p_0$ in the sense of \cite[Theorem 8.2]{BrondstedAnIntroductionToConvexPolytopes}. Since $0 \in Q - p_0$, note that $\lambda_j \leq 0$ for $j = 1, \ldots, n$, and we may reorder to ensure that 
	\[
	\textrm{$\lambda_1, \ldots, \lambda_s = 0$ while $\lambda_{s+1}, \ldots, \lambda_n < 0$}.
	\] 
	Note that
	\begin{equation}
		\label{eq:cone-intersection-half-spaces}
	\sigma_{H'} = \{ v \in \nu^\perp : \psi_1(v) \geq 0, \ldots, \psi_s(v) \geq 0 \}.
	\end{equation}
	Indeed, as $Q - p_0$ generates $\sigma_{H'}$, certainly $\sigma_{H'} \subseteq \{ v \in \nu^\perp : \psi_1(v) \geq 0, \ldots, \psi_s(v) \geq 0 \}$. On the other hand, if $v \in \nu^\perp$ with $\psi_1(v) \geq 0, \ldots, \psi_s(v) \geq 0$, we can take $0 < \epsilon \ll 1$ sufficiently small so that $\psi_{s+1}(\epsilon v) > \lambda_{s+1}, \ldots, \psi_n(\epsilon v) > \lambda_n$. This gives that $\epsilon v \in \bigcap_{j=1}^n \{ v \in \nu^\perp : \psi_j(v)\geq \lambda_j \} = Q - p_0$, so that $v \in \frac{1}{\epsilon}(Q - p_0) \subseteq \sigma_{H'}$. It follows that $\sigma_{H'} = \{ v \in \nu^\perp : \psi_1(v) \geq 0, \ldots, \psi_s(v) \geq 0 \}$ as claimed.

	Since $P = Q \cap \varphi^\perp$ is a (proper) face of $Q$, we have that $P - p_0$ is a face of $Q - p_0$, and thus $g - p_0$ is the smallest face of $Q - p_0$ containing $\hat\rho - p_0$ (recall, $g$ was the smallest face of $P$ containing $\hat\rho$). Since $g - p_0$ is a proper face, it is necessarily the intersection of the facets of $Q - p_0$ that contain it (equivalently, that contain $\hat\rho - p_0$), that is,
	\begin{equation*}
		g - p_0 = \bigcap_{\substack{j = 1, \ldots, n \mbox{ \scriptsize satisfying } \\ \psi_j(\hat\rho - p_0) = \lambda_j}} (Q - p_0) \cap \{ v \in \nu^\perp : \psi_j(v) = \lambda_j \}
	\end{equation*}
	(see, for example, \cite[Theorem 5.6, Theorem 8.3]{BrondstedAnIntroductionToConvexPolytopes} and \cite[Lemma 2.9, Exercises 2.14 \& 2.15]{ZieglerLecturesonPolytopes}).
	Since $0 = p_0 - p_0 \in g - p_0$, we must have that $0 = \psi_j(0) = \lambda_j$ whenever $\psi_j(\hat\rho - p_0) = \lambda_j$ is satisfied, that is, the index $j$ must be in the set $\{1, \ldots, s\}$ whenever $f_j$ is a facet of $Q - p_0$ that contains $g - p_0$. Note that there exists $\psi_j$ such that $\psi_j(\hat\rho - p_0) > \lambda_j$. Otherwise, $\hat\rho - p_0$ is in the intersection of all facets of $Q - p_0$, which is the empty set. Then consider the open set
	\begin{equation*}
		U = \bigcap_{\substack{j = 1, \ldots, n \mbox{ \scriptsize satisfying } \\ \psi_j(\hat\rho - p_0) > \lambda_j}} \{ v \in \nu^\perp : \psi_j(v) > \lambda_j \}
	\end{equation*}
	of $\nu^\perp$, which is nonempty as $\hat\rho - p_0 \in U$. Its translate $p_0 + U \subseteq H$ is an open neighborhood of $\hat\rho$ in $H$, and the proposition follows provided we can show that $p_0 + U$ intersects at most finitely many extremal rays of $\sigma$. Assume $\tau$ is an extremal ray of $\sigma$ intersecting $p_0 + U$, so that $\hat\tau = H \cap \tau \in p_0+U$. This gives $\hat\tau - p_0 \in U$ and also $\hat\tau - p_0 \in ((\sigma \cap H) - p_0) \subseteq \sigma_{H'}$. Let us now check that $\hat\tau - p_0 \in Q - p_0 = \bigcap_{j=1}^n \{ v \in \nu^\perp : \psi_j(v)\geq \lambda_j \}$. For any $j \in \{1, \ldots, n \}$ with $\psi_j(\hat\rho - p_0) > \lambda_j$, we have $\psi_j(\hat\tau-p_0) > \lambda_j$ as $\hat\tau - p_0 \in U$. On the other hand, for any $j \in \{1, \ldots, n \}$ with $\psi_j(\hat\rho - p_0) = \lambda_j$ we must have $j \in \{1, \ldots, s \}$, from which it follows by \autoref{eq:cone-intersection-half-spaces} that  $\psi_j(\hat\tau - p_0) \geq 0 = \lambda_j$ as $\hat\tau - p_0 \in \sigma_{H'}$. Thus, we see $\hat\tau - p_0 \in Q - p_0$ or $\hat\tau \in Q$. Since $Q$ is the convex hull of $p_0, \ldots, p_m$, we have
	\begin{equation*}
		\hat\tau = \xi_0 p_0 + \cdots + \xi_m p_m
	\end{equation*}
	for some $ \xi_0, \ldots, \xi_m \in \mathbb{R}_{\geq 0}$ with $\sum_{i=0}^m \xi_i = 1$. As the sum of $\xi_0p_0, \ldots, \xi_m p_m \in \sigma$ belongs to the extremal ray $\tau$, we must have that all $\xi_0p_0, \ldots, \xi_m p_m \in \tau$. Not all of these can be zero as they sum to $\hat\tau \neq 0$, so we must have $0 \neq \xi_i p_i \in \tau$ for some $i$ giving $\tau = \mathbb{R}_{\geq 0} p_i$ (and $\hat\tau = p_i$). Thus, $\tau$ is among the rays $\mathbb{R}_{\geq 0}p_0, \ldots, \mathbb{R}_{\geq 0}p_m$, so we see only finitely many extremal rays of $\sigma$ can intersect $p_0 + U$ giving our desired conclusion.
\end{proof}

We will now prove \autoref{thm:SFR-cone-finite-generation}

\begin{proof}[Proof of \autoref{thm:SFR-cone-finite-generation}]
By Gordan's Lemma (see, for instance, \cite[Lem.\ A.6.1]{deFernexEinMustataBook}), $k[\sigma \cap M]$ will be a finite type $k$-algebra if we can show that $\sigma$ is a closed cone generated (as a convex cone) by finitely many elements of $\sigma \cap M$. We proceed by strong induction on the dimension of $M_{\mathbb R}$, the ambient vector space of the cone. The case when $\dim_{\mathbb R}(M_{\mathbb R}) = 1$ is clear. So assume $\dim_{\mathbb R}(M_{\mathbb R}) \geq 2$.

Note $\sigma$ is always closed by \autoref{thm:split-F-regular-arbitrary-cone-closed}. If $L$ is the lineality space of $\sigma$, then $L$ is defined over $M_{\mathbb Q}$ by \autoref{extremal-rays-are-rational}\ref{extremal-rays-are-rational.lineality}. In particular, $L$ is generated as a cone by finitely many elements of $\sigma \cap M$ since it is a finite dimensional vector space. Let $\pi \colon M_{\mathbb R} \twoheadrightarrow M_{\mathbb R}/L$ be the projection map. Since $\sigma$ is closed and the linear span of $\sigma$ equals $M_{\mathbb R}$, it follows that $\pi(\sigma)$ is a closed strongly convex cone (see \autoref{prop:lineality-quotient-strong-convex}) whose linear span is $M_{\mathbb R}/L$. In other words, $\pi(\sigma)$ satisfies the hypotheses of \autoref{setup}. Moreover, the monoid algebra $k[\pi(\sigma) \cap \pi(M)]$ is split-$F$-regular with respect to the lattice $\pi(M)$ of $M_{\mathbb R}/L$ by \autoref{lem:cone-ses}\ref{lem:cone-ses.splitF-regular}. If $\dim_{\mathbb R}(L) \geq 1$, then $\pi(\sigma)$ is generated by finitely many elements of $\pi(\sigma) \cap \pi(M)$ by the induction hypothesis applied to $\dim_{\mathbb R}(M_{\mathbb R}/L)$. Then $\sigma$ is generated as a convex cone by a finite subset of $\sigma \cap M$ by \autoref{lem:cone-ses}\ref{lem:cone-ses.generators}.

Thus, we reduce to the case where $L = 0$, or equivalently, when $\sigma$ is strongly convex. We already know by \autoref{extremal-rays-are-rational}\ref{extremal-rays-are-rational.extremalRay} that the extremal rays of $\sigma$ are all rational. Assume for contradiction that $\sigma$ has infinitely many extremal rays. Then by a compactness argument one sees that some sequence of extremal rays of $\sigma$ accumulate to a limiting ray $\rho$. Since the linear span of $\sigma$, which by assumption is $M_{\mathbb R}$, has dimension $\geq 2$, every extremal ray of $\sigma$ lies on a proper face of $\sigma$ by \autoref{extremal-ray-properties}\ref{extremal-ray-properties.dimlinealitygeq2}. Since the union of the proper faces of $\sigma$ is the topological boundary of $\sigma$ in $M_{\mathbb R}$ by \autoref{rel-int-as-comp} (recall $\Int(\sigma)$ is open in $M_{\mathbb R}$ since $\sigma$ is full-dimensional), it follows that $\rho$ is also in the topological boundary of $\sigma$, that is, $\rho$ is in some proper face of $\sigma$. Let $\varphi \in \sigma^\vee$ such that 
\[
\rho \subseteq \sigma \cap \varphi^\perp.	
\]
We now argue that we may choose $\varphi$ to be in $M^*$, that is, $\varphi$ can be chosen to be rational.
Recall that the extremal rays of $\sigma^\vee$ are rational by \autoref{thm:dual-cone-rational-SFR}. Let $M^*$ be the dual lattice of $M$ in $M_{\mathbb R}^*$. Then there exist $\{\varphi_1, \dots, \varphi_n\} \subseteq \sigma^\vee \cap M^*$ and real numbers $\lambda_1,\dots,\lambda_n > 0$ such that 
\[
\varphi = \sum_{i = 1}^n \lambda_i\varphi_i.	
\]
Since $\rho \subseteq \Supp(\varphi_i)$ for all $i$, the only way for $\rho$ to be a subset of $\varphi^\perp$ is for $\rho$ to be contained in $\varphi_i^\perp$ for all $i$. Replacing $\varphi$ by one of the $\varphi_i$'s, this shows that there exists $\varphi \in \sigma^\vee \cap M^*$ such that $\rho \subseteq \varphi^\perp$. Note that $\varphi^\perp$ is defined over $M_{\mathbb Q}$ because $\varphi \in M^*$. 

\autoref{prop:extremal-rays-cannot-accumulate} then leads to a contradiction because by strong induction it follows that for any subgroup $M'$ of $M$ of strictly smaller rank, and for any submonoid $S$ of $M'$, if $k[S]$ is split-$F$-regular, then $S$ is a finitely generated monoid. Indeed, let $\mathbb{Z}S$ be the subgroup of $M'$ generated by $S$, let $\mathbb{R}S$ be the subspace of $M_{\mathbb R}$ generated by $S$ and let $\sigma_S$ be the convex cone in $\mathbb{R}S$ that is generated by $S$. Then $\mathbb{Z}S$ is a lattice of $\mathbb{R}S$, and by \autoref{normal-semigroup}\ref{normal-semigroup.3} and \autoref{integrally-closed-cones} it follows that if $k[S]$ is split-$F$-regular, then
\[
S = \sigma_S \cap \mathbb{Z}S.	
\] 
Thus, $k[\sigma_S \cap \mathbb{Z}S]$ is split-$F$-regular.
Now $\dim_{\mathbb R}(\mathbb{R}S)$ is the rank of the free Abelian group $\mathbb{Z}S$. Since $\mathbb{Z}S \subseteq M'$ and $M'$ is a subgroup of $M$ of strictly smaller rank, it follows that $\dim_{\mathbb R}(\mathbb{R}S) < \dim_{\mathbb R}(M_{\mathbb R})$, because the latter coincides with the free rank of $M$. Furthermore, $\sigma_S$ is full-dimensional in $\mathbb{R}S$ and $0 \in \sigma_S$, that is, $\sigma_S$ satisfies the hypotheses of \autoref{setup}. Therefore, by induction we get that $k[\sigma_S \cap \mathbb{Z}S]$ is a finite type $k$-algebra, and consequently, that $S$ is a finitely generated monoid.
\end{proof}

Finally, we can deduce the main result about split-$F$-regular monoid algebras, that is, we obtain an affirmative answer to \autoref{ques-split-monoid-algebra}.

\begin{theorem}
	\label{thm:split-F-regular-monoid-finite-generation}
	Let $S$ be a submonoid of a free Abelian group of finite rank. If the monoid algebra $k[S]$ is split-$F$-regular, then $k[S]$ is a finitely generated $k$-algebra. Equivalently, $S$ is a finitely generated Abelian monoid.
\end{theorem}

\begin{proof}
	Both in \autoref{subsec:monoids-to-cones} and in the last paragraph of the proof of \autoref{thm:SFR-cone-finite-generation} we have explained how this Theorem follows from \autoref{thm:SFR-cone-finite-generation}.
\end{proof}

We showed in \autoref{cor:splitF-regular-equiv-Fpure-regular} that the property of being split-$F$-regular is independent of the ground field as long as we are working in a fixed positive characteristic. We now get a better result because of finite generation.

\begin{corollary}
	\label{cor:split-F-regular-independence-field-characteristic}
	Let $S$ be a submonoid of a free Abelian group $L$ of finite rank. Let $k$ be a field of characteristic $p > 0$ such that $k[S]$ is split-$F$-regular. Then for any prime number $q > 0$ and for any field $K$ of characteristic $q$, $K[S]$ is Noetherian and split-$F$-regular.
\end{corollary}

\begin{proof}
	$S$ is a finitely generated monoid by \autoref{thm:split-F-regular-monoid-finite-generation}, and $S$ is normal by \autoref{normal-semigroup}\ref{normal-semigroup.3}.

	Thus, $K[S]$ is a finite type $K$-algebra for any field $K$. Let $K^a$ be the algebraic closure of $K$. Since $K[S]$ is a direct summand of $K^a[S]$, it suffices to show that if $K$ is an algebraically closed field of characteristic $q$, then $K[S]$ is split-$F$-regular (or, equivalently, that $K[S]$ is strongly $F$-regular since it is Noetherian and $F$-finite).

	 Let $S_0$ be the group of invertible elements in $S$. By \cite[Prop.\ 6.1.3(a)]{BrunsHerzog} we have a decomposition
	\[
	S = S_0 \oplus S',	
	\]
	where $S'$ is a \emph{positive} normal monoid (i.e. $0$ is the only invertible element of $S'$). Note that since $S$ is finitely generated, both $S_0$ and $S'$ are finitely generated monoids as well. Thus, $S_0 \cong \mathbb{Z}^{\oplus d}$ and we also have by \cite[Prop.\ 6.1.3(b)]{BrunsHerzog} that
	\[
	K[S] \cong K[\mathbb{Z}^{\oplus d}] \otimes_K K[S'].	
	\]
	Since $K[\mathbb{Z}^{\oplus d}] \cong K[x_1,\dots,x_d]\big{[}\frac{1}{x_1\cdots x_d}\big{]}$, for indeterminates $x_1,\dots,x_d$ over $K$, it follows that 
	\[
	K[S'] \hookrightarrow K[S]
	\] 
	is a smooth extension of Noetherian $F$-finite domains. Thus, in order to show that $K[S]$ is strongly $F$-regular, it suffices to show by ascent of strong $F$-regularity along regular maps of Noetherian $F$-finite rings \cite[Lem.\ 4.11]{LyubeznikSmithStrongWeakFregularityEquivalentforGraded} (see also \cite{HochsterHunekeFRegularityTestElementsBaseChange}) that $K[S']$ is strongly $F$-regular. 
	
	In other words, we can assume $S$ is a positive finitely generated normal commutative monoid. Then by \cite[Prop.\ 1(3) \& Lem.\ 1]{HochsterRingsInvariantsofTori} (see also, \cite[Ex.\ 6.1.10]{BrunsHerzog}), there exists an integer $n > 0$ such that $S$ can be realized as a full sub-semigroup of $\mathbb{N}^n$ and such that the inclusion $K[S] \hookrightarrow K[\mathbb{N}^n]$ splits. But $K[\mathbb{N}^n]$ is a polynomial ring over $K$, which is strongly $F$-regular. Thus, $K[S]$ is strongly $F$-regular as well.
\end{proof}

\subsection{An application to graded rings associated with valuations}

Let $K/k$ be a finitely generated field extension and $\nu$ be a valuation of $K/k$ with valuation ring $(\mathcal{O}_\nu, \mathfrak{m}_\nu, \kappa_\nu)$. Suppose $\nu$ is \emph{centered} on a (possibly non-Noetherian) local $k$-algebra domain $(R,\mathfrak{m}_R, \kappa_R)$ with fraction field $K$. This means the local ring $(\mathcal{O}_\nu, \mathfrak{m}_\nu, \kappa_\nu)$ dominates the local ring $(R,\mathfrak{m}_R, \kappa_R)$. Then 
\[
\Phi_R \coloneqq \nu(R - \{0\})
\]
is a sub-monoid of the value group $\Gamma_\nu$ of $\nu$. Since the fraction field of $R$ is $K$, it follows that $\Phi_R$ generates $\Gamma_\nu$ as a group. For any $m \in \Phi_R$, the \textbf{valuation ideal of order $m$ along $R$} is the ideal
\[
\mathfrak{a}_{\geq m}(R) \coloneqq \{r \in R \colon \nu(r) \geq m\} \cup \{0\}.	
\]
Additionally, we define $\mathfrak{a}_{>m}(R) \coloneqq \{r \in R \colon \nu(r) > m\} \cup \{0\}$. This is also an ideal of $R$ and it is contained in $\mathfrak{a}_{\geq m}(R)$. The \textbf{associated graded ring of $R$ along $\nu$} is the $\Phi_R$-graded $\kappa_R$-algebra
\[
\gr_{\nu}(R) \coloneqq \bigoplus_{m \in \Phi_R} \mathfrak{a}_{\geq m}(R)/\mathfrak{a}_{>m}(R).	
\]
This ring is an integral domain. Moreover, if $(R,\mathfrak{m}_R) \hookrightarrow (S, \mathfrak{m}_S)$ is a local extension of centers of $\nu$, then we get an induced extension of associated graded rings 
\[
\gr_\nu(R) \hookrightarrow \gr_\nu(S)	
\]
because for $m \in \Phi_R \subseteq \Phi_S$, we have $\mathfrak{a}_{\geq m}(R) = \mathfrak{a}_{\geq m}(S) \cap R$ (resp. $\mathfrak{a}_{>m}(R) = \mathfrak{a}_{> m}(S) \cap R$). In particular, we always have an extension
\[
\gr_\nu(R) \hookrightarrow \gr_\nu(\mathcal{O}_\nu)	
\]
whose image is contained in the subring 
\[
	\bigoplus_{m \in \Phi_R} \mathfrak{a}_{\geq m}(\mathcal{O}_\nu)/\mathfrak{a}_{> m}(\mathcal{O}_\nu)
\]
of $\gr_\nu(\mathcal{O}_\nu)$.

If the center $R$ is essentially of finite type over $k$, then by \cite[Prop.\ 3.5]{TeissierValuationsToric} (see also \cite[Rem.\ 3.7]{SpivakovskyLinearArtinApproximation}), 
\[
\td{\Frac(\gr_\nu(R))/\kappa_R} = \dim_{\mathbb Q}(\mathbb Q \otimes_{\mathbb Z} \Gamma_\nu) + \td{\kappa_\nu/\kappa_R}.	
\]
The right-hand side of the above equality is finite because by \textbf{Abhyankar's inequality} over Noetherian local centers \cite[Thm.\ 1]{AbhyankarValuationsLocalDomain}, we have
\begin{equation}
	\label{eq:Abhyankar-inequality}
	\dim_{\mathbb Q}(\mathbb Q \otimes_{\mathbb Z} \Gamma_\nu) + \td{\kappa_\nu/\kappa_R} \leq \dim(R).
\end{equation}
Since $\dim(R) + \td{\kappa_R/k} = \td{K/k}$, \autoref{eq:Abhyankar-inequality} is also equivalent to the inequality
\[
	\dim_{\mathbb Q}(\mathbb Q \otimes_{\mathbb Z} \Gamma_\nu) + \td{\kappa_\nu/k} \leq \td{K/k}.
\]
A valuation $\nu$ for which equality holds in either of the above two inequalities is called an \textbf{Abhyankar valuation of $K/k$}. Such valuations are higher rational rank analogs of valuations associated to prime divisors on normal models of $K/k$, and they enjoy many of the good properties of their divisorial counterparts. For example, if $\nu$ is an Abhyankar valuation of $K/k$, then $\Gamma_\nu$ is a finitely generated Abelian group (hence it is free since it is torsion-free) and $\kappa_\nu$ is a finitely generated extension of $k$ \cite[Thm.\ 1]{AbhyankarValuationsLocalDomain}. 

Of particular interest to us will be Abhyankar valuations of maximal rational rank (see \autoref{thm:finite-generation-SFR-valuation-algebra-max-rank}). By Abhyankar's inequality, the maximal rational rank a valuation $\nu$ of $K/k$ can have is $\td{K/k}$, in which case $\nu$ is an Abhyankar valuation of $K/k$ with the property that $\kappa_\nu$ is an algebraic extension of $k$. Note that if a maximal rational rank valuation $\nu$ is centered on a Noetherian local ring $(R, \mathfrak{m}_R, \kappa_R)$ that is essentially of finite type over $k$ with fraction field $K$, then by equality in Abhyankar's inequality, $\dim(R) = \dim_{\mathbb Q}(\mathbb Q \otimes_{\mathbb Z} \Gamma_\nu) = \td{K/k}$, that is, $R$ is the local ring of a closed point of a model of $K/k$.

An important problem in valuation theory as well as birational geometry is to understand the structure of $\gr_\nu(R)$. For instance, one can ask when this graded ring is a finitely generated $\kappa_R$-algebra. Indeed, a positive answer to the finite generation problem for certain special Abhyankar valuations arising in  algebraic $K$-stability \cite{LiuXuZhuangFiniteGeneration} was the main ingredient in the recent proof of the projectivity of $K$-moduli spaces and the solution of the Yau-Tian-Donaldson conjecture for log-Fano pairs. Further applications of $\gr_\nu(R)$ arise in Teissier's approach to the long-standing local uniformization problem in positive characteristic \cite{TeissierMonomialIdealsBinomialIdeals,TeissierOvereightDeformations}. In general, finite generation of $\gr_\nu(R)$ fails even when $\nu$ is a divisorial valuation of $K/k$; see \cite{CossartGalindoPiltantNonFiniteGenerationDivisorial} for an example when $\td{K/k} = 3$. However, divisorial valuations have the smallest rational rank. Thus, one may wonder if $\gr_\nu(R)$ is a finitely generated $\kappa_R$-algebra when $\nu$ is a maximal rational rank valuation. One reason for this is that if $R$ is a regular center of $\nu$ with a regular system of parameters $x_1,\dots,x_d$ whose valuations $\nu(x_1),\dots,\nu(x_d)$ (freely) generate $\Gamma_\nu$ (i.e. $R$ is a monomial center of $\nu$), then $\gr_\nu(R)$ is isomorphic to a polynomial ring in $d$ variables over $\kappa_R$. Unfortunately, finite generation of $\gr_\nu(R)$ also fails in general for maximal rational rank valuations. Examples exist in function fields of transcendence degree $2$ when $R$ is a normal hypersurface singularity by work of Cutkosky \cite[Ex.\ 3.1]{CutkoskyFiniteGeneration} and for function fields of transcendence degree $\geq 3$ for $R$ regular \cite[Ex.\ 3.2]{CutkoskyFiniteGeneration}.

Nevertheless, assuming split-$F$-regularity of $\gr_\nu(R)$, we obtain the following finite generation result for maximal rational rank valuations with arbitrary centers. 

\begin{theorem}
	\label{thm:finite-generation-SFR-valuation-algebra-max-rank}
	Let $K/k$ be a finitely generated field extension, where $k$ is an algebraically closed field of characteristic $p > 0$. Let $\nu$ be a valuation of $K/k$ of maximal rational rank, that is, 
	\[
		\dim_{\mathbb Q}(\mathbb{Q} \otimes_{\mathbb Z} \Gamma_\nu) = \td{K/k}.
	\]
	Suppose $\nu$ is centered on a Noetherian local $k$-subalgebra $(R, \mathfrak{m}_R, \kappa_R)$ of $K$ such that the fraction field of $R$ is $K$. If $\gr_\nu(R)$ is split-$F$-regular, then $\gr_\nu(R)$ is a finitely generated $k$-algebra.
\end{theorem}

The proof of \autoref{thm:finite-generation-SFR-valuation-algebra-max-rank} will be a straightforward consequence of the following result that applies to non-maximal rank valuations with trivial residue field extensions.

\begin{proposition}
	\label{prop:graded-pieces-valuation-algebra}
	Let $K/k$ be a finitely generated field extension and $\nu$ be a valuation of $K/k$ with valuation ring $(\mathcal{O}_\nu, \mathfrak{m}_\nu, \kappa_\nu)$. Suppose $\nu$ is centered on a (possibly non-Noetherian) local $k$-subalgebra $(R, \mathfrak{m}_R, \kappa_R)$ of $K$ such that the fraction field of $R$ is $K$. Let $m \in \nu(R - \{0\})$ and let $f \in \mathcal{O}_\nu$ such that $\nu(f) = m$. Then we have the following:
	\begin{enumerate}[label=(\arabic*)]
		\item The map 
		\begin{align*}
			\phi_f \colon \mathfrak{a}_{\geq m}(R) &\to \kappa_\nu\\
			g &\mapsto (g/f) + \mathfrak{m}_\nu
			\end{align*}
			is $R$-linear with kernel $\mathfrak{a}_{> m}(R)$.\label{prop:graded-pieces-valuation-algebra.1}

		\item If $R = \mathcal{O}_\nu$, then $\phi_f$ induces an isomorphism of $\kappa_\nu$-vector spaces
		\[
		\mathfrak{a}_{\geq m}(\mathcal{O}_\nu)/\mathfrak{a}_{> m}(\mathcal{O}_\nu) \xrightarrow{\simeq} \kappa_\nu.	
		\]
		\label{prop:graded-pieces-valuation-algebra.2}

		\item Suppose $\Gamma_\nu$ is the value group of $\nu$ and $\Phi_{\mathcal{O}_\nu} \coloneqq \nu(\mathcal{O}_\nu - \{0\})$ is the sub-monoid of non-negative elements of $\Gamma_\nu$. If $\Gamma_\nu$ is finitely generated, then
		\[
		\gr_\nu(\mathcal{O}_\nu) \cong	\kappa_\nu[\Phi_{\mathcal{O}_\nu}],
		\]
		where $\kappa_\nu[\Phi_{\mathcal{O}_\nu}]$ is the monoid algebra over $\kappa_\nu$ associated with $\Phi_{\mathcal{O}_\nu}$. Furthermore, $\gr_\nu(\mathcal{O}_\nu)$ is integrally closed in its fraction field.\label{prop:graded-pieces-valuation-algebra.3}

		\item Let $k$ be a field of characteristic $p > 0$. Suppose $\kappa_\nu = \kappa_R$ and $\Gamma_\nu$ is finitely generated. If $\gr_\nu(R)$ is split-$F$-regular, then $\gr_\nu(R)$ is a finitely generated $\kappa_R$-algebra.\label{prop:graded-pieces-valuation-algebra.4}
	\end{enumerate}
\end{proposition}

\begin{proof}[Proof of \autoref{prop:graded-pieces-valuation-algebra}]
	\ref{prop:graded-pieces-valuation-algebra.1} The $R$-linearity of $\phi_f$ is clear. For $g \in \mathfrak{a}_{\geq m}(R)$, we have $\phi_f(g) = 0$ precisely when $g/f \in \mathfrak{m}_\nu$, or equivalently, when $\nu(g) > \nu(f) = m$. Thus, $\ker(\phi_f) = \mathfrak{a}_{> m}(R)$.

	\ref{prop:graded-pieces-valuation-algebra.2} By \ref{prop:graded-pieces-valuation-algebra.1}, we have an injective $\mathcal{O}_\nu$-linear (hence also $\kappa_\nu$-linear) map 
	\[
		\mathfrak{a}_{\geq m}(\mathcal{O}_\nu)/\mathfrak{a}_{> m}(\mathcal{O}_\nu) \hookrightarrow \kappa_\nu.
	\]
	Since $\mathfrak{a}_{\geq m}(\mathcal{O}_\nu)/\mathfrak{a}_{> m}(\mathcal{O}_\nu)$ is a non-trivial $\kappa_\nu$-vector space, this map must be an isomorphism. Note that it makes sense that $\mathfrak{a}_{\geq m}(\mathcal{O}_\nu)/\mathfrak{a}_{> m}(\mathcal{O}_\nu)$ is a $1$-dimensional $\kappa_\nu$-vector space because the ideal $\mathfrak{a}_{\geq m}(\mathcal{O}_\nu)$ of $\mathcal{O}_\nu$ is principal; it is generated by any $f \in \mathcal{O}_\nu$ such that $\nu(f) = m$.

	\ref{prop:graded-pieces-valuation-algebra.3} We note that \ref{prop:graded-pieces-valuation-algebra.3} appears in \cite[Rem.\ 2.1]{CutkoskyFiniteGeneration} and in greater generality in \cite{BarnabeNovacoskiSpivakovskyChoiceFunction}. Nevertheless, we include a proof here for the reader's convenience, especially since the argument is straightforward in our setting. Totally ordered abelian groups are torsion-free. Thus, $\Gamma_\nu$ is a finitely generated free Abelian group by our hypothesis. Let $y_1,\dots,y_n \in K$ such that $\nu(y_1), \dots, \nu(y_n)$ freely generate $\Gamma_\nu$. Then for any $m \in \Phi_{\mathcal{O}_\nu}$, there exist unique $a_1,\dots,a_n \in \mathbb{Z}$ such that 
	\[
	\nu(y^{a_1}_1\cdots y^{a_n}_n) = m.	
	\]
	Define $y_m \coloneqq y^{a_1}_1\cdots y^{a_n}_n$.
	Note that $y_m \in \mathfrak{a}_{\geq m}(\mathcal{O}_\nu)$, and by \ref{prop:graded-pieces-valuation-algebra.2}, the class of $y_m$ in $\mathfrak{a}_{\geq m}(\mathcal{O}_\nu)/\mathfrak{a}_{> m}(\mathcal{O}_\nu)$ generates the latter as a $\kappa_\nu$-vector space. Furthermore, the elements $y_m$ are compatible with the monoid structure of $\Phi_{\mathcal{O}_\nu}$ in the sense that if $m_1, m_2 \in \Phi_{\mathcal{O}_\nu}$, then
	\[
	y_{m_1 + m_2} = y_{m_1} \cdot y_{m_2}.	
	\]
	Thus, by the universal property of $\kappa_\nu[\Phi_{\mathcal{O}_\nu}]$, we get an induced $\kappa_\nu$-algebra map 
	\[
		\kappa_\nu[\Phi_{\mathcal{O}_\nu}] \to \gr_\nu(\mathcal{O}_\nu),
	\]
	that for $m \in \Phi_{\mathcal{O}_\nu}$ sends the corresponding monoid algebra element $X^m \mapsto$ the class of $y_m$ in $\mathfrak{a}_{\geq m}(\mathcal{O}_\nu)/\mathfrak{a}_{> m}(\mathcal{O}_\nu)$. By \ref{prop:graded-pieces-valuation-algebra.2} it follows that this $\kappa_\nu$-algebra map is an isomorphism. Finally, $\Phi_{\mathcal{O}_\nu} = \{\gamma \in \Gamma_\nu \colon \gamma \geq 0\}$ is a normal sub-monoid of $\Gamma_\nu$ because if $\gamma \in \Gamma_\nu$ and $n > 0$ is an integer such that $n\gamma \geq 0$, then $\gamma \geq 0$. Then $\gr_\nu(\mathcal{O}_\nu) = \kappa_\nu[\Phi_{\mathcal{O}_\nu}]$ is integrally closed in its fraction field by \cite[Cor.\ 12.11]{GilmerCommutativeSemigroupRings}.

	\ref{prop:graded-pieces-valuation-algebra.4} As a consequence of \ref{prop:graded-pieces-valuation-algebra.2} and our assumption that $\kappa_R = \kappa_\nu$, we have that for $m \in \nu(R - \{0\})$, the $\kappa_R$-linear inclusion 
	\[
	\mathfrak{a}_{\geq m}(R)/\mathfrak{a}_{>m}(R) \hookrightarrow \mathfrak{a}_{\geq m}(\mathcal{O}_\nu)/\mathfrak{a}_{> m}(\mathcal{O}_\nu),	
	\]
	is an isomorphism. It then follows that the inclusion of associated graded rings 
	$
	\gr_\nu(R) \hookrightarrow \gr_\nu(\mathcal{O}_\nu)	
	$
	induces an isomorphism of $\gr_\nu(R)$ with the subring 
	\[
	\bigoplus_{m \in \Phi_R} \mathfrak{a}_{\geq m}(\mathcal{O}_\nu)/\mathfrak{a}_{> m}(\mathcal{O}_\nu)
	\]
	of $\gr_\nu(\mathcal{O}_\nu)$, where $\Phi_R$ is the monoid $\nu(R - \{0\})$. The isomorphism 
	\[
		 \gr_\nu(\mathcal{O}_\nu) \cong \kappa_\nu[\Phi_{\mathcal{O}_\nu}] = \kappa_R[\Phi_{\mathcal{O}_\nu}]
	\]
	of \ref{prop:graded-pieces-valuation-algebra.3} then gives us that 
	\[
		\bigoplus_{m \in \Phi_R} \mathfrak{a}_{\geq m}(\mathcal{O}_\nu)/\mathfrak{a}_{> m}(\mathcal{O}_\nu) \cong \kappa_R[\Phi_R].	
	\]
	Thus, $\gr_\nu(R) \cong \kappa_R[\Phi_R]$. Moreover, $\Phi_R$ is a sub-monoid of the free Abelian group $\Gamma_\nu$ of finite rank. Hence, the split-$F$-regularity of $\gr_\nu(R)$ implies that $\gr_\nu(R)$ is a finitely generated $k_R$-algebra by \autoref{thm:split-F-regular-monoid-finite-generation}.
\end{proof}

\begin{remark}
	\label{rem:finite-generation-value-group}
	Suppose we are in the situation of \autoref{prop:graded-pieces-valuation-algebra}. Then a necessary condition in order for $\gr_\nu(R)$ to be a finitely-generated $\kappa_R$-algebra is for the value group $\Gamma_\nu$ of $\nu$ to be a finitely generated group. Thus, our assumption that $\Gamma_\nu$ is finitely generated in \ref{prop:graded-pieces-valuation-algebra.3} and \ref{prop:graded-pieces-valuation-algebra.4} is harmless. Indeed, finite generation of  $\gr_\nu(R)$ implies that the associated monoid $\Phi_R \coloneqq \nu(R - \{0\})$ is finitely generated. For if we choose $\gamma_1, \dots, \gamma_n \in \Phi_R$ to be the degrees of a finite set of homogeneous generators of $\gr_\nu(R)$ as a $\kappa_R$-algebra, then $\gamma_1, \dots, \gamma_n$ generate $\Phi_R$ as a commutative monoid. Since $\Gamma_\nu$ is the group generated by $\Phi_R$, it follows that $\Gamma_\nu$ is a finitely generated group.
\end{remark}

\begin{proof}[Proof of \autoref{thm:finite-generation-SFR-valuation-algebra-max-rank}]
	By Abhyankar's inequality, we have $\td{\kappa_\nu/k} = 0$. Since $k$ is algebraically closed and $k \hookrightarrow \kappa_R \hookrightarrow \kappa_\nu$, it follows that $k =  \kappa_\nu = \kappa_R$. Moreover, since $\nu$ has maximal rational rank, $\nu$ is an Abhyankar valuation of $K/k$. Thus, the value group $\Gamma_\nu$ is a finitely generated Abelian group. We are then in the situation of \ref{prop:graded-pieces-valuation-algebra.4} of \autoref{prop:graded-pieces-valuation-algebra}, and the theorem follows.
\end{proof}

\bibliographystyle{skalpha}
\bibliography{main}

\end{document}